\documentclass[12pt]{amsart}
\usepackage{amsmath,amssymb}
\usepackage[all]{xy}

\newcommand\base{H/{\Lambda}{\times}M}

\newcommand\cover{H{\times}M}

\newcommand\quotient{Z{\backslash}H}

\newcommand\inv{^{-1}}
\newcommand\fiber{Z{\times}M}

\newcommand\ft{\mathfrak t}
\newcommand\fh{\mathfrak h}\newcommand\fu{\mathfrak u}
\newcommand\fl{\mathfrak l}
\newcommand\fL{\mathfrak L}

\newcommand\fz{\mathfrak z}
\newcommand\fn{\mathfrak n}
\newcommand\fm{\mathfrak m}
\newcommand\ff{\mathfrak F}

\newcommand\sw{\mathcal W}
\newcommand\sg{\mathcal G}
\newcommand\sF{\mathcal F}
\newcommand\tsw{\tilde {\mathcal W}}
\newcommand\tsv{\tilde {\mathcal V}}
\newcommand\sv{\mathcal V}\newcommand\su{\mathcal U}
\newcommand\sL{\mathcal L}

\newcommand\mP{\mathcal{P}}

\newcommand\br{\bar \rho}
\newcommand\swc{{\sw}_{[\chi]}}

\newcommand\svc{{\sv}_{[\chi]}}
\newcommand\tsvc{{\tsv}_{[\chi]}}
\newcommand\G{\Gamma}
\newcommand\g{\gamma}

\newcommand\Ga{\mathbb G}

\newcommand\Ta{\mathbb T}

\newcommand\Ra{\mathbb R}

\newcommand\Za{\mathbb Z}

\newcommand\fff{\mathfrak f}
\newcommand\Na{\mathbb N}

 \DeclareMathOperator{\Id}{Id}
 
\DeclareMathOperator{\Diff}{Diff} 
\DeclareMathOperator{\Aff}{Aff} \DeclareMathOperator{\Vect}{Vect}

\DeclareMathOperator{\Homeo}{Homeo}
\DeclareMathOperator{\Int}{Int} 
\DeclareMathOperator{\Hom}{Hom} \DeclareMathOperator{\Aut}{Aut}
\DeclareMathOperator{\Ad}{Ad} \DeclareMathOperator{\supp}{supp}
\DeclareMathOperator{\pol}{pol} \DeclareMathOperator{\Emb}{Emb}
\DeclareMathOperator{\Isom}{Isom}

\newtheorem{theorem}{Theorem}[section]
\newtheorem{proposition}[theorem]{Proposition}
\newtheorem{lemma}[theorem]{Lemma}
\newtheorem{corollary}[theorem]{Corollary}
\newtheorem{defn}[theorem]{Definition}

\begin{document}

\title[Local rigidity of group actions]{Local rigidity of affine
actions of higher rank groups and lattices}
\author[D. Fisher and G.Margulis]{David Fisher and Gregory Margulis}
\thanks{First author partially supported by NSF grant DMS-0226121 and a PSC-CUNY grant. Second author partially
supported by NSF grant DMS-0244406. The authors would also like to
thank the FIM at ETHZ for hospitality and support.}

\begin{abstract}
Let $J$ be a semisimple Lie group with all simple factors of real
rank at least two. Let $\Gamma<J$ be a lattice. We prove a very
general local rigidity result about actions of $J$ or $\Gamma$.
This shows that almost all so-called ``standard actions" are
locally rigid. As a special case, we see that any action of
$\Gamma$ by toral automorphisms is locally rigid. More generally,
given a manifold $M$ on which $\Gamma$ acts isometrically and a
torus $\Ta^n$ on which it acts by automorphisms, we show that the
diagonal action on $\Ta^n{\times}M$ is locally rigid.

This paper is the culmination of a series of papers and depends
heavily on our work in \cite{FM1,FM2}.  The reader willing to
accept the main results of those papers as ``black boxes" should
be able to read the present paper without referring to them.
\end{abstract}

\maketitle


\section{{\bf Introduction}}
\label{sec:intro}

Throughout this paper $J$ is a (connected) semisimple Lie group
with no compact factors and all simple factors of real rank at
least two, and $\Gamma<J$ is a lattice.  The purpose of this paper
is to prove the following:

\begin{theorem}
\label{theorem:main} Let $\rho$ be a quasi-affine action of $J$ or
$\Gamma$ on a compact manifold $X$.   Then the action is
$C^{\infty,\infty}$ and $C^{3,0}$ locally rigid. Furthermore,
there exists an integer $k_0$, such that the action is $C^{k,k-n}$
locally rigid for all $k>k_0$, where $n=\frac{1}{2}\dim{X}+3$.
\end{theorem}

\noindent {\bf Remark on regularity:} The number $k_0=\max(k_1,n)$
where $k_1$ is determined by properties of certain foliations
associated to the dynamics of $\rho(g_i)$, for a specific finite
set of choices of $g_1,\ldots,g_k$ in $J$ or $\G$.  If $k$ is
even, we can let $n=\frac{1}{2}\dim{X}+2$ instead.

 We proceed to define the terms in the
theorem. We say $H$ is a {\em connected real algebraic group} if
it is the connected component of the real points $H(\mathbb R)^0$
of an algebraic group defined over $\mathbb R$.

\begin{defn}
\label{definition:affine} {\bf a)} Let $H$ be a connected real
algebraic group, $\Lambda<H$ a cocompact lattice. Assume a
topological group $G$ acts continuously on $H/{\Lambda}$. We say
that the $G$ action on $H/{\Lambda}$ is {\em affine} if every
element of $G$ acts via an affine diffeomorphism.

\noindent {\bf b)} More generally, let $M$ be a compact manifold.
Assume a group  $G$ acts affinely on $H/{\Lambda}$. Choose a
Riemannian metric on $M$ and a cocycle
$\iota:G{\times}H/{\Lambda}{\rightarrow}\Isom(M)$. We call the
skew product action of $G$ on $H/{\Lambda}{\times}M$ defined by
$d{\cdot}(x,m)=(d{\cdot}x, \iota(d,x){\cdot}m)$ a {\em
quasi-affine action}.
\end{defn}

\noindent We always write $X=H/{\Lambda}{\times}M.$ Recall that an
{\em affine diffeomorphism} $d$ of $H/{\Lambda}$ is one covered by
a diffeomorphism $\tilde d$ of $H$ where $\tilde d=A{\circ}T_h$
where $A$ is an automorphism of $H$ such that $A(\Lambda)=\Lambda$
and $T_h$ is left translation by $h{\in}H$. The full group of
affine diffeomorphisms of $H/{\Lambda}$ is a finite dimensional
Lie group which we write as $\Aff(H/{\Lambda})$.  The definition
of acting affinely given above is equivalent to saying the action
is given by a homomorphism
$\pi:G{\rightarrow}{\Aff(H/{\Lambda})}$. See \cite[Section
6.1]{FM1} for a description of $\Aff(H/{\Lambda})$ and a
classification of affine actions of $J$ or $\Gamma$ as above. Note
also that the case of quasi-affine actions as defined here
includes products of affine actions with trivial actions.  Another
class of examples give the following:

\begin{corollary}
\label{corollary:torus} Let $J$ be as above and $\Gamma<J$ a
lattice.  Then any action of $\Gamma$ by automorphisms of
${\mathbb T}^m$ is $C^{\infty,\infty}$ and $C^{3,0}$ locally
rigid. Furthermore there exists a positive integer $k_0{\geq}3$,
depending on the action, such that the action is
$C^{k,k-\frac{1}{2}m-3}$ locally rigid for all
$k{\geq}\min(\frac{1}{2}m+3,k_0)$.
\end{corollary}

We now formally define local rigidity in this context.

\begin{defn}
Given a topological group $G$ and a continuous $C^{\infty}$
action, ${\rho}:D{\times}X{\rightarrow}X$, by diffeomorphisms on a
manifold $X$, we say that the action is {\em $C^{k,r}$ locally
rigid}, where $r{\leq}k$,  if any continuous action $\rho'$ by
$C^k$ diffeomorphisms, that is sufficiently $C^k$ close to $\rho$
is conjugate to $\rho$ by a small $C^r$ diffeomorphism.  We say
that a continuous action $\rho$ is {\em $C^{\infty,\infty}$
locally rigid} if any continuous $C^{\infty}$ action which is
sufficiently $C^{\infty}$ close to $\rho$ is conjugate to $\rho$
by a small $C^{\infty}$ diffeomorphism.
\end{defn}

The special case of $C^{k,k}$ local rigidity says exactly that the
homomorphism $\rho:G{\rightarrow}\Diff^k(X)$ is locally rigid.  In
other words that any homomorphism close to $\rho$ is conjugate to
$\rho$ by a small element of $\Diff^k(X)$. Since the $C^{\infty}$
topology is defined as the inverse limit of the $C^k$ topologies,
two $C^{\infty}$ diffeomorphisms are $C^{\infty}$ close if they
are $C^k$ close for some large $k$.  Our proof shows explicitly
that a $C^{\infty}$ perturbation $\rho'$ of $\rho$ which is $C^k$
close to $\rho$ is conjugate to $\rho$ by a $C^{\infty}$
diffeomorphism which is $C^{k-n}$ close to the identity where $n$
is as in Theorem \ref{theorem:main}.  The topology we take on
$\Hom(G, \Diff^k(X))$ to define close above is the compact-open
topology.

Gromov in \cite{G} and Zimmer in \cite{Z2} suggested that one
might be able to ``essentially classify" all volume preserving
smooth actions of higher rank semisimple groups and their lattices
on compact manifolds.  This would be, in a sense, a ``non-linear"
analogue of the second author's superrigidity theorems, since one
of the consequences of the superrigidity theorems is a
classification of all finite dimensional linear representations of
higher rank lattices (modulo issues concerning finite image
representations). In \cite{Z2}, Zimmer also proposed the study of
local rigidity of known actions of higher rank lattices on compact
manifolds, as a ``non-linear" analogue of the classical local
rigidity theorems of Calabi-Vesentini, Selberg and Weil. These
show that any cocompact lattice $\Gamma$ in any simple Lie group
$J$ is locally rigid, as long as $J$ is not locally isomorphic to
$SL_2(\mathbb R)$ \cite{CV,S,W}. I.e. any embedding of $\Gamma$ in
$J$ close to the defining one $i:{\Gamma}{\rightarrow}J$ is simply
a conjugate of $i$ by a small element of $J$.  Since $J$ acts
transitively on $J/{\Gamma}$, our theorem can be taken to be a
generalization of Weil's result in the case when $J$ is a higher
rank simple group. A perturbation $\Gamma'$ of $\Gamma$ in $J$
defines a perturbation of the original $J$ action on $J/{\Gamma}$
since $J/{\Gamma}$ and $J/{\Gamma'}$ are diffeomorphic.  The
conjugacy between these actions can easily be seen to give a
conjugacy between $\Gamma$ and $\Gamma'$.


 Many results have been proven concerning local rigidity of affine
 actions of higher rank lattices and Lie groups, particularly when
  the action is assumed to satisfy some strong hyperbolicity condition.
    The first results of this kind are due to Hurder \cite{H1}. He proved
that the standard action on ${\mathbb T}^n$ of any finite index
subgroup in $SL_n(\mathbb Z)$ is deformation rigid for $n{\geq}3$.
(This involves assuming a path of nearby actions and obtaining a
path of conjugacies.) The same actions were shown to be locally
rigid in \cite{KL1} and \cite{KLZ}.  Many other results along
these lines were obtained by many authors, we refer to the
introduction of \cite{MQ} for a more detailed discussion. Here we
mention that all standard Anosov actions on tori and nilmanifolds
were proven to be locally rigid in \cite{KS} and all so-called
weakly hyperbolic actions were proven to be locally rigid in
\cite{MQ}.

For isometric actions, there are also results.  In \cite{B1},
Benveniste shows that any isometric action of any cocompact
lattice in a group $J$ as above is $C^{\infty, \infty}$ locally
rigid. The interested reader should refer to the introduction to
\cite{FM2} for a discussion of earlier, weaker results by Zimmer
concerning (certain) perturbations of isometric actions of groups
with property $(T)$.  In our previous paper \cite{FM2}, we have
proven:

\begin{theorem}
\label{theorem:isomrigid} Let $G$ be locally compact, compactly
generated group with property $(T)$. Let $(X,g)$ be a compact
Riemannian manifold, and let $\rho$ be an action of $G$ on $X$ by
isometries. Then the action  $C^{\infty,\infty}$ locally rigid and
is $C^{k,k-\kappa}$ locally rigid for any $k{\geq}2$ and any
$\kappa>0$.
\end{theorem}

\noindent We remark that Theorem \ref{theorem:isomrigid} holds for
a much broader class of groups than Theorem \ref{theorem:main}.

 The proof of
Theorem \ref{theorem:main} uses a foliated generalization of
Theorem \ref{theorem:isomrigid} also proven in \cite{FM2}.  This
result is recalled below in section \ref{section:conjugacy} where
it is applied in the course of our proof.

For actions which are neither weakly hyperbolic nor isometric all
previous results, due to Nitica and Torok, concern affine actions
which are products of Anosov actions and trivial actions
\cite{NT,NT2,T}. For example, take the standard action of
$SL_n(\mathbb Z)$ on ${\mathbb T}^n$ and let $\rho$ denote the
action obtained on ${\mathbb T}^n{\times}S^1$ by taking
${\gamma}(t,s)=(\gamma{t},s)$.  Then Nitica and Torok show that,
given $k>0$, any $C^{\infty}$ action $\rho'$ that is sufficiently
$C^2$ close to $\rho$ is conjugate to $\rho$ by a $C^0$ small,
$C^k$ diffeomorphism. (This result does not imply $C^{\infty}$
local rigidity because the size of the perturbation must be made
smaller to obtain more derivatives in the conjugacy.) Their full
result is more general, allowing one to replace the standard
action of $SL_n(\mathbb Z)$ on ${\mathbb T}^n$ by any so-called
$TNS$ action of a higher rank lattice on a torus. They also prove
some more general results for deformation rigidity, but always for
products of $TNS$ and trivial actions.

We note here that non-locally rigid volume preserving actions of
higher rank semisimple groups and their lattices on compact
manifolds have been constructed, first in \cite{KL2} and later and
more generally in \cite{B2}.  Those in \cite{B2} are even shown to
have smooth volume preserving deformations. See also \cite{F} for
a more general construction and another proof that the
deformations are non-trivial. A weaker result is shown in
\cite{KL2}.

\section{\bf Affine actions, perturbations and quotients}
\label{section:resultsfromfm1}

This section primarily recalls results from \cite{FM1}.  Let
$\tilde J$ be the universal cover of $J$ and $\tilde \G$ the
pre-image of $\G$ under the covering map $\tilde J{\rightarrow}J$.
Any action of $J$ or $\G$ can be viewed as an action of $\tilde J$
or $\tilde \G$ respectively, so we assume, without loss of
generality, that $J$ is simply connected.

\subsection{Describing affine actions}
\label{subsection:describingaffineactions}

In this section we recall from \cite{FM1} another description of
the actions we are considering.  This description provides an
extremely simple description of the derivative cocycle for the
action which allows a simple description of the dynamical
foliations for elements of the acting group, as well as of the
central foliation of the entire group.

Throughout this section $H$ will be a connected real algebraic
group and $\Lambda<H$ will be a cocompact lattice.  We now recall
three technical results from \cite{FM1}.

\begin{theorem}
\label{theorem:describingactionsG} Let $\rho$ be an affine action
of $J$ on $H/{\Lambda}$. Then the action $\rho$ is given by
$\rho(j)[h]=[\pi_0(j)h]$ where $\pi_0:J{\rightarrow}H$ is a
continuous homomorphism.
\end{theorem}

\noindent This is a special case of \cite[Theorem 6.4]{FM1}.  As
indicated there, the result holds with the weaker assumption that
$J$ has no compact simple factors.

The analogous result for $\G$ actions is more complicated and can
require that we view $H/{\Lambda}$ as a homogeneous space for a
different Lie group.  The following is a rearrangement of
\cite[Proposition 6.3]{FM1}.  Given a Lie group $L$, we denote
it's automorphism group by $\Aut(L)$.  Since $\Aut(L)$ is a closed
subgroup of $GL(\dim(L))$ it is a Lie group.

\begin{proposition}
\label{proposition:Hreplacement} Given a real algebraic group $H$
there is a connected cover $p:H'{\rightarrow}H$ and a realization
of $H'$ as a connected real algebraic group, such that
\begin{enumerate}

\item the connected component $\Aut(H')^0$ of $\Aut(H')$ has the
structure of a connected real algebraic group,

\item $\Aut(H')^0<\Aut(H')$ is a finite index subgroup,

\item $Aut(H')^0$ acts rationally on $H'$.

\end{enumerate}
\end{proposition}

\noindent The key point is to choose the algebraic structure on
$H'$ so that the connected component of the center of $H'$ is
contained in the unipotent radical.  It follow that
$\Aut(H')^0{\ltimes}H'$ is a connected real algebraic group.

Let $\Lambda'$ be $p{\inv}(\Lambda)$.  It follows from the
construction given in \cite{FM1} that, possibly after passing to a
finite index subgroup $\G'$ in $\G$, any affine action $\rho$ of
$\G$ on $H/{\Lambda}$ remains affine when we view $H/{\Lambda}$ as
$H'/{\Lambda'}$.  For the remainder of this paper, we assume that
we have replaced our group $H$ with a group $H'$ as described in
Proposition \ref{proposition:Hreplacement}.

Before giving the analogous description of affine $\Gamma$
actions, we need to recall a consequence of the superrigidity
theorems \cite{M1,M2}. In \cite{M1,M2} these are only stated for
$J$ algebraic, but the extension to $J$ as assumed here is
sketched in \cite[Theorem 7.12]{F2}, see \cite{FM1} for detailed
proofs. We will use the notation introduced here in the statements
below. If $J$ is as above and $\Gamma<J$ is a lattice, and $L$ is
an algebraic group, we call a homomorphism
$\pi:\Gamma{\rightarrow}L$ {\em superrigid} if it almost extends
to a homomorphism of $J$. This means that there is a continuous
homomorphism $\pi^E:J{\rightarrow}L$ and a homomorphism
$\pi^K:\Gamma{\rightarrow}L$ with bounded image such that
$\pi(\gamma)=\pi^E(\gamma)\pi^K(\gamma)$ and $\pi^E(\Gamma)$
commutes with $\pi^K(\Gamma)$. The superrigidity theorems imply
that any continuous homomorphism of $\Gamma$ into an algebraic
group is superrigid. This can be deduced easily from Lemma VII.5.1
and Theorems VII.5.15 and VII.6.16 of \cite{M2}.

\begin{theorem}
\label{theorem:describingactionsGamma} Let $\rho$ be an affine
action of $\Gamma$ on $H/{\Lambda}$. Then there is a finite index
subgroup $\Gamma'<\Gamma$ and a homomorphism
$\pi_0:\G'{\rightarrow}\Aff(H)=\Aut(H){\ltimes}H$ such that
$\rho(\gamma)[h]=[\pi_0(\gamma)h]$.  Furthermore, we can assume
that $\pi_0(\Gamma')$ is contained in $\Aut(H)^0{\ltimes}H$ and
that $\pi_0(\g)=\pi_0^E(\g)\pi_0^K(\g)$ where
$\pi^E_0:J{\rightarrow}\Aff(H)$ is a homomorphism and
$\pi^K_0:\G'{\rightarrow}\Aff(H)$ is a homomorphism with bounded
image, and the images of $\pi^E_0$ and $\pi^K_0$ commute.
\end{theorem}

\noindent This is a rephrasing of \cite[Theorem 6.5]{FM1}.  The
final conclusion concerning the fact that $\pi_0$ is the product
of a the restriction of a homomorphism of $J$ and a homomorphism
with bounded image follows from the superrigidity theorems
discussed above.

We can now describe the {\em central foliation} for a quasi-affine
action $\rho$ of either $J$ or $\G$. We will denote the central
foliation $\ff$.  If $\rho$ is a $J$ action, and $M$ is trivial,
then the central foliation is just the orbit foliation for the
left action of $Z=Z_H(\pi_0(J))$ on $H/{\Lambda}$.  If $M$ is
non-trivial, we have a projection $\base{\rightarrow}H/{\Lambda}$
and the central foliation is given by the pre-images in $\base$ of
$Z$ orbits in $H/{\Lambda}$.

Let $\Gamma'$ be the subgroup of finite index given by Theorem
\ref{theorem:describingactionsGamma} and further assume that
$\Gamma'$ is normal. Let $A$ be the connected component of
$\Aut(H)$, and let $L=A{\ltimes}H$. Note that $L$ is an algebraic
group. In this case, we let $Z=Z_L(\pi_0^E(\G')){\cap}H$. If $M$
is trivial, the {\em central foliation} for the action, which we
denote $\ff$, is then defined to be the foliation given by orbits
of $Z$ on $H/{\Lambda}$.  If $M$ is non-trivial, we have a
projection $\base{\rightarrow}H/{\Lambda}$ and the central
foliation is given by the pre-images  in $\base$ of $Z$ orbits in
$H/{\Lambda}$.

We will refer to the tangent space of the central foliation $\ff$
as {\em the central distribution for the group action}.

We now want to define a Riemannian metric on $\base$ so that the
$\rho$ action is isometric along leaves of $\ff$.  Since $M$ is
assumed in Definition \ref{definition:affine} to be Riemannian
manifold with Riemannian metric $g_M$, and $\rho$ is defined to be
isometric along $M$ fibers, it suffices to define a Riemannian
metric on $H/{\Lambda}$ for which the affine action on
$H/{\Lambda}$ is isometric along $Z$ orbits. Let $\fh$ be the Lie
algebra of $H$. An inner product on $\fh$ defines a right
invariant Riemannian metric on $H$ and therefore a Riemannian
metric on $H/{\Lambda}$.  For the case of $J$ actions, we have
that the derivative action on $H/{\Lambda}$ is given by
$D\rho(j)([h],v)=([\pi_0(j)h, \Ad_H(\pi_0(j))v)$ for $j{\in}J$.
Since $\Ad_H{\circ}\pi_0|_{\fz}$ is the trivial representation of
$J$, it is clear that any inner product on $\fh$ defines a
Riemannian metric with the desired property.  Let $\fm$ be an
$\Ad(\pi_0(J))|_{\fh}$ invariant complement to $\fz$.  For
simplicity in arguments below, we choose a metric on $\fh$ such
that $\fm$ is orthogonal to $\fz$.

For $\G$ actions, we need to be slightly more careful.  Let $\fl$
be the Lie algebra of $L$, which contains $\fh$ as an ideal.  We
denote by $\Ad_{\fh}$ the restriction of the adjoint action of $L$
on $\fl$ to $\fh$. Recall that $\pi^K(\G')<C$ where $C<L$ is
compact and take an $\Ad_{\fh}(C)$ invariant metric on $\fh$. This
then defines a Riemannian metric on $\base$ for which $\rho(\G')$
is isometric along the central foliation defined above.

We perform one further modification to the metric to guarantee
that the action of all of $\G$, and not just $\G'$, is isometric
along the central foliation. Since the image of $\pi_0$ of $\G'$
is semisimple, we can choose an $\Ad_{\fh}(\pi_0(\G'))$ invariant
subspace $\fm<\fh$ orthogonal to $\fz$ such that
$\fh=\fz{\oplus}\fm$, and we can choose an inner product on $\fh$
for which $\fz$ is orthogonal to $\fm$. There are corresponding
subbundles of $T(H/{\Lambda})$ which we can write as
$H/{\Lambda}{\times}\fm$ and $H/{\Lambda}{\times}\fz$. Note that
because these bundles are $\Gamma'$ invariant and $\Gamma'$ can be
chosen to be normal in $\Gamma$, they are also $\Gamma$ invariant.
Given a vector space $V$, denote by $S_+^2(V)$ the cone of
positive definite symmetric two tensors on $V$. The Riemannian
metric on $H/{\Lambda}$ is a section $g_{\fh}$ of
$H/{\Lambda}{\times}S_+^2(\fh)$ which lies in the subbundle given
by $H/{\Lambda}{\times}S_+^2(\fz){\oplus}S_+^2(\fm)$ or
equivalently as the sum in $H/{\Lambda}{\times}S_+^2(\fh)$ of a
section of $g_{\fz}{\in}H/{\Lambda}{\times}S_+^2(\fz)$ and a
section of $g_{\fm}{\in}H/{\Lambda}{\times}S_+^2(\fm)$.  Since
$g_{\fz}$ is $\G'$ invariant, and $\G'<\G$ is finite index, we can
average $g_{\fz}$ over coset representatives for $\G/\G'$ to
obtain a $\G$ invariant section $g'_{\fz}$ in
$H/{\Lambda}{\times}S_+^2(\fz)$.  Replacing $g_{\fh}$ by
$g'_{\fz}{\oplus}g_{\fm}$ we have a Riemannian metric on
$H/{\Lambda}$ such that the entire $\G$ action is isometric along
the central foliation.

\subsection{\bf Semiconjugacy}
\label{subsection:semiconjugacy}

Let $H,\Lambda,\G$ and $J$ be as in the preceding subsection, let
$G=J$ or $\G$ and let $\rho$ be a quasi-affine action of $G$ on
$H/{\Lambda}{\times}M$.   Then by Theorems
\ref{theorem:describingactionsG} and
\ref{theorem:describingactionsGamma} there is a subgroup $G'<G$
such that $\rho|_{G'}$ is defined by a continuous homomorphism
$\pi:G'{\rightarrow}\Aut(H){\ltimes}H$. It follows that the $G'$
action lifts to $H{\times}M$. As explained in \cite{FM1},
following the statement of Theorem $6.7$, for any small enough
$C^0$ perturbation $\rho'$ of $\rho$, the $G'$ action defined by
$\rho'$ also lifts to $H{\times}M$. (We note that this is
trivially true for actions of connected groups.) By the discussion
in subsection \ref{subsection:describingaffineactions} there is a
unique subgroup $Z$ in $H$ which is the maximal subgroup of $H$
such that the derivative of $\rho$ on $Z$ cosets is an isometry
for an appropriate choice of metric on $H/{\Lambda}$. The
description given there shows that the lift of $\rho(G')$ to
$H{\times}M$ descends to an action $\bar \rho$ of $G'$ on
$Z{\backslash}H$.  We denote by $p$ the natural projection
$H{\times}M{\rightarrow}Z{\backslash}H$.

\begin{theorem}
\label{theorem:semiconjugacy}   Let
$H/{\Lambda}{\times}M,\rho,G',Z$ and $\br$ be as in the preceding
paragraph. Given any action $\rho'$ sufficiently $C^1$ close to
$\rho$, there is a continuous $G'{\times}{\Lambda}$ equivariant
map $f:(H{\times}M, \rho'){\rightarrow}(Z{\backslash}H, \br)$, and
$f$ is $C^0$ close to $p$.  Furthermore if
$\rho'{\rightarrow}\rho$ in the $C^1$ topology then
$f{\rightarrow}p$ in the $C^0$ topology.
\end{theorem}

\noindent This is \cite[Theorem $1.8$]{FM1}. We note that
$f{\rightarrow}p$ in the $C^0$ topology means that
$d(f(x),p(x)){\rightarrow}0$ uniformly on $H{\times}M$. There is
some ambiguity in this, since there is no $\Lambda$ invariant
metric on $Z{\backslash}H$, but it is true that
$d(f(x),p(x)){\rightarrow}0$ in the metric on $Z{\backslash}H$
which makes $p$ a Riemannian submersion.  For the remainder of
this subsection, we assume that $G'=G$.

The map $f$ defines a partition $\tilde \ff'$ of $H{\times}M$ into
sets of the form $f{\inv}(x)$ where $x$ is in $Z{\backslash}H$.
Since $f$ is $\Lambda$ equivariant, this partition defines to a
partition $\ff'$ of $H/{\Lambda}{\times}M$.  We will show that, as
a consequence of Theorem \ref{theorem:semiconjugacy} there is a
$\Lambda$ equivariant map
$\theta:H{\times}M{\rightarrow}H{\times}M$ mapping $\ff'$ to $\ff$
and intertwining the actions of $G$ on $\ff$ and $\ff'$, but first
we need some definitions.

 If $\mP$ is a partition of a topological space $X$
and $\rho$ is an action of a group $D$ on $X$, then we say $\rho$
preserves $\mP$, if for any set $V{\in}\mP$, the set $\rho(d)V$ is
in $\mP$ for any $d$ in $D$.

Given two actions $\rho$ and $\rho'$ of a group $D$ on a
topological space $X$ and two partitions $\mP$ and $\mP'$ of the
space $X$ where $\rho$ preserves $\mP$ and $\rho'$ preserves
$\mP'$, we call a map $\phi:X{\rightarrow}X$ a {\em partition
semi-conjugacy from $(X,\rho,\mP)$ to $(X,\rho'\mP')$} if for any
subset $V{\in}\mP$ we have
\begin{enumerate}
\item $\phi(V)$ is an element of the partition $\mP'$ and,

\item $\phi(\rho(d)V)=\rho'(d)\phi(V)$ for any $d$ in $D$.
\end{enumerate}

\noindent If $\phi$ is a homeomorphism, we call $\phi$ a {\em
partition conjugacy between $(X,\rho,\mP)$ and $(X,\rho'\mP')$}.
Similarly one can refer to actions as being {\em partition
(semi-)conjugate}.

To be consistent with the vocabulary of \cite{HPS}, when we are
given two actions $\rho$ and $\rho'$ of a group $D$ on a
topological space $X$ where $\rho$ (resp. $\rho'$) preserves a
foliation $\ff$ of $X$ (resp. a foliation $\ff'$ of $X$), a
partition (semi-)conjugacy from $(X,\rho,\ff)$ to $(X,\rho',\ff')$
will be called a {\em leaf (semi-)conjugacy}. Similarly, when we
do not want to make explicit reference to the (semi-)conjugacy, we
will say that two actions are  {\em leaf (semi-)conjugate}.

We now construct a map $\tilde \theta:\cover{\rightarrow}\cover$
using $f$ and $p$. The space $\cover$ is a smooth locally trivial
fiber bundle over $\quotient$ with fiber $\fiber$, so given
$x{\in}Z{\backslash}H$, we can find a neighborhood $U$ of $x$ such
that $p{\inv}(U)$ is diffeomorphic to $U{\times}\fiber$.  We can
therefore introduce coordinates on $p{\inv}(U)$ of the form
$(x,y)$ where $x$ is in $U$ and $y{\in}\fiber$.  In these
coordinates, $p(x,y)=x$.  We can further write $y=(z,m)$ where
$z{\in}Z$ and $m{\in}M$.

Moreover, if we let $\fm$ be the $G$ invariant complement of $\fz$
in $\fh$, then the tangent space to $\cover$ at any point
$(x,z,m)$ can be written as
$T(\cover)_{(x,z,m)}=\fm\oplus\fz\oplus{TM_m}$. We can further
choose the local product structure on $p{\inv}(U)$ such that
$(U,y)=\exp_yW$ where $W$ is the product of a fixed small ball in
$\fm$ with a small ball in $TM_m$. By shrinking $W$ and therefore
$U$ slightly, we obtain a trivialization of $p{\inv}(U)$ that
extends to a trivialization of $p{\inv}(U')$ for $U'$ an open set
strictly containing $U$.

By choosing $\rho'$ close enough to $\rho$, we can arrange for $f$
to be arbitrarily $C^0$ close to $p$, uniformly on $\cover$. This
implies that given any compact set $K$ in $H{\times}M$, by
restricting to sufficiently small $C^1$ perturbations $\rho'$ of
$\rho$, we can make the Hausdorff distance between
$f{\inv}(x){\cap}K$ and $p{\inv}(x){\cap}K$ as small as desired
for every $x$ in $f(K){\cap}p(K)$. Since $f$ and $p$ are $\Lambda$
equivariant and the $\Lambda$ action on $\cover$ is cocompact, for
 small enough perturbations we have
$f{\inv}(x){\subset}p{\inv}(U)$. Then for a point $(x,y)$, we let
$\tilde \theta(x,y)=(U,y){\cap}p{\inv}(f(x,y))$. Therefore $\tilde
\theta(f{\inv}(x))=p{\inv}(x)$ for any $x{\in}Z{\backslash}H$.
Since $\tilde \theta$ is $\Lambda$ equivariant by construction, we
have a map $\theta:\base{\rightarrow}\base$ and have established
the following:

\begin{corollary}
\label{corollary:semiconjugacy} Let $\rho'$ be an action of $G$ on
$H/{\Lambda}{\times}M$ sufficiently $C^1$ close to $\rho$.  Then
there is a $C^0$ small map
$\theta:H/{\Lambda}{\times}M{\rightarrow}H/{\Lambda}{\times}M$
with the following properties.
\begin{enumerate}

\item If ${\rho'}{\rightarrow}\rho$ in the $C^1$ topology then
$\theta{\rightarrow}\Id$ in the $C^0$ topology.

\item $\theta(\ff')=\ff$

\item  the map $\theta$ is a partition semi-conjugacy from
$(\base,\rho,\ff')$ to $(\base,\rho,\ff)$.
\end{enumerate}
\end{corollary}

\noindent{\bf Remarks:}\begin{enumerate} \item One can deduce
Corollary \ref{corollary:semiconjugacy} directly from the proof of
Theorem \ref{theorem:semiconjugacy}. \item The argument there
proves more.  It shows that the set of maps
$\{\theta{\circ}\rho'(g)|g{\in}G\}$ defines a $G$ action on
$\base$ that is $C^0$ close to $\rho$. \item Remark $(2)$ can be
deduced easily from \cite[Theorem 6.7]{FM1}, but to avoid
introducing additional notations and defintions, we do not do this
here. \item The conclusion of Remark $(2)$ will follow once we
show, in subsection \ref{section:fibers}, that $\theta$ is a $C^0$
small homeomorphism.
\end{enumerate}

\section{\bf Hyperbolic dynamics and stability modulo central directions}
\label{section:structuralstability}

In this section, we show that the map
$\theta:\base{\rightarrow}\base$ defined in the last subsection is
a homeomorphism.  Since $\theta(\ff')=\ff$, this implies that
$\ff'$ is a foliation.  We further show that $\ff'$ is a foliation
by $C^r$ leaves where $r$ depends only on the $C^1$ size of the
perturbation.  The map $\theta$ is then easily seen to be $C^r$
along leaves of $\ff$.  For technical reasons involving the last
steps of our proof, once we have shown that $\theta$ is a
homeomorphism, we prefer to work with $\phi=\theta{\inv}$.

We now briefly outline the argument of the section.  Some of the
terminology used here is defined below in subsection
\ref{section:prelim3}.  First in subsection \ref{section:prelim3}
we recall some definitions concerning partially hyperbolic
diffeomorphisms and a theorem of Hirsch, Pugh and Shub.  In
subsection \ref{subsection:dynamicsofaffine}, we prove some basic
facts concerning the dynamics of the affine actions of $G$ that we
are considering and produce a finite subset $\Phi$ of $G$ such
that the intersections of the central foliations of $\rho(g)$ for
$g$ in $\Phi$ is the central foliation for $G$ defined above in
subsection \ref{subsection:describingaffineactions}.  In
subsection \ref{section:fibers}, we show that the map
$\theta:\base{\rightarrow}\base$ defined in subsection
\ref{subsection:semiconjugacy} is a homeomorphism and we let
$\phi=\theta{\inv}$. It then follows that $\ff'=\phi(\ff)$ is
$\rho'$ invariant foliation. Finally in subsection
\ref{subsection:HPS}, we show that any leaf $\fL$ of $\ff'$ is the
transverse intersection of central leaves $\sw^c_{\rho'(g)}$ for
the diffeomorphisms $\rho'(g)$ where $g$ is in $\Phi$.  Since the
theorem of Hirsch, Pugh and Shub implies that each foliation
$\sw^c_{\rho'(g)}$ is by $C^r$ leaves, it follows that $\ff'$ is a
foliation by $C^r$ leaves, where $r$ depends only on the $C^1$
size of the perturbation $\rho'$.

\subsection{\bf Hyperbolic dynamics and foliations}
\label{section:prelim3}

The use of the word {\em foliation} varies with context. Here a
{\em foliation by $C^k$ leaves} will be a continuous foliation
whose leaves are $C^k$ injectively immersed submanifolds that vary
continuously in the $C^k$ topology in the transverse direction. To
specify transverse regularity we will say that a foliation is
transversely $C^r$.  A foliation by $C^k$ leaves which is
tranversely $C^k$ is called simply a $C^k$ foliation. Note our
language does not agree with that in the reference \cite{HPS}
where our foliation by $C^k$ leaves is a $C^k$ unbranched
lamination and sometimes a $C^k$ injective leaf immersion.  Given
a foliation $\ff$, we denote the leaf through a point $x$ by
$\ff(x)$.

Given a foliation by $C^1$ leaves, $\ff$, of a manifold $X$, a
diffeomorphism $f$ is said to be {\em $r$-normally hyperbolic} to
the foliation $\ff$ if there exists a continuous $f$ invariant
splitting $TX=E_{f}^u{\oplus}T\ff{\oplus}E_{f}^s$ such that for
every $x{\in}X$,
\begin{enumerate}
\item $\|Df_x|_{E_f^u}{\inv}\|{\inv}>\|Df_x|_{\ff}\|^r$ and,

\item $\|Df_x|_{E_f^s}\|<\|Df_x|_{\ff}{\inv}\|^{-r}$.
\end{enumerate}

\noindent For any invariant subbundle $V$ of $TX$ and any fixed
Riemannian metric on $X$, the norm above is the operator norm of
$Df_x|_{V_x}$. See \cite{HPS}, chapter 1, for a more detailed
discussion of normal hyperbolicity. There $r$-normally hyperbolic
is also called ``immediately, relatively $r$-normally hyperbolic."
Also the definition given there is slightly different, and applies
also to non-invertible $f$.  That the definitions are equivalent
for $f$ invertible is the content of the remark following
Definition $1$ in the introduction to \cite{HPS}.

We note that $f$ being $r$-normally hyperbolic to $\ff$ is
equivalent to the existence of constants $a,b>1$ with $a>b^r$ and
a continuous $f$ invariant splitting
$TX=E_{f}^u{\oplus}T\ff{\oplus}E_{f}^s$ such that

\begin{enumerate}

\item $\|Df^n(v^u)\|{\geq}a^n\|v^u\|$ for all $v^u{\in}E_{f}^u$,

\item $\|Df^n(v^s)\|{\leq}a^{-n}\|v^s\|$ for all $v^s{\in}E_{f}^s$
and

\item  ${\inv}b^{-n}\|v^0\|<\|Df^n(v^0)\|{\leq}{b^n}\|v^0\|$ for
all $v^0{\in}T{\ff}$ and all integers $n$.

\end{enumerate}

The definition of $r$-normally hyperbolic is motivated by the
theory of partially hyperbolic diffeomorphisms. Given an
automorphism $f$ of a vector bundle $E{\rightarrow}X$ and
constants $a>b{\geq}1$, we say $f$ is {\em $(a,b)$-partially
hyperbolic} or simply {\em partially hyperbolic} if there is a
metric on $E$ and a constant and $C{\geq}1$ a continuous $f$
invariant
 non trivial splitting $E=E_{f}^u{\oplus}E_{f}^c{\oplus}E_{f}^s$
 such that:
\begin{enumerate}

\item $\|f^n(v^u)\|{\geq}Ca^n\|v^u\|$ for all $v^u{\in}E_{f}^u$,

\item $\|f^n(v^s)\|{\leq}C{\inv}a^{-n}\|v^s\|$ for all
$v^s{\in}E_{f}^s$ and

\item  $C{\inv}b^{-n}\|v^0\|<\|f^n(v^0)\|{\leq}C{b^n}\|v^0\|$ for
all $v^0{\in}E_{f}^c$ and all integers $n$.

\end{enumerate}

\noindent A $C^1$ diffeomorphism $f$ of a manifold $X$ is {\em
$(a,b)$-partially hyperbolic} if the derivative action $Df$ is
$(a,b)$-partially hyperbolic on $TX$. We remark that for any
partially hyperbolic diffeomorphism, there always exists an {\it
adapted metric} for which $C=1$. Note that $E_{f}^c$ is called the
{\em central distribution} of $f$, $E_{f}^u$ is called the {\em
unstable distribution} of $f$ and $E_{f}^s$ the {\em stable
distribution} of $f$. We will also refer to the sums
$E_{f}^{cu}=E_{f}^u{\oplus}E_{f}^c$ and
$E_{f}^{cs}=E_{f}^s{\oplus}E_{f}^c$ as the weak unstable and weak
stable distributions, respectively.

Integrability of various distributions for partially hyperbolic
dynamical systems is the subject of much research.  The stable and
unstable distributions are always tangent to invariant foliations
which we call the stable and unstable foliations and denote by
$\sw_f^s$ and $\sw_f^u$.  If the central distribution is tangent
to an $f$ invariant foliation, we call that foliation a {\em
central foliation} and denote it by $\sw^c_f$. If there is a
unique foliation tangent to the central distribution we call the
central distribution {\em uniquely integrable}. For smooth
distributions unique integrability is a consequence of
integrability, but the central distribution is usually not smooth.
For general partially hyperbolic diffeomorphisms, unique
integrability of central foliations is difficult to establish.  If
the central distribution of an $(a,b)$-partially hyperbolic
diffeomorphism $f$ is tangent to an invariant foliation $\sw^c_f$,
then $f$ is $r$-normally hyperbolic to $\sw^c_f$ for any $r$ such
that $a>b^r$.

Given $C^k$ partially hyperbolic diffeomorphism $f$ which is
$l$-normally hyperbolic to a central foliation, for $k,l{\geq}1$,
it follows from \cite[Theorem 6.8]{HPS} that there are foliations
tangent to the weak (un)stable distribution, which we call the
{\em weak (un)stable foliation} and denote by $\sw_f^{cu}$ and
$\sw_f^{cs}$.


In subsection \ref{subsection:HPS}, we need to use the work of
Hirsch-Pugh-Shub on perturbations of partially hyperbolic actions
of $\mathbb Z$. We state a special case of some of their results
from \cite{HPS}.

\begin{theorem}
\label{theorem:hps} Let $f$ be an  $(a,b)$-partially hyperbolic
$C^k$ diffeomorphism of a compact manifold $M$ which is
$k$-normally hyperbolic to a $C^k$ central foliation $\sw^c_f$.
Then for any $\delta>0$, if $f'$ is a $C^k$ diffeomorphism of $M$
which is sufficiently $C^1$ close to $f$ we have the following:
\begin{enumerate}
\item $f'$ is $(a',b')$-partially hyperbolic, where
$|a-a'|<\delta$ and $|b-b'|<\delta$,  and  the splitting
$TM=E_{f'}^u{\oplus}E_{f'}^c{\oplus}E_{f'}^s$ for $f'$ is $C^0$
close to the splitting for $f$;

\item there exist $f'$ invariant foliations by $C^k$ leaves
$\sw^{cs}_{f'}$ tangent to $E_{f'}^c{\oplus}E_{f'}^s,
\sw^{cu}_{f'}$ tangent to $E_{f'}^c{\oplus}E_{f'}^u,\sw^{s}_{f'}$
tangent to $E_{f'}^s,\sw^{u}_{f'}$ tangent to $E_{f'}^u$ and
$\sw^c_{f'}$ tangent to $E^c_{f'}$, and each foliation is close in
the natural topology on foliations by $C^k$ leaves to the
corresponding foliation for $f$.


\end{enumerate}
\end{theorem}

\noindent Statement $(1)$ is standard. Statement $(2)$ follows
from \cite{HPS} Theorem $6.1$ statement $f$, see also Theorem
$6.8$ of that book for more details.  
(The exact results in \cite{HPS} are more general.)

\subsection{Dynamics of affine actions}
\label{subsection:dynamicsofaffine}

For the remainder of this section the group $G$ is either a
connected, simply connected, semisimple Lie group, $J$, with all
simple factors of real rank at least two or a lattice $\G<J$. We
fix a manifold $M$, a real algebraic group $H$, and a cocompact
lattice $\Lambda<H$ and fix a quasi-affine action $\rho$ of $G$ on
$\base$.  We recall from subsection
\ref{subsection:describingaffineactions} that there is a finite
index subgroup $G'<G$ such that $\rho$ is defined by a
homomorphism $\pi:G'{\rightarrow}<\Aut(H){\ltimes}H$ where
$\pi=\pi^E\pi^K$, where $\pi^E$ is a continuous homomorphism of
$J$, the homomorphism $\pi^K$ has bounded image, and the images of
$\pi^E$ and $\pi^K$ commute. For the remainder of this section, we
assume
 that $G'=G$.  As above we let $L$ be the product of the Zariski
 closure of $\pi(G')$ and $H$ and let
 $Z=Z_L(\pi^E(J))\cap{H}$ and $\fz$ be the Lie algebra of $Z$.
(Note that when $G=J$, $L$ is always just $H$.) According to the
discussion in section \ref{subsection:describingaffineactions}, we
can fix a Riemannian metric on $H/{\Lambda}{\times}M$ such that
$\rho(g)$ is an isometry of the metric restricted to the tangent
space of $\ff$. Given $g{\in}G$  there is a natural choice of
$\rho(g)$ invariant sub-bundles of $T(H/{\Lambda})$ with respect
to which $\rho(g)$ is partially hyperbolic whenever
$\Ad(\pi(g))|_{\fh}$ has an eigenvalue off the unit circle. We
first describe the case of affine actions. Writing coordinates on
$T(H/{\Lambda}){\cong}H/{\Lambda}{\times}\fh$ as $([h],v)$ with
$[h]$ in $H/{\Lambda}$ and $v{\in}\fh$ we have
$D\rho(g)=(\rho(g)([h]), \Ad_{\fh}(\pi(g))v)$. We let
$\fff^s_{\rho(g)}$ (resp. $\fff^u_{\rho(g)}$) be the subspace of
$\fh$ for which $\Ad(\pi(g))$ has all eigenvalues of modulus less
than one (resp. all eigenvalues of modulus greater than one) and
$\fff^c_{\rho(g)}$ be the subspace of $\fh$ where
$\Ad_{\fff}(\pi(g))$ has all eigenvalues of modulus one. We can
then define sub-bundles of $T(H/{\Lambda})$ as
$E^s_{\rho(g)}=H/{\Lambda}{\times}\fff^s_{(\rho(g))},
E^u_{\rho(g)}=H/{\Lambda}{\times}\fff^u_{(\rho(g))}$ and
$E^c_{\rho(g)}=H/{\Lambda}{\times}\fff^c_{(\rho(g))}$.  It is
straightforward to verify that $\rho(g)$ is partially hyperbolic
with respect to this splitting whenever this splitting is
non-trivial.  For the remainder of this paper, whenever we refer
to $\rho(g)$ as a partially hyperbolic diffeomorphism, we mean
partially hyperbolic with respect to this choice of splitting. We
collect here some basic consequences for the dynamics of the
action $\rho$.

\begin{proposition}
\label{proposition:explicitfoliations} For any affine action
$\rho$ of $G$ on $H/{\Lambda}$ and any $g{\in}G$ there are Lie
subgroups $F^s_{\rho(g)},F^u_{\rho(g)}$ and $F^c_{\rho(g)}$ in $H$
such that the foliations $\sw^s_{\rho(g)},\sw^u_{\rho(g)}$ and
$\sw^c_{\rho(g)}$ consist of orbits of the corresponding group
acting on the left on $H/{\Lambda}$.  Furthermore
\begin{enumerate}

\item the groups $F^s_{\rho(g)}$ and $F^u_{\rho(g)}$ are
nilpotent,

\item $Z<F^c_{\rho(g)}$ and
$Z{\cap}F^s_{\rho(g)}=Z{\cap}F^u_{\rho(g)}=1$,

\item for every point in $H/{\Lambda}$ the orbit map for
$F^s_{\rho(g)}$ and $F^u_{\rho(g)}$ are injective immersions.
\end{enumerate}
\end{proposition}

\begin{proof}
That $F^s_{\rho(g)},F^u_{\rho(g)}$ and $F^c_{\rho(g)}$ are
subgroups, as well as claims $1$ and $2$ are consequences of the
fact $\fff^s_{\rho(g)},\fff^u_{\rho(g)}$ and $\fff^c_{\rho(g)}$
are Lie subalgebras of $\fh$.  This is true since if $v$ and $w$
are eigenvectors of $\Ad_L{\circ}{\pi_0}|_{\fh}(g)$ with
eigenvalues $\lambda$ and $\mu$, then $[v,w]$ is an eigenvector
with eigenvalue $\lambda\mu$.

We prove $3$ for $\sw^s_{\rho(g)}$, the proof is identical for
$\sw^u_{\rho(g)}$. Assume $3$ is false, then there is an element
of $f{\in}\Lambda{\cap}{h{\inv}}F^s_{\rho(g)}h$.  Since $f$ is in
$\Lambda$, $f$ is an element of $\pi_1(H/{\Lambda})$, which we can
represent by a curve $\bar f$ lying entirely in
${h{\inv}}F^s_{\rho(g)}h$. Since $\rho(g)$ is a contraction on
$F^s_{\rho(g)}$ and therefore ${h{\inv}}F^s_{\rho(g)}h$, for some
large $n$, the curve $\rho^n(g)\bar f$ is small and therefore
contractible, a contradiction, since $\rho^n(g)$ is a
diffeomorphism and $\bar f$ is not contractible in $H/{\Lambda}$.
\end{proof}

\noindent We now discuss the case of a quasi-affine action $\rho$.
We denote by $\hat \rho$ the affine action of $G$ from which
$\rho$ is defined.  We can define a splitting of $T(\base)$ as
$E^s_{\rho(g)}=H/{\Lambda}{\times}\fff^s_{(\hat \rho(g))},
H/{\Lambda}{\times}E^u_{\hat \rho(g)}=\fff^u_{(\hat \rho(g))}$ and
$E^c_{\hat \rho(g)}=H/{\Lambda}{\times}\fff^c_{(\hat
\rho(g))}{\times}TM$. Again it is easy to see that $\rho(g)$ is
partially hyperbolic with respect to this splitting whenever this
splitting is non-trivial. These sub-bundles are tangent to
foliations where $\sw^s_{\rho(g)}$ consists of $F^s_{\hat
\rho(g)}$ orbits, $\sw^u_{\rho(g)}$ consists of $F^u_{\hat
\rho(g)}$ orbits, and $\sw^c_{\rho(g)}$ consists of products of
$F^c_{\hat \rho(g)}$ orbits with $M$. It follows that all
dynamical foliations for any quasi-affine action are smooth.

 We define $E^0_{\rho}(G)$ to be the distribution
$H/{\Lambda}{\times}\fz{\times}TM$ which is tangent to the
foliation $\ff$. We state a lemma here which says that there are
finitely many elements in the acting group the sum of whose
(un)stable directions is the complement of $E^0_{\rho}(G)$ and
that, therefore, the intersection of their central distributions
is exactly $E^0_{\rho}(G)$.

\begin{lemma}
\label{lemma:basicdynamics}  There exits a finite set $\Phi$ of
elements in $G$ such that
$$T(\base)=H/{\Lambda}{\times}\big(\sum_{g{\in}\Phi}E^s_{\rho(g)}\big)\times\fz\times{TM}.$$
\end{lemma}

\begin{proof}
The derivative of $\rho$ on $T(\base)$ leaves invariant $TM$ and
$\fh$.   We let $\Ad$ be the adjoint representation of $L$. It
follows from the description of affine actions in subsection
\ref{subsection:describingaffineactions}, that $\fz$ is invariant
under $\Ad|_{\fh}(\pi(G))$, that there is an $\Ad|_{\fh}(\pi(G))$
invariant complement $\fm$ to $\fz$ and that $\Ad|_{\fm}(\pi(g))$
is  $\Ad|_{\fm}(\pi^E(g)\pi^K(g))$ where the $J$ representation
$\pi^E$  does not contain the trivial representation.

Recall from \cite[Subsection $3.8$]{FM1} that for any element $g$
of $GL_n(\mathbb R)$, there is a unique decomposition of $g=us=su$
where $u$ is unipotent and $s$ is semisimple. Further, we have a
unique decomposition $s=cp=pc$ where all eigenvalues of $p$ are
positive and all eigenvalues of $c$ have modulus one. We refer to
$p$ as the {\em polar part} of $g$ and denote it by $\pol(g)$. As
remarked there, one can define the polar part of an element for
elements of any real algebraic group and this definition is
independent of the realization of the group as an algebraic group.

By \cite[Lemma $3.21$]{FM1} there is a finite collection $\Psi$ of
elements in $G'$ whose polar parts are Zariski dense in $J$.
Combined with the fact that the representation
$\Ad{\circ}\pi^E|_{\fm}$ of $J$ does not contain invariant
vectors, this implies that $\fm{\cap}{\bigcap}\fff^c_{\rho(g)}=0$.
Letting $\Phi=\Psi{\cap}\Psi{\inv}$ completes the proof.
\end{proof}

As above, we let $\bar \rho$ be the action on $Z{\backslash}H$
defined by lifting $\rho$ to an action on $H{\times}M$ and looking
at the action on the leaves of the central foliation there.

To describe some further properties of the dynamics, we recall the
local product structure on $H{\times}M$ as a bundle over
$Z{\backslash}H$ from subsection \ref{subsection:semiconjugacy}.
Recall that the sub-bundle $H{\times}M{\times}\fm$ of
$T(H{\times}M)$ is a $G$ invariant complement to
$H{\times}\fz{\times}TM$.  Letting $\exp$ be the exponential map
for our fixed metric on $H{\times}M$ and letting
$B_{\fm}(0,\varepsilon)$ be the ball of radius $\varepsilon$  in
$\fm$, by choosing $\varepsilon$ small enough, we can guarantee
that $\exp_x(B_{\fm}(0,\varepsilon)$ defines a family of manifolds
transverse to the fibers of $p$.  Furthermore if we write
$W_x=\exp_x(B_{\fm}(0,\varepsilon))$ then, for small enough
$\varepsilon$, we have a local product structure on $H{\times}M$
given by
$$p{\inv}(W_x){\cong}B_{\fm}(0,\varepsilon){\times}p{\inv}{(x)}$$
$${\cong}B_{\fm}(0,\varepsilon){\times}(Z{\times}M).$$
We define a Riemannian metric on $Z{\backslash}H$ so the map
$p:\cover{\rightarrow}\quotient$ is a Riemannian submersion. The
next lemma says that for any small enough perturbation $\rho'$ of
$\rho$, points on the same transversal to $p{\inv}(x)$ can be
moved apart.

\begin{lemma}
\label{lemma:dynamics} There exists $c>0$ depending only on $\rho$
such that for $x$ in $H{\times}M$ and $y,z{\in}W_x$,  there is
$g{\in}\Phi$ and a non-negative integer $n$ such that
$d(p(\rho'(g)^nz,\rho'(g)^ny)>c$.
\end{lemma}

\begin{proof}
As this is a special case of a fairly standard fact from
(partially) hyperbolic dynamics, we merely sketch the proof.  It
suffices to consider points $y{\in}W_x$ with $d(x,y)<c$, since
otherwise the lemma is true for any $g{\in}\Phi$ and $n=0$.  We
assume $c$ is small enough that $B(x,c)$ is a convex, normal
neighborhood of $x$. Therefore the geodesic between $z$ and $y$ is
entirely contained in $B(x,c)$ and we can pull the geodesic back
to the $T_x(\base)$ where it can be approximated to first order by
a segment on a straight line which we denote by $V_{z,y}$.

Since $W_x$ is transverse to $\ff$, by Lemma
\ref{lemma:basicdynamics}, we can choose $g{\in}\Phi$ such that
the angle between $V_{z,y}$ and $\sw^{cs}_{\rho'(g)}(x)$ is
bounded away from zero for all sufficiently small $\rho'$.  The
dynamics of $\rho'(g)$ then force the angle between
$D\rho'(g)^nV_{z,y}$ and $\sw^{cs}_{\rho'(g)}(\rho'(g)^nx)$ to be
uniformly bounded away from zero.  This implies that
$\|D\rho'(g)^nV_{z,y}\|$ grows at an exponential rate controlled
by the uniform lower bound on the angle, on the constants $a,b$
for which $\rho(g)$ is $(a,b)$-partially hyperbolic and the $C^1$
size of the perturbation $\rho'$.  By choosing $c$  small enough,
the first order behavior of $d(\rho'(g)^nz,\rho'(g)^ny)$ is given
by $\|D\rho'(g)^nV_{z,y}\|$, so possibly after shrinking $c$
again, we can assume that when $\|D\rho'(g)^nV_{z,y}\|>2c$ then
$d(\rho'(g)^nx,\rho'(g)^ny)>c$.
\end{proof}

\begin{corollary}
\label{corollary:separatingpoints} There is a constant $c$
depending only on $\rho$ such that for any $x,y{\in}H{\times}M$
with $p(x){\neq}p(y)$ and any $\rho'$ sufficiently $C^1$ close to
$\rho$, there exists $g{\in}\Phi$ and a non-negative integer $n$
such that $$d(p(\rho(g)^n(x)),p(\rho(g)^n(y)))>c$$.
\end{corollary}

\begin{proof} Since we are only concerned with the distance
between projections under $p$, to prove the corollary, it suffices
to consider $y{\in}W_x$.  This case is immediate from Lemma
\ref{lemma:dynamics}.
\end{proof}

\subsection{\bf Fiber structure}
\label{section:fibers}

Throughout this subsection, we keep the notations and assumptions
of subsection \ref{subsection:dynamicsofaffine}. In this
subsection we show that the map $\theta$ is a homeomorphism. As an
immediate consequence, $\ff'$ is a foliation and $\theta$ is a
leaf conjugacy between $(\base,\rho',\ff')$ and
$(\base,\rho,\ff)$.

\begin{theorem}
\label{theorem:niceleaves} The map
$\theta:\base{\rightarrow}\base$ defined in subsection
\ref{subsection:semiconjugacy} is a homeomorphism.
\end{theorem}

\noindent Recall that our assumptions, stated in subsection
\ref{subsection:dynamicsofaffine}, imply that the $G$ action on
$\base$ defined by $\rho$ lifts to $\cover$.  As remarked in
subsection \ref{subsection:semiconjugacy}, if $\rho$ lifts to
$\cover$ then all small enough perturbations of $\rho$ also lift
to $\cover$.

\begin{proof}

We show that $\theta$ is a homeomorphism by showing it is a
homeomorphism on each fiber $f{\inv}(x)$.

We first show $\theta|_{f{\inv}(x)}$  is injective by a dynamical
argument. Assume $\theta|_{f{\inv}(x)}$  is not injective, then
there are two distinct points, $w,z{\in}f{\inv}(x)$ such that
$\theta(z)=\theta(w)$. This forces $p(z){\neq}p(w)$ since
otherwise $z=w$.  Since, $w,z{\in}W_x$ by definition, it follows
from Lemma \ref{lemma:basicdynamics} that for some $g{\in}G$ the
distance between $d(p({\rho'}(g)(z)), p({\rho'}(g)w))$ is greater
than some constant $c$ depending only on $\rho$.  Since $f$ can be
made arbitrarily $C^0$ close to $p$ by restricting to small enough
$C^1$ perturbations $\rho'$ of $\rho$, we know that
$d(f({\rho'}(g)z),f({\rho'}(g)w){\geq}\frac{c}{2}$ so
$f({\rho'}(g)z){\neq}f({\rho'}(g)w)$.  By equivariance of $f$ we
have ${\bar \rho}(g)f(z)\neq{\bar \rho}(g)f(w)$ which implies that
$f(z)\neq{f(w)}$, a contradiction.

We now show that $\theta|_{f{\inv}(x)}$ is surjective.  Let
$\pi_{U'}$ be composition of the restriction of $\theta$ to
$p{\inv}(U)$ composed with projection on the second coordinate. We
can then look at the set $V_y=\pi_{U'}{\inv}(y)$.   Given
$(x,y){\in}p{\inv}(U)$, we show that the map
$f:V_y{\rightarrow}{U}$ is onto. This has the desired implication,
since if $\theta|_{f{\inv}(x)}$ is not surjective then there
exists $y$ such that $\pi_{U'}{\inv}(y){\cap}f{\inv}(x)$ is empty
and therefore $x{\notin}f(V_y)$.  Since $f$ is $C^0$ close to $p$,
the map $\psi_U'$ is $C^0$ close to projection on the second
coordinate. Therefore, after identifying $V_y$ with a subset of
$(U',y)$ by a vertical projection, $f:V_y{\rightarrow}{U'}$ is
$C^0$ close to the identity map. Our result now follows from the
following general topological lemma.

\begin{lemma}
\label{lemma:browder} Let $B$ be the ball or radius $r$ about zero
in a Euclidean space $E$.  Let $F$ be any continuous map from $B$
into $E$ such that $d(F(x),x)<{\varepsilon}$ for all $x{\in}B$.
Then $F(B)$ contains the ball of radius $r-2{\varepsilon}$ about
$0$.
\end{lemma}

\begin{proof}
This generalizes the key point in the proof of the Browder fixed
point theorem: we assume the map is not surjective and use this to
construct a deformation retract from a closed ball onto it's
boundary. Let $B'$ be the ball of radius $r-2{\varepsilon}$ and
$x$ a point in $B'$.  Assume $x{\notin}F(B)$. Let
$B_{\varepsilon}=B(x,2{\varepsilon})$ and look at
$F(B_{\varepsilon}){\subset}B$.  Let $S^{k-1}$ be the boundary of
$B(x, 2{\varepsilon})$ Define a map from $B_{\varepsilon}$ to
$S^{k-1}$ by taking $y{\in}B_{\varepsilon}$ to $F(y)$ and then
projecting to $S^{k-1}$ along the ray from $x$ to $F(y)$. This
gives a continuous map $\bar F$ from $B^k$ to $S^{k-1}$ which,
when restricted to $S^{k-1}$ is $C^0$ close to the identity.  The
map $\bar F$ is $C^0$ close to the identity on $S^{k-1}$ and so is
of degree one and homotopic to the identity. Therefore we can
define a map from $B^k$ to $S^{k-1}$ which is the identity on
$S^{k-1}$ as follows. Take $B^k$ and embed it in a larger closed
ball $B_1^k$.  Let $S^{k-1}$ denote the boundary of $B^k$ and
$S_1^{k-1}$ denote the boundary of $B_1^{k}$.  Our map is defined
by first taking $\bar F$ to get a map from $B_1^{k}$ to
$B_1^{k}{\backslash}{\Int}(B^{k})$ and then compose with the
deformation retract from $B_1^{k}{\backslash}{\Int}(B^{k})$ to
$S_1^{k-1}$ described by the homotopy from $\bar F$ to the
identity. This then gives a new map $\bar F'$ which is a
deformation retract from $B_1^k$ to $S_1^{k-1}$. This is
impossible since ${\pi}_{k-1}(S^{k-1})={\mathbb Z}$ and
$\pi_{k-1}(B)=0$. Therefore $x{\in}F(B)$.
\end{proof}

\end{proof}

\subsection{The leaves of $\ff'$ are smooth}
\label{subsection:HPS}

This subsection is essentially a proof that $\ff'$ is a foliation
by $C^k$ leaves, as defined at the beginning of subsection
\ref{section:prelim3}. The degree of smoothness will depend on the
$C^1$ size of the perturbation $\rho'$ or, more precisely, on the
largest $r$ such that $\rho'(g)$ is $r$-normally hyperbolic for
all $g{\in}{\Phi}$ where $\Phi$ is as in Lemma
\ref{lemma:basicdynamics}.  For technical reasons both here and in
later sections, we let $\psi=\theta{\inv}$ and work with $\psi$
rather than $\theta$.

\begin{theorem}
\label{theorem:leafconjugacy} Given $k$, there is a neighborhood
$U$ of $\rho$  in the space $\Hom(G,\Diff^1(\base))$ such that for
any $\rho'{\in}U$ the homeomorphism
$\psi:H/{\Lambda}{\times}M{\rightarrow}H/{\Lambda}{\times}M$
defined above has the following properties.

\begin{enumerate}
\item The central foliation $\ff$ of $\base$ for the action $\rho$
is mapped by $\psi$ to a foliation $\ff'$ of $\base$ that is
central for the action $\rho'$.

\item The leaves of $\ff'$ are $C^k$ and $\psi$ is $C^k$ along
leaves with $k$-jet depending continuously on $\base$.

\item  the homeomorphism $\psi$ is a leaf conjugacy between
$(\base,\ff,\rho)$ and $(\base,\ff',\rho')$.

\item The map $\psi$ is $C^0$ close to the identity and also $C^k$
small along leaves,


\end{enumerate}
\end{theorem}

The main point is the improvement in the regularity of leaves of
$\ff'$, and all other conclusions follow quickly from this one.
The key fact is:

\begin{lemma}
\label{lemma:subfoliated} If $\psi$ is as defined above and
$\rho(g)$ is partially hyperbolic, then
$\psi(\sw^c_{\rho(g)})=\sw^c_{\rho'(g)},\psi(\sw^{cu}_{\rho(g)})=\sw^{cu}_{\rho'(g)}$
and $\psi(\sw^{cs}_{\rho(g)})=\sw^{cs}_{\rho'(g)}$.
\end{lemma}

\begin{proof}
The proof proceeds in two steps.  First we show that any leaf
$\sv'$ of $\sw^{cs}_{\rho'(g)}$ (resp. $\sw^{cu}_{\rho'(g)}$) is a
union of leaves of $\ff'$. Second, we show that for any leaf
$\sv'$ of $\sw^{cs}_{\rho'(g)}$ (resp. of
$\sw^{cu}_{\rho'(g)}$)there is a leaf $\sv$ of
$\sw^{cs}_{\rho(g)}$ (resp. of $\sw^{cu}_{\rho(g)}$) such that
$\sv'{\subset}\psi(\sv)$.  Interchanging the roles of $\rho'$ and
$\rho$ in the argument proves the reverse inclusion, forcing
$\psi(\sv)=\sv'$. Since any leaf of $\sw^c_{\rho(g)}$ is a
transverse intersection of leaves of $\sw^{cu}_{\rho(g)}$ and
$\sw^{cs}_{\rho(g)}$ and \cite[Theorem 6.8]{HPS} implies that any
leaf of $\sw^c_{\rho'(g)}$ is a transverse intersection of leaves
of $\sw^{cu}_{\rho'(g)}$ and $\sw^{cs}_{\rho'(g)}$ this
immediately implies that for any $\sv'$ of $\sw^{c}_{\rho'(g)}$,
there is a leaf $\sv$ of $\sw^{c}_{\rho(g)}$ such that
$\sv'=\psi(\sv)$.   To prove all of these statements, we will use
the construction of leaves of $\sw^{cs}_{\rho'(g)}$ (resp.
$\sw^{cu}_{\rho'(g)}$) from \cite{HPS}, which we recall in the
following paragraph.


Following \cite{HPS}, section $6$, we pick a smooth local
transversal $\eta$ to the tangent bundle $E^c_{\rho(g)}$ to the
foliation $\sw^c_{\rho(g)}$. As noted there, this can be chosen to
be a smooth approximation to $E^s_{\rho(g)}{\oplus}E^u_{\rho(g)}$.
Since in our setting, $E^s_{\rho(g)}{\oplus}E^u_{\rho(g)}$ is
smooth, we let $\eta=E^s_{\rho(g)}{\oplus}E^u_{\rho(g)}$.  We
denote by $V$ the manifold which is the disjoint union of all
leaves of $\sw^c_{\rho(g)}$.  Note $V$ does not have a countable
base and may not be separable, see \cite[Examples 2 and 2', page
68]{HPS} and following for related discussion. Let
$i:V{\rightarrow}\base$ be the inclusion and pull $\eta$ back to a
bundle $i^*\eta$. Note that there is a metric on $\eta$ and
therefore $i^*\eta$ defined by our choice of Riemannian metric on
$\base$ and let $i^*\eta(l)$ be the bundle of discs of radius $l$.
Then as described in \cite{HPS}, there are numbers $r>0$ and
$\varepsilon_0>0$ such that
$\exp{\circ}i_*:i^*\eta(\varepsilon_0){\rightarrow}\base$ is a
diffeomorphism when restricted to $i^*\eta|_{B(x,r)}$ where
$x{\in}V$ and $B(x,r)$ is a ball of radius $r$ in a leaf of
$\sw^c_{\rho(g)}$. (In \cite{HPS}, the sets $B(x,r)$ are replaced
by plaques of a plaquation of $V$.  In our context, i.e. when the
action of $\rho(g)$ is isometric along $V$, it is easy to see that
one can find a plaquation by small enough balls.)  Now for any
$\varepsilon<\varepsilon_0$, we can pull back the action of
$\rho(g)$ (resp. $\rho'(g)$) on $\base$ to a (partially defined or
overflowing)  $i^*\rho(g)$ (resp. $i^*\rho'(g)$) action on
$i^*\eta(\varepsilon)$, see \cite[pages 94-95]{HPS} for details.
As in \cite{HPS}, we define a submanifold $\tilde
\sw^{cs}_{\rho'(g)}$ of $i^*\eta(\varepsilon)$ by $\tilde
\sw^{cs}_{\rho'(g)}=\cap_{n{\leq}0}i^*\rho'(g)^n(i^*\eta(\varepsilon))$.
By \cite[p.107]{HPS} this is a $C^k$ submanifold of
$i^*\eta(\varepsilon)$ such that $\exp{\circ}i_*(\tilde
\sw^{cs}_{\rho'(g)})$ is the foliation $\sw^{cs}_{\rho'(g)}$.
Replacing $\rho'(g)$ by $\rho'(g{\inv})$, gives $\tilde
\sw^{cu}_{\rho'(g)}=\cap_{n{\geq}0}i^*\rho'(g)^n(i^\eta(\varepsilon))$
a $C^k$ submanifold such that $\exp{\circ}i_*(\tilde
\sw^{cu}_{\rho'(g)})$ is the foliation $\sw^{cu}_{\rho'(g)}$. As
in \cite[Theorem 6.8]{HPS} the intersection $\tilde
\sw^{cu}_{\rho'(g)}{\cap}\tilde \sw^{cs}_{\rho'(g)}$ is transverse
and the image of a section
$\sigma_{\rho'(g)}:V{\rightarrow}i^*{\eta(\varepsilon)}$ such that
$\exp{\circ}i_*{\circ}\sigma_{\rho'(g)}(\sv)$ is a leaf of
$\sw^{c}_{\rho'(g)}$.

Note that there is a  foliation $i^*\ff$ (resp. $i^*\ff'$) of
$i^*\eta(\varepsilon)$ defined on each component of $V$ by pulling
back $\ff{\cap}\exp(\eta(\varepsilon)|_{\sv})$ where
$\sv{\in}\sw^c_{\rho(g)}$.  Note that we consider the leaves of
these foliations to be connected components of pre-images of
leaves rather than entire pre-images of leaves. This foliation is
preserved by $i^*{\rho(g)}$ (respectively $i^*{\rho'(g)}$) for any
$g$ in $G$.

Since $\psi$ is a $C^0$ small homeomorphism, for any leaf
$\sv{\in}\sw^c_{\rho(g)}$ and any $x{\in}\sv$, we have
$$\psi(\exp(\eta(\frac{\varepsilon}{2})|_{B_{\sv}(x,\frac{r}{2})}))\subset\exp(\eta(\varepsilon)|_{B_{\sv}(x,r)})$$
and so we can pull back $\psi$ to a map
$i^*{\psi}:i^*\eta(\frac{\varepsilon}{2}){\rightarrow}i^*\eta(\varepsilon)$.

We now have the following diagram of $\Za$ actions:

$$\xymatrix{
(i^*\eta(\frac{\varepsilon}{2}),i^*\rho(g))
\ar[d]_{i^*\psi}\ar[r]_{\exp{\circ}i_*} &
({\base},\rho(g))\ar[d]_{\psi}
\\
(i^*(\eta)(\varepsilon),i^*{\rho}'(g))\ar[r]_{\hskip .2in
\exp{\circ}i_*}&(\base,\rho'(g))\\}$$

\noindent where the horizontal arrows are equivariant and the
vertical arrows are leaf conjugacies.

We first show that each leaf of $\sw^{cs}_{\rho'(g)}$ is a union
of leaves of $\ff'$.  Given a leaf $\su$ of $\ff$, we can find a
leaf $\sL$ of $\sw^c_{\rho(g)}$ such that $\su{\subset}\sL$.  Note
that for the leaf $\sL$ of $\sw^c_{\rho(g)}$, we have that $\tilde
\sw^{cs}_{\rho(g)}{\cap}i^*\eta(\varepsilon)|_{\sL}=(\exp{\circ}i_*){\inv}(\sv){\cap}i^*\eta(\varepsilon)|_{\sL}$
where $\sv$ is the unique leaf of $\sw^{cs}_{\rho(g)}$ containing
$\sL$.
Furthermore, from the definition of $\tilde \sw^{cs}_{\rho(g)}$ it
then follows that
$$(i^*{\circ}{\exp}){\inv}(\su){\cap}i^*\eta(\varepsilon)|_{\sL}{\subset}i^*\rho(g)^n(i^*\eta(\frac{\varepsilon}{2}))$$
for all $n<0$.  Since $i^*\psi$ is a leaf conjugacy, possibly
after shrinking $\varepsilon$, we have
$$i^*\psi((i^*{\circ}{\exp}){\inv}(\su){\cap}i^*\eta(\varepsilon)|_{\sL}{\subset}i^*\rho'(g)^n(i^*\eta({\varepsilon}))$$
for all $n$.  This implies that
$$i^*\psi((i^*{\circ}{\exp}){\inv}(\su){\cap}i^*\eta(\varepsilon)|_{\sL})=(i^*{\circ}\exp){\inv}(\psi(\su)){\cap}i^*\eta(\varepsilon)|_{\sL}$$
is contained in $\tilde
\sw^{cs}_{\rho'(g)}{\cap}i^*\eta(\varepsilon)|_{\sL}$.  This then
implies that $\psi(\su)$ is contained in $\sw^{cs}_{\rho'(g)}$ as
desired.

Now the fact that $i^*\psi$ is a leaf conjugacy, the definition of
$\tilde \sw^{cs}_{\rho'(g)}$ and $\tilde \sw^{cs}_{\rho(g)}$ and
the fact that each leaf of $\sw^{cs}_{\rho'(g)}$ (resp.
$\sw^{cs}_{\rho(g)}$) is a union of leaves of $\ff'$ (resp. $\ff$)
implies that $i^*{\psi}(\tilde \sw^{cs}_{\rho(g)}){\subset}\tilde
\sw^{cs}_{\rho'(g)}$.  This implies that for any
$\sv{\in}\sw^{cs}_{\rho}(g)$ there is
$\sv'{\in}\sw^{cs}_{\rho'}(g)$ such that $\psi(\sv){\subset}\sv'$.
Interchanging the roles of $\rho$ and $\rho'$ and replacing $\psi$
by $\psi{\inv}$, the same argument proves that for any
$\sv'{\in}\sw^{cs}_{\rho'}(g)$ there is
$\sv''{\in}\sw^{cs}_{\rho}(g)$ such that
$\psi{\inv}(\sv'){\subset}\sv''$.  These two facts then imply that
for any $\sv{\in}\sw^{cs}_{\rho}(g)$ there is
$\sv'{\in}\sw^{cs}_{\rho'}(g)$ such that $\psi(\sv)=\sv'$.

A similar argument using $g{\inv}$ in place of $g$ implies that
for any $\sv{\in}\sw^{cu}_{\rho(g)}$ there is
$\sv'{\in}\sw^{cu}_{\rho'(g)}$ such that $\psi(\sv)=\sv'$. As
remarked above, it then follows that for any leaf $\sv$ of
$\sw^{c}_{\rho(g)}$ we have $\psi(\sv)={\sv'}$ where $\sv'$ is a
leaf of $\sw^{c}_{\rho'(g)}$ .


\end{proof}

\noindent{\bf Remark:} The proof of Lemma \ref{lemma:subfoliated}
does not depend on all of our assumptions and the statement could
be made more axiomatic.  All we require is that $\psi$ is a leaf
conjugacy and that any leaf of  $\sw^c_{\rho(g)}$ is a union of
leaves of $\ff$ for any $g{\in}\Phi$.

We recall two definitions and a lemma from \cite{MQ} page 145-6:

\begin{defn}
\begin{enumerate}
\item Let $N$ be a smooth Riemannian manifold and $N_1,N_2$ two
immersed $C^k$ manifolds.
 We say that $N_1$ and $N_2$ intersect {\it $s$-transversely} if $N_1{\cap}N_2$ a manifold $N'$
of dimension $\dim(TN_1(x){\cap}TN_2(x))$ for any
$x{\in}N_1{\cap}N_2$.

\item Let $N$ be a smooth Riemannian manifold and
$N_1,{\ldots},N_l$ a collection of $C^k$ immersed submanifolds. We
say that the family $N_1,{\ldots},N_l$ intersects {\it
$s$-transversely} if $\cap_{i=1}^{j-1}N_i$ intersects
$s$-transversely with $N_j$ for $j=2,...l$.
\end{enumerate}
\end{defn}

\begin{defn}
Let $N$ be a smooth manifold and $N_1,{\ldots},N_l$ a collection
of $C^k$ manifolds.  We call a collection $N_1',{\ldots},N_l'$ of
$C^k$ submanifolds of $N$ a {\em topologically trivial
$\varepsilon$-perturbation} of $N_1,{\ldots},N_l$ if there is a
homeomorphism $h:N{\rightarrow}N$ such that
\begin{enumerate}
\item $N_i'=h(N_i)$ \item $d(h(x),x)<\varepsilon$ for all
$x{\in}N$, and \item  for any $x{\in}N$ and any $i$, let $B_i(x)$
(resp. $B'_i(h(x))$) be the unit balls in the tangent space
$TN_i(x)$ (resp. $TN'_i(h(x))$). Then $B_i(x_i)$ and $B'_i(h(x))$
are within an $\varepsilon$ neighborhood of each other in $TN$.
\end{enumerate}
\end{defn}

\noindent{\bf Remark:} On \cite[p.145]{MQ}, the same notion is
called an $\varepsilon$-perturbation.  We choose to modify the
terminology, since the original terminology is somewhat deceptive.

\begin{lemma}
\label{lemma:stransverse}  Let $N$ be a compact manifold and
$N_1,{\ldots},N_l$ be $C^k$ submanifolds intersecting
$s$-transversely, then

\begin{enumerate}
\item $\cap_{i=1}^{l}N_i$ is a $C^k$ submanifold, and

\item there exists $\varepsilon>0$ depending only on $N$, such
that if $N'_1,{\ldots},N'_l$ is any topologically trivial
$\varepsilon$-perturbation $N'_1,{\ldots},N'_l$ of
$N_1,{\ldots},N_l$, then $N'_1,{\ldots},N'_l$ intersect
$s$-transversely.
\end{enumerate}
\end{lemma}

Point $(1)$ for $l=2$ is Lemma 5.5(1) of \cite{MQ}, where $C^k$
replaces the word smooth. The proof is the same.  As noted on page
$146$ of \cite{MQ}, the case of $l>2$ follows by induction.
Similarly part $(2)$ follows from \cite[Lemma 5.5(2)]{MQ} and
induction.  The proof of \cite[Lemma 5.5(2)]{MQ} implicitly uses
that if $N_1,N_2$ are $C^k$ submanifolds of $N$ and $N_1',N_2'$
are a topologically trivial $\varepsilon$-perturbation of
$N_1,N_2$ then $\dim(TN_1'(x){\cap}TN_2'(x))=\dim(N_1'{\cap}N_2')$
for every $x$ in $N$.  A priori $\dim(TN_1'(x){\cap}TN_2'(x))$
could drop, but it is in fact bounded below by
$\dim(N_1'{\cap}N_2')$.  This can be deduced from standard facts
about transversality using an argument similar to the proof of
\cite[Lemma 5.5(1)]{MQ}.  This is not noted explicitly in
\cite{MQ}.

Since we will need to know not just that the leaves of $\ff'$ are
$C^k$ submanifolds but that $\ff'$ is a foliation by $C^k$ leaves,
we require a slight strengthening of Lemma
\ref{lemma:stransverse}, also remarked on \cite[page $146$]{MQ}.
Let $\ff_1,{\ldots}\ff_l$ be foliations by $C^k$ leaves of a
compact manifold $N$.  We say that the $\ff_i$ intersect
$s$-transversely, if for each $x{\in}N$, leaves $\ff_i(x)$
intersect $s$-transversely and the dimension of the intersections
$\cap_{i=1}^j{\ff_i(x)}$ is independent of $x$ for any $j$ from
$2$ to $l$.  We say that a collection of foliations
$\ff_1',\ldots,\ff_l'$ of $N$ is an $\varepsilon$-perturbation of
$\ff_1,\ldots,\ff_l$  if:
\begin{enumerate}
\item there exists a homeomorphism $h:N{\rightarrow}N$ with
$h(\ff_i)=\ff_i'$ for $i$ from $1$ to $l$ and
$d(h(x),x)<\varepsilon$ for ever $x{\in}N$; \item for any
$x{\in}N$ and any $i$, let $B_i(x)$ (resp. $B'_i(h(x))$) be the
unit balls in the tangent space $T\ff_i(x)$ (resp.
$T\ff'_i(h(x))$). Then $B_i(x_i)$ and $B'_i(h(x))$ are within an
$\varepsilon$ neighborhood of each other in $TN$.
\end{enumerate}

As remarked in \cite{MQ}, a slight modification of the proof of
Lemma \ref{lemma:stransverse} shows that
\begin{enumerate}
\item if $\ff_1,{\ldots},\ff_l$ are $s$-transverse foliations by
$C^k$ leaves, then the foliation defined by intersections of
leaves of $\ff_1,{\ldots},\ff_l$ is a foliation by $C^k$ leaves,
and;

\item if $N$ is compact, there exists $\varepsilon>0$ such that
any $\varepsilon$-perturbation of $\ff_1,{\ldots},\ff_l$ is
$s$-transverse.
\end{enumerate}

\begin{lemma}
\label{lemma:smoothleaves} Given $k$, if $\rho'$ is a sufficiently
$C^1$ small, $C^k$ perturbation of $\rho$, the leaves of the
foliation $\ff'$ are the $s$-transverse intersections of leaves of
$\sw^c_{\rho'(g)}$ for $g{\in}\Phi$. Therefore $\ff'$ is a
foliation by $C^k$ leaves. Furthermore, the foliation $\ff'$ is
close to $\ff$ in the natural topology on foliations by $C^k$
leaves.
\end{lemma}

\begin{proof}
We fix a neighborhood $U$ of $\rho$ in $\Hom(D,\Diff^1(M))$ such
that

\begin{enumerate}
\item for $g{\in}{\Phi}$, $\rho(g')$ is close enough to $\rho(g)$
to satisfy the hypotheses of Theorem \ref{theorem:hps},

\item the map $\psi=\theta{\inv}$ constructed in subsection
\ref{subsection:semiconjugacy} satisfies
$d(\psi(x),x)<\varepsilon$ (or equivalently
$d(\theta(x),x)<\varepsilon$) for $\varepsilon$ as in Lemma
\ref{lemma:stransverse}.
\end{enumerate}

Let $\fL'$ be an arbitrary leaf of $\ff'$ and let
$\fL=\psi{\inv}(\fL')$.  Note that the leaf $\fL$ of $\ff$ is the
$s$-transverse intersection of leaves $\sv^c_{\rho(g)}$ of
$\sw^c_{\rho(g)}$ for $g{\in}\Phi$. By Lemma
\ref{lemma:subfoliated}, for every $g{\in}{\Phi}$, we know that
$\psi(\sv^c_{\rho(g)})$ is a leaf $\sv^c_{\rho'(g)}$ of
$\sw^c_{\rho'(g)}$ and therefore by Theorem \ref{theorem:hps} a
$C^k$ submanifold of $\base$ which is $C^k$ close to some leaf of
$\sw^c_{\rho(g)}$. Since $\psi$ is a homeomorphism and can be made
arbitrarily small by choosing $\rho'$ close enough to $\rho$,
Lemma \ref{lemma:stransverse}$(2)$ implies that
$\psi(\sv^c_{\rho(g)})=\sv^c_{\rho'(g)}$ intersect
$s$-transversely in a $C^k$ manifold.  Since
$$\cap_{d{\in}\Phi}\psi(\sv^c_{\rho(g)})=\psi(\cap_{d{\in}\Phi}
\sv^c_{\rho(g)})$$
$$=\psi(\fL)=\fL'$$
it follows that every leaf $\fL'$ of $\ff'$ is a $C^k$ submanifold
of $\base$.  The remarks following Lemma \ref{lemma:stransverse}
the imply that $\ff'$ is a foliation  by $C^k$ leaves.
\end{proof}

\begin{proof}[Proof of Theorem \ref{theorem:leafconjugacy}]
The homeomorphism $\psi{\inv}=\theta$ is constructed in subsection
\ref{subsection:semiconjugacy} and shown to be a homeomorphism in
Theorem \ref{theorem:niceleaves}. Since $\psi{\inv}$ is given by
projecting from leaves of $\ff'$ to leaves of $\ff$ via a smooth
transversal, and leaves of $\ff'$ are $C^k$ by Lemma
\ref{lemma:smoothleaves}, the map $\psi$ is $C^k$ and $C^k$ small
along fibers.

The remaining conclusions follow from Corollary
\ref{corollary:semiconjugacy} and Theorem
\ref{theorem:niceleaves}.

\end{proof}

We will eventually need one additional fact concerning $\psi$
which is now straightforward. To state this fact about $\psi$, we
need to define some additional dynamical foliations. Let
$E^0_{\rho'}(G)$ be the distribution tangent to $\ff'$.  Given
$g{\in}\Phi$, we take the distributions
$E^0_{\rho}(G){\oplus}E^s_{\rho(g)}$ and
$E^0_{\rho'}(G){\oplus}E^s_{\rho'(g)}$.   Recall that $\ff$ is
tangent to $E^0_{\rho}(G)$.  For the $\rho$ action there is a
smooth foliation tangent to $E^0_{\rho}(G){\oplus}E^s_{\rho(g)}$,
which we denote by ${\ff}{\oplus}\sw^s_{\rho(g)}$.  To see this,
one notes that the group $Z$ normalizes the group $F^s_{\rho(g)}$
and so the product $ZF^s_{\rho(g)}$ is a subgroup of $H$.  For the
case of affine actions, the foliation
${\ff}{\oplus}\sw^s_{\rho(g)}$ is just the orbit foliation for the
left action of $ZF^s_{\rho(g)}$ on $H/{\Lambda}$.  For
quasi-affine actions, we recall that there is a natural projection
$H/{\Lambda}{\times}M{\rightarrow}M$ and the foliation
${\ff}{\oplus}\sw^s_{\rho(g)}$ is given by pre-images in $\base$
of the $ZF^s_{\rho(g)}$ orbits in $H/{\Lambda}$. We note that the
lift of any leaf of ${\ff}{\oplus}\sw^s_{\rho(g)}$ to $\cover$ is
of the form of $p{\inv}(\sv)$ where $\sv$ is a leaf of
$\sw^s_{\bar \rho(g)}$.

\begin{proposition}
\label{proposition:funnyleaves} For every $g{\in}\Phi$, there is a
$\rho'$ invariant foliation $\ff'{\oplus}\sw^s_{\rho'(g)}$ of
$\base$ tangent to $E^0_{\rho'}(G){\oplus}E^s_{\rho'(g)}$ such
that
$\psi({\ff}{\oplus}\sw^s_{\rho(g)})={\ff'}{\oplus}\sw^s_{\rho'(g)}$.
\end{proposition}

\begin{proof}
We can define the foliation $\ff'{\oplus}\sw^s_{\rho'(g)}$ to be
$\psi(\ff{\oplus}\sw^s_{\rho(g)})$.  A leaf $\sv$ of the foliation
$\ff{\oplus}\sw^s_{\rho(g)}$ is given by sets of points sharply
forward asymptotic to a leaf $\sL$ of $\ff$. Here, as in
\cite{HPS}, $x$ is sharply forward asymptotic to $\fL$ if
$d(\rho^n(g)x,\rho^n(\sL))$ goes to zero at least as fast as
$\exp(-\lambda{n})$ for some $\lambda>0$ depending on the
dynamics. Since $\psi$ is a leaf conjugacy it follows that a leaf
$\sv'$ of $\ff'{\oplus}\sw^s_{\rho'(g)}$ is the given by sets of
points sharply forward asymptotic to a leaf $\sL'$ of $\ff'$. By
\cite[Theorem 6.8 (e) and (f)]{HPS}, the leaves of
$\sw^s_{\rho'(g)}$ are exactly sets of points sharply forward
asymptotic to the orbit of a point on a leaf of
$\sw^c_{\rho'(g)}$.  Therefore, a leaf $\sv'$ of
$\ff'{\oplus}\sw^s_{\rho'(g)}$ is the union of all leaves of
$\sw^s_{\rho'(g)}$ through a leaf $\fL'$ of $\ff'$.  This
immediately implies that $\ff'{\oplus}\sw^s_{\rho'(g)}$ is tangent
to $E^0_{\rho'}(G){\oplus}E^s_{\rho'(g)}$ and completes the proof.
\end{proof}

\section{\bf Property T and conjugacy}
\label{section:conjugacy}

In this section, we modify the leaf conjugacy obtained at the end
of the last section to obtain a semiconjugacy.  The a priori
regularity of this semiconjugacy will be somewhat bad.  In section
\ref{section:smoothing} we show it is a homeomorphism, in section
\ref{section:katokspatzier} we show it is differentiable along
many foliations and in section \ref{section:finalarguments} we
show it is differentiable and even $C^{\infty}$ when $\rho'$ is
$C^{\infty}$. The key ingredient in the arguments of this section
is \cite[Theorem $2.16$]{FM2}, so we begin by recalling some
notation and definitions from subsection $2.3$ of that paper. For
most of this section $G$ will be a compactly generated topological
group, though for our applications, $G$ will be $J$ or $\G$ as
above.

Throughout this section $X$ will be a second countable, compact,
Hausdorff manifold and $\ff$ will be a foliation of $X$ by $C^k$
leaves. For background on foliated spaces, their tangent bundles,
and transverse invariant measures, the reader is referred to
\cite{CC} or \cite{MS}.

We let $\Diff^k(X,\ff)$ be the set of homeomorphisms of $X$ which
preserve $\ff$ and restrict to $C^k$ diffeomorphisms on each leaf
with derivatives depending continuously on $x$ in $X$.  For
$1{\leq}k{\leq}{\infty}$, there is a natural $C^k$ topology on
$\Diff^k(X,\ff)$. The definition of this topology is
straightforward and is recalled in \cite[Subsection $2.3$]{FM2}.

We now define a special class of perturbations of actions.

\begin{defn}
\label{definition:foliatedperturbation} Let $G$ be a compactly
generated, topological group and $\rho$ an action of $G$ on $X$
defined by a homomorphism from $G$ to $\Diff^{\infty}(X,\ff)$. Let
$\rho'$ be another action of $G$ on $X$ defined by a homomorphism
form $g$ to $\Diff^k(X,\ff)$. Let $U$ be a (small) neighborhood of
the identity in $\Diff^k(X,\ff)$ and $K$ be a compact generating
set for $G$. We call $\rho'$ a {\em $(U,C^k)$-foliated
perturbation} of $\rho$ if:
\begin{enumerate}
\item  for every leaf $\fL$ of $\ff$ and every $g{\in}G$, we have
$\rho(g)\fL=\rho'(g)\fL$ and,

\item $\rho'(g)\rho(g){\inv}$ is in $U$ for every $g$ in $K$.

\end{enumerate}
\end{defn}

We fix a continuous, leafwise smooth Riemannian metric $g_{\ff}$
on $T{\ff}$, the tangent bundle to the foliation and note that
$g_{\ff}$ defines a volume form and corresponding measure on each
leaf $\fL$ of $\ff$, both of which we denote by $\nu_{\ff}$.
(Metrics $g_{\ff}$ exist by a standard partition of unity
argument.)  Let $G$ be a group and $\rho$ an action of $G$ on $X$
defined by a homomorphism from $G$ to $\Diff^k(X,\ff)$. We say the
action is {\em leafwise isometric} if $g_{\ff}$ is invariant under
the action.  When $G=\Za$ and $\Za=<f>$, we will call $f$ a {\em
leafwise isometry}.

For the remainder of this section, we will assume that the
foliation has a transverse invariant measure $\nu$. By integrating
the transverse invariant measure $\nu$ against the Riemannian
measure on the leaves of $\ff$, we obtain a measure $\mu$ on $X$
which is finite when $X$ is compact. In this case, we normalize
$g_{\ff}$ so that $\mu(X)=1$. We will write $(X,\ff,g_{\ff},\mu)$
for our space equipped with the above data, sometime leaving one
or more of $\ff, g_{\ff}$ and $\mu$ implicit.  We will refer to
the subgroup of $\Diff^k(X,\ff)$ which preserves $\nu$ as
$\Diff_{\nu}^k(X,\ff)$.  Note that if $\rho$ is an action of $G$
on $X$ defined by a homomorphism into $\Diff_{\nu}^k(X,\ff)$ and
$\rho$ is leafwise isometric, then $\rho$ preserves $\mu$.
Furthermore if $\rho$ is an action of $\G$ on $X$ defined by a
homomorphism into $\Diff_{\nu}^k(X,\ff)$ and $\rho'$ is a
$(U,C^k)$-leafwise perturbation of $\rho$, then it follows easily
from the definition that $\rho'$ is defined by a homomorphism into
$\Diff^k_{\nu}(X,\ff)$ since the induced map on transversals is
the same.

Before stating one of the main results of \cite{FM2}, we will need
a coarse quantitative measure of the $C^k$ size of $C^k$ map. We
denote by $B_{\ff}(x,r)$ the ball in $\fL_x$ about $x$ of radius
$r$. For a sufficiently small value of $r>0$, we can canonically
identify each $B_{\ff}(x,2r)$ with the ball of radius $2r$ in
Euclidean space via the exponential map from $T\ff_x$ to $\fL_x$.
We first consider the case when $k$ is an integer, where we can
give a pointwise measure of size. Recall that a $C^k$ self map of
a manifold $Z$ acts on $k$-jets of $C^k$ functions on $Z$. Any
metric on $TZ$ defines a pointwise norm on each fiber of the
bundle of $J^k(Z)$ of $k$-jets of functions on $Z$. For any $C^k$
diffeomorphism $f$ we can define $\|j^k(f)(z)\|$ as the operator
norm of the map induced by $f$ from $J^k(Z)_z$ to $J^k(Z)_{f(z)}$.
For a more detailed discussion on jets and an explicit
construction of the norm on $J^k(Z)_z$, see \cite[Section
$4$]{FM2}.  We say that a map $f$ has $C^k$ size less than
$\delta$ on a set $U$ if $\|j^k(f)(z)\|<\delta$ for all $z$ in
$U$. If $k$ is not an integer, we say say that $f$ has $C^k$ size
less than $\delta$ on $U$ if $f$ has $C^{k'}$ size less than
$\delta$ on $U$ where $k'$ is the greatest integer less than $k$
and $j^{k'}(f)$ satisfies a (local) H\"older estimate on $U$. See
\cite[Section $4$]{FM2} for a more detailed discussion of H\"older
estimates.

\noindent{\bf Remark:} This notion of $C^k$ size is not very
sharp.  The size of the identity map will be $1$, as will be the
size of any isometry of the metric.  We only use this notion of
size to control estimates on a map at points where the map is
known to be ``fairly large" and where we only want bounds to show
it is ``not too large".

For the following theorem, we assume that the holonomy groupoid of
$(X,\ff)$ is Hausdorff. This is a standard technical assumption
that allows us to define certain function spaces on ``pairs of
points on the same leaf of $(X,\ff)$".  See \cite[Subsection
$6.1$]{FM2}, \cite{CC} and \cite{MS} for further discussion. All
the foliations considered in this paper for the proof of Theorem
\ref{theorem:main} are covered by fiber bundles, and in that case
 the holonomy groupoid is Hausdorff.  We now recall \cite[Theorem
$2.16$]{FM2}.

\begin{theorem}
\label{theorem:almostconjugacygen} Let $G$ be a locally compact,
$\sigma$-compact group with property $(T)$.  Let $\rho$ be a
continuous leafwise isometric action of $G$ on $X$ defined by a
homomorphism from $G$ to $\Diff_{\nu}^{\infty}(X,\ff)$. Then for
any $k{\geq}3,\kappa>0$ and any $\varsigma>0$ there exists a
neighborhood $U$ of the identity in $\Diff^k(X,\ff)$ such that for
any continuous $(U,C^k)$-foliated perturbation $\rho'$ of $\rho$
there exists a measurable map $\phi:X{\rightarrow}X$ such that:

\begin{enumerate}
\item $\phi{\circ}\rho(g)=\rho'(g){\circ}\phi$ for all $g{\in}G$,

\item $\phi$ maps each leaf of $\ff$ into itself,


\item there is a subset $S{\subset}X$ with $\mu(S)=1-\varsigma$
and $\G{\cdot}S$ has full measure in $X$, and a constant
$r{\in}\Ra^+$, depending only on $X,\ff$ and $g_{\ff}$, such that,
for every $x{\in}S$, the map $\phi:B_{\ff}(x,r){\rightarrow}\fL_x$
is $C^{k-1-\kappa}$-close to the identity; more precisely, with
our chosen identification of $B_{\ff}(x,2r)$ with the ball or
radius $2r$ in Euclidean space,
$\phi-\Id:B_{\ff}(x,r){\rightarrow}B_{\ff}(x,2r)$ has
$C^{k-1-\kappa}$ norm less than $\varsigma$ for every $x{\in}S$,
and

\item there exists $0<t<1$ depending only on $G$ and $K$ such that
the set of $x{\in}X$ where the $C^{k-1-\kappa}$ size of $\phi$ on
$B_{\ff}(x,r)$ is not less than $(1+\varsigma)^{l+1}$ has measure
less than $t^l\varsigma$ and,

\end{enumerate}
Furthermore, for any $l{\geq}k$, if $\rho'$ is a $C^{2l-k+1}$
action, then by choosing $U$ small enough, we can choose $\phi$ to
be $C^{l}$ on $B_{\ff}(x,r)$ for almost every $x$ in $X$.  In
particular, if $\rho'$ is $C^{\infty}$ then for any $l{\geq}k$, by
choosing $U$ small enough, we can choose $\phi$ to be $C^l$ on
$B_{\ff}(x,r)$ for almost every $x$ in $X$.
\end{theorem}

\noindent {\bf Remarks:} \begin{enumerate} \item Since $\rho'$ is
a foliated perturbation of $\rho$, the transverse measure $\nu$ is
$\rho'$ invariant.  This is because $\rho'$ defines the same
action on transversals as $\rho$.  \item The map $\phi$
constructed in the theorem is not even $C^0$ close the identity on
$X$. However, the proof of the theorem shows that for every
$1{\leq}q{<}{\infty}$, possibly after changing $U$ depending on
$q$, we have $\int_X(d(x,\phi(x))^qd\mu{\leq}\varsigma$.
\end{enumerate}

We now proceed to show how Theorem
\ref{theorem:almostconjugacygen} can be applied in the proof of
Theorem \ref{theorem:main}.  As before, we fix a semisimple Lie
group $J$ with all simple factors of real rank at least two and a
lattice $\G$ in $J$ and let $G$ be one of $J$ and $\G$.  We also
fix an algebraic group $H$, a cocompact lattice $\Lambda<H$, a
compact manifold $M$ and a quasi-affine action $\rho$ of $G$ on
$\base$. Once again, we assume that $\rho$ lifts to $\cover$. We
fix the foliation $\ff$ of $\base$ by central leaves for $\rho(G)$
as in subsection \ref{subsection:describingaffineactions}. We
further assume that the $G$ action defined by $\rho$ lifts to an
action on $\cover$. We note that there is a transverse invariant
measure $\nu$ to $\ff$ defined by lifting to $\cover$ and
identifying local transversals with their projections to
$\quotient$.

\begin{proposition}
\label{proposition:foliatedperturbation} Let
$U\subset\Hom(G,\Diff_{\nu}^k(\base,\ff))$ be a neighborhood of
$\rho$, there is a neighborhood $V$ of $\rho$ in
$\Hom(G,\Diff^1(\base))$ such that if $\rho'{\in}V$ is a $C^k$
action and $\psi$ is the homeomorphism from Theorem
\ref{theorem:leafconjugacy}, then
$\psi{\inv}{\circ}{\rho}{\circ}{\psi}$ is in $U$.  Furthermore,
given $m{\geq}k$, by assuming $\rho'$ is $C^{m}$ and possibly
after shrinking $V$, we can also guarantee that
$\psi{\inv}{\circ}{\rho}{\circ}{\psi}$ is in
$\Hom(G,\Diff^l_{\nu}(\base,\ff))$.
\end{proposition}

\begin{proof}
This is immediate from the definitions and Theorem
\ref{theorem:leafconjugacy}.
\end{proof}

\begin{theorem}
\label{theorem:candidateconjugacy} For every $k{\geq}3, \kappa>0$
and $\varsigma>0$, there is a neighborhood $V$ of $\rho$ in
$\Hom(G,\Diff^k(\base))$ such that if $\rho'{\in}V$ then there
exists a measurable map $\varphi:\base{\rightarrow}\base$ such
that:
\begin{enumerate}
\item $\varphi{\circ}\rho(g)=\rho'(g){\circ}\varphi$ for all
$g{\in}G$,

\item $\varphi$ maps each leaf of $\ff$ into a leaf of $\ff'$,

\item
$\varphi(\ff{\oplus}\sw^{s}_{\rho(g)})=\ff{\oplus}\sw^{s}_{\rho'(g)}$
on a set of full measure in $\base$ for any $g{\in}G$,

\item there is a subset $S{\subset}X$ with $\mu(S)=1-\varsigma$
and $\G{\cdot}S$ is of full measure in $\base$, and a constant
$r{\in}\Ra^+$, depending only on $X,\ff$ and $g_{\ff}$, such that,
for every $x{\in}S$, the map
$\varphi:B_{\ff}(x,r){\rightarrow}\ff(x)$ is
$C^{k-1-\kappa}$-close to the identity; more precisely, with our
chosen identification of $B_{\ff}(x,2r)$ with the ball or radius
$2r$ in Euclidean space,
$\varphi-\Id:B_{\ff}(x,r){\rightarrow}B_{\ff}(x,2r)$ has
$C^{k-1-\kappa}$ norm less than $\varsigma$ for every $x{\in}S$,
and

\item there exists $0<t<1$ depending only on $\G$ and $K$ such
that the set of $x{\in}X$ where the $C^{k-1-\kappa}$ size of
$\varphi$ on $B_{\ff}(x,r)$ is not less than $(1+\varsigma)^{l+1}$
has measure less than $t^l\varsigma$ and,

\end{enumerate}
Furthermore, for any $l{\geq}k$, if $\rho'$ is a $C^{2l-k+1}$
action, then by choosing $U$ small enough, we can choose $\varphi$
to be $C^{l}$ on $B_{\ff}(x,r)$ for almost every $x$ in $X$.  In
particular, if $\rho'$ is $C^{\infty}$ then for any $l{\geq}k$, by
choosing $U$ small enough, we can choose $\phi$ to be $C^l$ on
$B_{\ff}(x,r)$ for almost every $x$ in $X$.
\end{theorem}

\noindent{\bf Remark:} Conclusion $(4)$ combined with equivariance
of $\varphi$ and the fact that the central foliation is the
quotient of the fibers of a bundle, imply that for almost every
$x$, the map $\varphi$ is $C^{k-1-\kappa}$ along $\ff(x)$ and that
the derivative $D\varphi:T\ff(x){\rightarrow}T\ff(x)$ is an
isomorphism at all points of $\ff(x)$.

\begin{proof}
By Proposition \ref{proposition:foliatedperturbation}, we can
apply Theorem \ref{theorem:almostconjugacygen} to the actions
$\rho$ and $\psi\inv{\circ}\rho'{\circ}{\psi}$.  This produces a
map $\phi$ satisfying the conclusions of Theorem
\ref{theorem:almostconjugacygen}.  We let
$\varphi=\psi{\circ}\phi$ which satisfies $(1)$ and $(2)$ by
construction. Since $\psi$ is uniformly $C^k$ small when
restricted to any leaf of $\ff$, the estimates in $(4)$ and $(5)$
follow from the estimates in Theorem
\ref{theorem:almostconjugacygen} $(3)$ and $(4)$.  Point $(3)$
follows from Proposition \ref{proposition:funnyleaves}, the fact
that $\phi$ maps almost every leaf of $\ff$ to itself, and the
fact that leaves of $\ff{\oplus}\sw^s_{\rho(g)}$ are unions of
leaves of $\ff$.
\end{proof}

The majority of the remainder of this paper is devoted to a proof
that the map $\varphi$ constructed above is a small diffeomorphism
with regularity depending on the regularity of $\rho'$. This
suffices to prove Theorem \ref{theorem:main} in the case when the
$\rho$ action lifts to $\cover$.  An additional argument in
subsection \ref{subsection:finiteindex} completes the proof.  With
this in mind, we fix:

\noindent {\bf Notation for the remainder of this paper:} As above
$J$ will be a semisimple Lie group with all simple factors of real
rank at least two and $\G<J$ will be a lattice.  We will fix $G$
to be one of $J$ and $\G$ and also fix a quasi-affine action
$\rho$ of $G$ on $\base$.  Until subsection
\ref{subsection:finiteindex}, we will assume that the $G$ action
$\rho$ lifts to $\cover$.  In subsection
\ref{subsection:finiteindex}, we explain how to remove this
assumption.   In addition, we will fix an integer $k$ and $\rho'$
will always denote a perturbation of $\rho$ which is sufficiently
$C^k$ small so as to be able to apply Theorem
\ref{theorem:candidateconjugacy}, and $\varphi$ will be the
resulting semi-conjugacy.  We allow the possibility that $\rho'$
is $C^l$ for some $l{>}k$, including $l=\infty$, so as to be able
to prove the $C^{\infty}$ case of Theorem \ref{theorem:main}.  We
also fix the maps $\psi$ and $\phi$ from Theorems
\ref{theorem:leafconjugacy} and \ref{theorem:almostconjugacygen},
the projection $p:H{\times}M{\rightarrow}Z{\backslash}H$ and the
map $f:H{\times}M{\rightarrow}Z{\backslash}H$ from Theorem
\ref{theorem:semiconjugacy}.

\section{\bf Continuity along dynamical foliations}
\label{section:smoothing}

In this section, we show that $\varphi$ is a homeomorphism by
showing that it is a homeomorphism when restricted to certain
dynamical foliations. To do this we show that for any
$g{\in}\Phi$, $\varphi$ maps contracting leaves for the action
$\rho(g)$ on $\base$ to contracting leaves for $\rho'(g)$ and
deduce from this that $\varphi$ is a homeomorphism along those
foliations. Throughout this section, all notation is as fixed at
the end of the section \ref{section:conjugacy} or as in subsection
\ref{section:prelim3}. Once we have shown that $\varphi$ is a
homeomorphism, we show, in subsection
\ref{subsection:finiteindex}, how to remove the assumption that
the action $\rho$ lifts to a $G$ action on $\cover$.

\subsection{Equivariance of contracting foliations}
\label{subsection:contractingleaves}

We use the equivariance of $\varphi$ to show that:

\begin{proposition}
\label{proposition:stable} For every $g{\in}\Phi$ and almost every
$x{\in}X$, the map $\varphi$ defined in Theorem
\ref{theorem:candidateconjugacy} maps a set of full measure in the
leaf $\sw_{\rho(g)}^{s}(x)$ into the leaf
${\sw}_{\rho'(g)}^{s}(\varphi(x))$.
\end{proposition}

The proof of this proposition takes up the rest of this
subsection. If $\varphi$ were continuous as well as equivariant,
this would follow easily from standard dynamical arguments. We
begin by introducing some terminology and notation. Fix a finite
set $\Phi$ of elements in $G$ as in Lemma
\ref{lemma:basicdynamics} for the remainder of this section. We
introduce a function which measures the extent to which $\varphi$
does not take stable leaves to stable leaves for a fixed element
$g{\in}\Phi$. We define this function on a
$\base{\times}F^s_{\rho(g)}$, where $F^s_{\rho(g)}<H$ is as in
Proposition \ref{proposition:explicitfoliations}. We denote the
identity in $F^s_{\rho(g)}$ by $e_F$.  We note that $\rho(g)$
induces a contracting automorphism of $F^s_{\rho(g)}$ which we
denote by $\varrho(g)$. We denote the diffeomorphism
$(\rho(g),\varrho(g))$ of $\base{\times}F^s_{\rho(g)}$ by $\tilde
\rho(g)$.  The projection
$\pi:\base{\times}F^s_{\rho(g)}{\rightarrow}\base$ is equivariant
for the $\Za$ action generated by $\tilde \rho(g)$ on
$\base{\times}{F^s_{\rho(g)}}$ and the $\Za$ action generated by
$\rho(g)$ on $\base$. First note that if we take the leaf of
$\ff'{\oplus}\sw^s_{\rho'(g)}$ through a point $x$ in $X$, this is
foliated by stable leaves, each of which intersects the leaf of
$\ff'$ through $x$ in exactly one point. Given a point $y$ on the
leaf of $\ff'{\oplus}\sw^s_{\rho'(g)}$ through $x$ we will look at
it's projection to the leaf $\ff'(x)$ through $x$ defined by this
unique intersection point, call this point $p_{\ff'}(y)$. We
denote the restriction of the fixed Riemannian metric on $\base$
to the foliation $\ff$ by $g_{\ff}$. Note that the bounds on the
derivatives of $\varphi$ along $\ff$ from Theorem
\ref{theorem:candidateconjugacy}$(4)$ and $(5)$ imply that for
almost every $x$ in $\base$, there is a small ball
$B_{\ff}(x,\varepsilon(x))$ such that $\phi$ is a $C^{k-1-\kappa}$
diffeomorphism when restricted to $B_{\ff}(x,\varepsilon(x))$,
where $\kappa$ depends only on the size of the perturbation. For
$x$ in $S$ as defined in Theorem \ref{theorem:candidateconjugacy}
point $(4)$, the number $\varepsilon$ is very close to the number
$r$ specified in that theorem.  For general $x$, the number
$\varepsilon$ depends on the bound from $(5)$ of Theorem
\ref{theorem:candidateconjugacy}. Regardless, whenever
$p_{\ff'}(\varphi(fx))$ is in $\varphi(B_{\ff}(x,\varepsilon(x))$
we define:

$$\delta(x,f)=d_{\ff'}((\varphi|_{B_{\ff}(x, \varepsilon(x))}){\inv}(p_{\ff'}(\varphi(fx))),x).$$

\noindent And let $\delta(x,f)=\infty$ otherwise. While the fact
that $\varepsilon(x)$ is not $\rho(g)$ invariant prevents us from
concluding that $\delta(x,f)$ is $\rho(g)$ invariant, we do have
the following weaker condition on $\delta(f,x)$.   Since $\varphi$
and $\pi$ are $G$ equivariant, $\rho(g)$ is isometric along $\ff$
and $p_{\ff'}$ is $\rho(g)$ equivariant, it follows that if
$\delta(x,f)<\infty$ and $\delta(\rho(g)x,\varrho(g)f)<\infty$,
then $\delta(x,f)=\delta(\rho(g)x,\varrho(g)f)$. It is clear that
$\delta(x,e_F)=0$ for almost every $x{\in}H/{\Lambda}{\times}M$.

We recall some basic facts concerning density points.  For more
discussion of the density points, including a proof of the density
point theorem, see \cite[IV.1]{M2}.  First, we need to specify a
$b$-metric on $\base{\times}F_{\rho(g)}^s$.  Recall that given a
number $b>1$, a $b$-metric on a topological space $Y$ is a map
$d:Y{\times}Y{\rightarrow}\Ra^{\geq0}$ satisfying the usual axioms
of a metric, except that the triangle inequality is replaced by
$d(x,y){\leq}b(d(x,z)+d(z,y))$. Our $b$-metric will be the sum of
the metric induced by our choice of Riemannian metric along
$\base$ with a metric on the fiber analogous to the one introduced
in the proof of \cite[Corollary IV.1.6]{M2}.  Given a ball $B$ in
$F_{\rho(g)}^s$, we define a right invariant $2$-distance function
on $E_{\rho(g)}^s$, by letting:
$$n_B(f_1,f_2)=\max\{n{\in}{\mathbb Z}|(f_1{\inv}f_2){\in}\varrho(g)^n(B)\}$$
$$d_B(f_1,f_2)=2^{-n_B(f_1,f_2)}.$$
This is an addition invariant $2$-distance function by the proof
of \cite[Corollary IV.1.6]{M2}. Given $d_B$ and the distance
$d_{\base}$ induced by adapted metric on $\base$, we define a
$2$-distance function on $\base{\times}E_{\rho(g)}^s$ by letting
$$d((x_1,f_1),(x_2,f_2))=d_{\base}(x_1,x_2)+d_B(f_1,f_2).$$
Whenever discussing density points in $\base{\times}F_{\rho(g)}^s$
(resp. $F_{\rho(g)}^s$ or $\base$) we mean density points with
respect to balls in the metric $d$ (resp. $d_B$ or $d_X$). For
$d_B,d_X$ and $d$, we will denote by $B_d(x,\varepsilon)$ (resp.
$B_{d_B}(x,\varepsilon)$ or $B_{d_{base}}(x,\varepsilon)$)
 the $d$ ball about $x$ of radius $\varepsilon$.

Given a topological space $Y$, a $b$-distance function $d$ on $Y$,
a measure $\mu$ on $Y$, and a measurable set $C{\subset}Y$, we
call a point $y{\in}Y$ a density point of $C$ if
$$\lim_{\varepsilon{\rightarrow}0}\frac{\mu(C{\cap}B_d(x,\varepsilon))}{\mu(B_d(x,\varepsilon))}=1.$$
Note that this conclusion is most reasonable in the case where
$\mu$ is a regular Borel measure which is $d$-finite dimensional
in the sense of \cite[IV.1]{M2}. The generalization of the
classical density point theorem as stated in  \cite[Theorem
IV.1.5]{M2} says that if $Y,d$ are as above and  that if $\mu$ is
$d$-finite dimensional, then the subset of $C$ consisting of
density points of $C$ is of full measure in $C$. We do not give a
more detailed discussion here, since we will use the density point
theorem only through the following consequence, which is a special
case of \cite[Corollary IV.1.6]{M2}.

\begin{proposition}
\label{proposition:densitypoints} Let $F$ be a locally compact,
compactly generated topological group, $\varrho:F{\rightarrow}F$ a
contracting automorphism of $F$ and $C{\subset}V$ a (Haar)
measurable subset.  Then if $e_F$ is a density point of $C$, the
sequence $\{\varrho(g)^{-n}(C)\}_{n{\in}N}$ converges in measure
to $F$.
\end{proposition}

\begin{proof}
In the proof of \cite[Corollary IV.1.6]{M2}, it is shown that the
sets $\{\varrho(g)^{-n}(f{\inv}C)\}_{n{\in}N}$ converge in measure
to $F$ whenever $f$ is a density point of $C$.  This implies the
desired conclusion.
\end{proof}

Let $V_{\beta}$ be the set of points of in
$\base{\times}F_{\rho(d)}^s$ such that $\delta(x,f){\leq}{\beta}$
and let $U_{\beta}$ be the set of points $x$ in $\base$ such that
$\varepsilon(x)<\beta$.  We now show that for any $\beta>0$, the
set $V_{\beta}$ is of full measure in
$U_{\beta}{\times}F_{\rho(d)}^s$.  Note that by the conclusions of
Theorem \ref{theorem:candidateconjugacy}, the set
$U_{\beta}{\rightarrow}\base$ in measure as $\beta{\rightarrow}0$.
The proof of the lemma below is complicated by the fact that we
need to work with points $x$ that are density points both in
$\base$ and along the leaf of $\sw^s_{\rho(g)}(x)$.

\begin{lemma}
\label{lemma:densitypoints} For every $\beta>0$, the set
$V_{\beta}$ is a set of full measure in
$U_{\beta}{\times}F_{\rho(g)}^s$.
\end{lemma}

\begin{proof}
We fix $\beta$ and $U_{\beta}$.  Fix a constant $\eta>0$.  Then by
Theorem \ref{theorem:candidateconjugacy}$(5)$ there is a number
$C_1$ and a set $U_1{\subset}{U_{\beta}}$ with
$\mu(U_1){\geq}(1-\eta)\mu(U_{\beta})$ where for any $x{\in}U_1$,
$\varphi|_{U_1}$ is differentiable and
$\|D\phi(x)|_{\ff}\|{\leq}C_1$.

By Luzin's theorem, we can can choose a set $U_2{\subset}\base$
with $\mu(U_2){\geq}1-\eta$ and a continuous map
$\theta:\base{\rightarrow}\base$ such that $\varphi=\theta$ on
$U_2$.  Let $U_3=U_1{\cap}U_2{\cap}U_{\beta}$ and note that
$\mu(U_3){\geq}(1-2\eta)\mu(U_{\beta})$.

Define the map $\Pi:\base{\times}F_{\rho(g)}^s{\rightarrow}\base$
 by $\Pi(x,f)=fx$.  We define $\tilde
\varphi=\varphi{\circ}\Pi$ and $\tilde \theta=\theta{\circ}\Pi$.
As a result, we have $\tilde \varphi=\tilde \theta$ on
$\Pi{\inv}(U_2)$.

The manifold $\base{\times}F^s_{\rho(g)}$ is equipped with a
product measure $\mu{\times}\nu$ where $\nu$ is Haar measure on
$F^s_{\rho(g)}$. Note that there is no difficulty in applying
Fubini's theorem to this product measure.  For all
$f{\in}F^s_{\rho(g)}$ we $\Pi(\base,f)=\base$, and
$\Pi(\cdot,f)_*\mu=\mu$ and therefore
$\mu(\Pi(\cdot,f){\inv}U_i)=\mu(U_i)$.


By Fubini's theorem and the density point theorem, the set $U_4$
of points which are density points for
$U_3{\cap}(x{\times}F^s_{\rho(g)})$ are of full $\nu$ measure in
$U_3{\cap}(x{\times}F^s_{\rho(g)})$ for almost all $x{\in}\base$.
Applying Fubini again, implies that
$U_4{\cap}({\base{\times}\{f\}})$ is of full $\mu$ measure in
$U_3{\cap}(\base{\times}\{f\})$ for almost every $f$ with
$(x,f){\in}\Pi{\inv}(U_3)$ for some $x{\in}\base$. By changing
basepoint by translating by $f$, we can assume that
$U_5=U_4{\cap}({\base{\times}\{e_F\}})$ is of full measure in
$U_3{\cap}(\base{\times}\{e_F\})$ and that
$\mu(U_5){\geq}(1-2\eta)\mu(U_{\beta})$.

Let $N_R(x)=\{i|\rho(g)^{-i}x{\in}U_5\}$.  The set $N_R(x)$ is
infinite for almost every $x{\in}U_5$ by the Poincar\'{e}
Recurrence Theorem. Given $f{\in}F^s_{\rho(g)}$ let
$N_R(x,f)=N_R(x){\cap}\{j|\varrho(g)^{-j}f{\in}(U_3{\cap}(\{\rho(g)^{-j}(x)\}{\times}F_{\rho(g)}^s)\}$.
Then for $\nu$ almost every $f{\in}F^s_{\rho(g)}$, Proposition
\ref{proposition:densitypoints} implies that $N_R(x,f)$ is
infinite for almost every $x$ in $U_5$.

For $x{\in}U_5, f{\in}F^s_{\rho(g)}, y=fx$ and $n{\in}N_R(x,f)$,
it follows that
$$d(\varphi(\rho(g)^{-n}x),\varphi(\rho(g)^{-n}y))=d(\theta(\rho(g)^{-n}x),\theta(\rho(g)^{-n}y))$$
\noindent since the definition of $N_R(x,f)$ implies that
$\rho(g)^{-n}x$ and $\rho(g)^{-n}y$ are in $U_2$.

The definition of $U_1$,  the fact that $\rho(g)^{-n}(x)$ is in
$U_1$ and compactness of $\base$ imply that there exists a
constant $C$ depending only on the geometry of $\base$ such that
$$d_{\ff}((\varphi|_{B_{\ff}(x,\varepsilon(x))}){\inv}p_{\ff'}(\varphi(y)),x)=$$
$$d_{\ff}((\varphi|_{B_{\ff}(x,\varepsilon(x))}){\inv}(p_{\ff'}(\varphi(\rho(g)^{-n}y))),\rho(g)^{-n}x){\leq}$$
$$CC_1d(\varphi(\rho(g)^{-n}x),\varphi(\rho(g)^{-n}y)=CC_1d(\theta(\rho(g)^{-n}x),\theta(\rho(g)^{-n}y))$$
whenever $x{\in}U_5$ and $n{\in}N_R(x,y)$.  Since $\theta$ is
uniformly continuous and
$d(\rho(g)^{-n}x,\rho(g)^{-n}y){\rightarrow}0$ as
$n{\rightarrow}{\infty}$, by choosing $n{\in}N_R(x,f)$ large
enough, we can guarantee that
$$CC_1d(\theta(\rho(g)^{-n}x),\theta(\rho(g)^{-n}y))<\eta.$$

\noindent Since the choice of $\eta$ is free, this proves the
Lemma.
\end{proof}

\begin{proof}[Proof of Proposition \ref{proposition:stable}]
Take the sequence $V_{\frac{1}{n}}$.  Then
$V=\bigcap_{n=1}^{\infty}V_{\frac{1}{n}}$ is a set of full measure
in $\base{\times}F_{\rho(g)}^s$ and is also a set of full measure
in almost every fiber.  By definition of $V_{\beta}$, for any
$x{\in}U_{\beta}$ such that $V$ is of full measure in the
$F_{\rho(g)}^s$ fiber over $x$, $\varphi$ takes a set of points of
full measure in $\sw_{\rho(g)}^s(x)$ to points in
$\sw_{\rho'(g)}^s(\varphi(x))$.
\end{proof}

\subsection{$\varphi$ is a homeomorphism}
\label{subsection:continuity}

In our setting, $\varphi=\phi{\circ}\psi$ is not a priori a
homeomorphism since $\phi$ is not even a priori continuous.
However, we will show that $\varphi$ is agrees almost everywhere
with a homeomorphism when restricted to the leaves of any of the
foliations $\sw^s_{\rho(g)}$ for $g{\in}\Phi$.  We will then use
this fact to prove that $\varphi$ is in fact a homeomorphism.  We
begin with some definitions.  Recall that if $X$ is a Riemannian
manifold with a foliation $\ff$, there is a natural volume on the
leaves of $\ff$ defined by the restriction of the Riemannian
metric to $T\ff$.  We will call a map of a non-compact space {\em
uniformly  small} if it is uniformly close to the identity on all
compact sets.  Similarly, we say that two homeomorphisms $h,g$ are
uniformly close if $hg{\inv}$ is uniformly small and we say that a
sequence of homeomorphisms $h_n$ on a non-compact {\em converge
uniformly} to a homeomorphism $h$ if the maps $h_nh{\inv}$ is
uniformly small.

\begin{defn}
\label{defn:continuousalongleaves} Given a Riemannian manifold $X$
equipped with a foliation $\ff$ by smooth manifolds and a map
$h:X{\rightarrow}X$, we say that $f$ is:
\begin{enumerate}
\item  {\em essentially continuous} along $\ff$ if for almost
every $x{\in}X$ the restriction of $h$ to $\ff(x)$ agrees almost
everywhere with a continuous map, \item  {\em essentially a
homeomorphism} along $\ff$ if for almost every $x{\in}X$ the
restriction of $h$ to $\ff(x)$ agrees almost everywhere with a
uniformly small homeomorphism and, \item {\em essentially
uniformly continuous} along $\ff$ if it is essentially a
homeomorphism along $\ff$ and  for every sequence
$x_n{\rightarrow}x$ with $h(x_n){\rightarrow}h(x)$, the maps
$h|_{\ff(x_n)}$ agree almost everywhere with maps which converge
uniformly to a homeomorphism $\bar
h:\ff(x){\rightarrow}\ff(h(x))$.
\end{enumerate}
\end{defn}

 The first step in proving continuity of
$\varphi$ is proving:

\begin{proposition}
\label{proposition:continuousalongleaves} For any $g{\in}\Phi$ the
map $\varphi$ is essentially uniformly continuous along
$\sw_{\rho(g)}^c$.
\end{proposition}

Before proving the proposition, we require a lemma that follows
immediately from the definition of $\psi$.  Recall that $\psi$ is
continuous and is covered by a map $\tilde \psi$ such that:

$${\xymatrix{{(\cover,\rho)}\ar[d]_{p}\ar[r]_{\tilde \psi}&{(\cover,\rho')}\ar[dl]^{f}\\
{(\quotient, {\bar \rho})}\\}}$$

\noindent where all maps are right $\Lambda$ equivariant and $p$
and $f$ are left $G$ equivariant.  As an immediate consequence of
this and the fact that $f$ and $p$ are uniformly $C^0$ close and
$\Lambda$ equivariant we have:

\begin{lemma}
\label{lemma:projectionsofstableleaves}
 Let $\sv'$ be a leaf of
$\sw^s_{\rho'(g)}$ and $\tilde \sv'$ the lift of $\sv'$ to
$\cover$. Then $f:\tilde \sv' {\rightarrow}f(\tilde \sv')$ is
homeomorphism onto a leaf of $\sw_{\bar \rho(g)}^s$. Furthermore
if $\sv$ is a leaf of $\sw^s_{\rho(g)}$ with lift $\tilde \sv$ to
$\cover$, and
\begin{enumerate}
\item  $\tilde \sv$ is close to $\tilde \sv'$ and,

\item $f(\tilde \sv')=p(\tilde \sv)$

\end{enumerate}
then $f:\tilde \sv'{\rightarrow}f(\tilde \sv')$ is uniformly close
to $p:\tilde \sv{\rightarrow}p(\tilde \sv)$.
\end{lemma}

Since $\phi$ preserves the foliation of $\base$ which is covered
by fibers of $p$, we have that $\varphi$ commutes locally with the
projections $f$ and $p$, i.e. that for any $U{\subset}\base$, we
have

$${\xymatrix{{U}\ar[d]_{p}\ar[r]_{\varphi}&{\varphi(U)}\ar[dl]^{f}\\
{p(U)}\\}}$$

\noindent where the arrows $p$ and $f$ are defined by viewing $U$
and $\varphi(U)$ as subsets of $\cover$. With this in mind we can
now proceed to prove Proposition
\ref{proposition:continuousalongleaves}.

\begin{proof}[Proof of Proposition
\ref{proposition:continuousalongleaves}] By Proposition
\ref{proposition:stable}, for almost every $x$, we have that
$$\varphi(\sw_{\rho(g)}^s(x)){\subset}\sw_{\rho'(g)}^s(\varphi(x)).$$
By  the commutative diagram above, for almost every $x$ there is a
neighborhood $U$ of $x$ such that
$f(\varphi(\sw_{\rho(g)}^s(x){\cap}U))=p(\sw_{\rho(g)}^s(\varphi(x)){\cap}\varphi(U))$
and therefore
$\varphi(\sw_{\rho(g)}^s(x){\cap}U)=\sw_{\rho'(g)}^s(\varphi(x)){\cap}{\varphi(U)}$.
Furthermore, since by Lemma \ref{lemma:projectionsofstableleaves}
$f|_{\sw_{\rho'(g)}^s(\varphi(x)){\cap}U}$ projects
$\sw_{\rho'(g)}^s(\varphi(x)){\cap}U$ homeomorphically onto
$f(\sw_{\rho'(g)}^s(\varphi(x)){\cap}\varphi(U))$, we can write
$\varphi|_{\tsw_{\rho(g)}^s(x)}$ as

$$p|_{(\sw_{\rho(g)}^s(x){\cap}\varphi(U))}{\circ}f|_{\sw_{\rho'(g)}^s(\varphi(x)){\cap}{\varphi(U)}}{\inv}$$

\noindent which is clearly a homeomorphism.  The fact that
$f|_{(\sw_{\rho'(g)}^s(\varphi(x)){\cap}{\varphi(U)})}$ is $C^0$
close to $p|_{\sw_{\rho(g)}^s(x)}{\inv}{\cap}U$ by Lemma
\ref{lemma:projectionsofstableleaves} implies that $\varphi$ is
essentially a homeomorphism along $\sw_{\rho(g)}^s$.

That $\varphi$ is essentially uniformly continuous along
$\sw^c_{\rho(g)}$ follows from the fact that $f$ and $p$ are $C^0$
and uniformly $C^0$ close on all of $H{\times}M$, are
homeomorphisms when restricted to leaves of $\sw_{\rho(g)}^s$ and
$\sw_{\rho'(g)}^s$ respectively and the fact that the foliations
$\sw_{\rho(g)}^s$ and $\sw_{\rho'(g)}^s$ are continuous.

\end{proof}

At this point we want to conclude that since $\varphi$ is a
essentially a homeomorphism and essentially uniformly continuous
along foliations whose tangent spaces span $T(\base)$ at each
point, $\varphi$ is a homeomorphism.  However it is unclear that
$\varphi$ should agree with a single well-defined global
homeomorphism.  Our proof of this uses that the foliations
involved are smooth, or at least absolutely continuous, in order
to use Fubini's theorem repeatedly.  We first give two general
lemmas from which we will deduce continuity of $\varphi$.  To
avoid technicalities concerning integrability, we will prefer to
work with $1$ dimensional foliations.  For our application we need
only the second statement in the following lemma, but we state and
prove the first statement since it makes the ideas involved
clearer.

\begin{lemma}
\label{lemma:continuityalongfoliations} \begin{enumerate} \item
Let $X$ be an $n$ dimensional compact Riemannian manifold and
$V_1,V_2, {\ldots}V_n$ smooth nowhere vanishing vector fields such
that $TX_x=\oplus_{i=1}^nV_i(x)$ for every $x{\in}X$. Let $\ff_i$
be the foliation tangent to $V_i$, and let  $h:X{\rightarrow}Y$ be
a measurable map which is essentially uniformly continuous along
each $\ff_i$. Then $h$ is a homeomorphism. \item Let $X$ be as
above, let $V_1{\ldots}V_k$ be nowhere vanishing vector fields,
and let $\ff^c$ be a smooth foliation of $X$ by manifolds of
dimension $n-k$ such that
$TX_x=\oplus_{i=1}^kV_i(x){\oplus}T\ff^c(x)$.  Let
$h:X{\rightarrow}Y$ be a measurable map that is essentially
uniformly continuous along each $\ff_i$. Further assume that we
can cover $X$ by foliation charts $U_i$ for $\ff^c$ such that $h$
is a small homeomorphism along most leaves of $\ff^c$ in $U_i$.
Then $h$ is a homeomorphism.
\end{enumerate}
\end{lemma}

\noindent{\bf Remark:} An examination of the proof indicates that
we could make slightly more general assumptions on the vector
fields $V_i$ provided we choose a collection of vector fields
which spans the tangent space of $X$ at every point and such that
the foliations $\ff_i$ are absolutely continuous. We only state
and prove the version needed for our applications to avoid
unnecessary technicalities.

\begin{proof}
We first prove $(1)$ and then explain how to modify the proof to
prove $(2)$.  We work in a chart $U$ that is a foliation chart for
each $\ff_i$ and can in fact assume that $U=\Ra^n$ and that the
foliation $\ff_i$ is given by lines parallel to the line $l_i$
where $x_j=0$ for $j{\neq}i$. We denote the line parallel to $l_i$
passing through the point $y$ by $\tilde l_i(y)$. Let
$W_j=l_1{\times}l_2{\times}{\cdots}{\times}l_j$ and let $\tilde
W_j$ be a $j$ plane parallel to $W_j$ specified by coordinates
$(x_{j+1},{\cdots},x_n)$. We prove by induction that $h$ agrees
almost everywhere with a homeomorphism $\tilde
h_j(x_j+1,{\cdots},x_n)$ along almost every $j$ plane $\tilde
W_j$. By assumption $h$ agrees almost everywhere with a small
homeomorphism along almost every line parallel to $l_1$. Assume
$h$ agrees almost everywhere with a continuous function $\tilde
h_j$ along almost every $j$ plane $\tilde W_j$ parallel to $W_j$.
Then by Fubini, for almost every such $j$ plane $\tilde W_j$, we
have that for almost every $y{\in}\tilde W_j$, the map $h$ agrees
almost everywhere on $\tilde l_{j+1}(y)$ with a small
homeomorphism $\tilde h_y$.  We define a map
$h_{j+1}(x_{j+2},\ldots,x_n)$ on $\tilde W_{j+1}$ by letting
$h_j=\tilde h_y$ where that map agrees almost everywhere with $h$.
This map extends continuously to a homeomorphism on $\tilde
W_{j+1}$ since $h$ is essentially uniformly continuous and
therefore the maps $\tilde h_y$ are uniformly continuous. The map
$h_n$ is a small homeomorphism that agrees almost everywhere on
$U$ with $h$.  That $h_n$ is independent of the chart chosen
follows easily from the definitions.

For $(2)$ we use induction to prove a slightly weaker statement
which still suffices. We re-index our vector fields as $V_{k+1},
{\ldots}, V_n$ and re-index the resulting foliations similarly. We
work in a foliation chart $U$ for $\ff^c$ and $\ff_i$ such that
$\ff^c$ is given by planes $\tilde W_k$ of the form
$\Ra^k{\times}(x_{k+1},{\ldots},x_n)$ and for $k+1{\leq}i{\leq}n$,
each foliation $\ff_i$ is given by lines $l_i$ as above.  To begin
our induction, we use that $h$ agrees with a small homeomorphism
along most of the planes $\tilde W_k$. The induction follows as
before, with the change that at each step we only assume that $h$
agrees almost everywhere with a homeomorphism on most planes
$\tilde W_j$ show that $h$ agrees almost everywhere with a
homeomorphism on most planes $\tilde W_{j+1}$
\end{proof}

We now need a lemma to show that in the setting of our
applications we can apply Lemma
\ref{lemma:continuityalongfoliations}.

\begin{lemma}
\label{lemma:1dimfoliations} Let $X$ be a compact Riemannian
manifold, $\ff$ a smooth foliation of $X$ and $V$ a nowhere
vanishing smooth vector field on $X$ such that $V(x){\in}T\ff_x$
for every $x{\in}X$.  Let $\ff_V$ be the foliation tangent to $V$.
If $h:X{\rightarrow}Y$ is essentially uniformly continuous along
$\ff$, then $h$ is essentially uniformly continuous along $\ff_V$.
\end{lemma}

\begin{proof}  It suffices to work in a chart $U$ which are foliation
charts for both $\ff$ and $\ff_V$.  We can choose coordinates on
such a chart so that $U=\Ra^n$, and the foliation $\ff$ is given
by $\Ra^k{\times}y$ where $y{\in}\Ra^{n-1}$ and the foliation
$\ff_V$ is given by $\Ra{\times}z$ where $z{\in}\Ra^{n-1}$.  We
are assuming that $h$ agrees with a homeomorphism $\tilde h$ of
$\Ra^k{\times}y$ for almost every $y$, and by Fubini's theorem
this implies that $\tilde h$ and $h$ agree almost everywhere on
$\Ra{\times}z$ for almost every $z$.
\end{proof}

We are now prepared to prove continuity of $\varphi$.

\begin{theorem}
\label{theorem:continuity} The map $\varphi$ constructed in the
proof of Theorem \ref{theorem:candidateconjugacy} agrees almost
everywhere with a $C^0$ small homeomorphism.
\end{theorem}

\begin{proof}
We will apply Lemma \ref{lemma:continuityalongfoliations} to
$\varphi$. In doing so, we let $\ff^c$ be $\ff$.  It follows from
Theorem \ref{theorem:candidateconjugacy}$(4)$ that by restricting
to small enough perturbations, we can cover $\base$ by foliation
charts where $\varphi$ is $C^{k-1-\kappa}$ small on most leaves of
$\ff$ for some $\kappa{\leq}1$, and therefore that $\ff$ satisfies
the hypotheses on $\ff^c$ in Lemma
\ref{lemma:continuityalongfoliations}.  We choose elements
$V_i{\in}\fh$ such that each $V_i{\in}\fff^s_{\rho(g)}$ for some
$g{\in}\Phi$ and such that $\oplus V_i \oplus \fz=\fh$.  Then each
$V_i$ defines a smooth non-vanishing vector field $\tilde V_i$ on
$\base$, and $T(\base)_x=\oplus V_i(x) \oplus T\ff$.  By Lemma
\ref{lemma:1dimfoliations} and Proposition
\ref{proposition:continuousalongleaves}, we have that $\varphi$ is
essentially uniformly continuous along the foliation $\ff_i$
tangent to $V_i$ for each $i$.  Therefore, we can apply Lemma
\ref{lemma:continuityalongfoliations} to $\varphi$ which implies
that $\varphi$ is a small homeomorpism.
\end{proof}

\subsection{Additional arguments in the case of discrete groups}
\label{subsection:finiteindex}

In the case of $\G$ actions, we have been assuming that the
unperturbed $\G$ action lifts to the cover $\cover$.  As remarked
above, this is always true on a finite index subgroup $\G'$ of
$\G$ which depends only on $\rho$.  We have constructed a
continuous $C^0$ small conjugacy for the $\G'$ actions with
additional regularity along $\ff$ and we now explain how to
replace this with a $C^0$ small conjugacy for the $\G$ actions
with the same additional regularity along $\ff$. The passage to
$\G'$ is required in the proof of Theorem
\ref{theorem:almostconjugacygen}. In that proof, when we conjugate
the $\rho'$ action by $\psi$  we only know that the $\G'$ action
defined by $\psi{\inv}{\circ}{\rho'}{\circ}\psi$ preserves $\ff$
and therefore is in small neighborhood of $\rho$ in
$\Hom(\G',\Diff^k_{\nu}(\base,\ff)$. Neither of these facts is
clear for the full $\G$ action.  In fact it suffices to show that
 $\psi{\inv}{\circ}{\rho'}{\circ}\psi(\G)$ preserves $\ff$, since
 closeness to $\rho$ in $\Hom(\G,\Diff^k_{\nu}(\base,\ff)$ then
 follows from the definition of $\psi$.  Therefore the remainder
 of this subsection is dedicated to a proof that $\psi{\inv}{\circ}{\rho'}{\circ}\psi(\G)$
 preserves $\ff$.  Without loss of generality, we assume $\G'$ is normal in $\G$.

Given two (closed) subsets $A,B$ of a metric space $(X,d)$, we let
$d_S(A,B)=\inf_{a{\in}A,b{\in}B}d(a,b)$.

\begin{defn}
\label{defn:leafwiseexpansive} Let a group $D$ act on a manifold
$X$ preserving a foliation $\ff$.  We call the action {\em
$c$-leafwise expansive} if there exists a constant $c$, such that
if $\fL$ and $\fL'$ are distinct leaves of $\ff$ and there is
$f{\in}F$ such that $d_S(\rho(g)\fL,\rho(g)(\fL'))>c$.
\end{defn}

\noindent Note that many foliations, e.g. any foliation with a
dense leaf, do not admit leafwise expansive actions. We will be
applying Definition \ref{defn:leafwiseexpansive} to the lift of
$\rho$ to $H{\times}M$, which is $c$-leafwise expansive by
Corollary \ref{corollary:separatingpoints}.

\begin{lemma}
\label{lemma:commutingwithexpansive} Let $c>0$ and $\rho$ be a
$c$-leafwise expansive action of a group $D$ on a foliated metric
space $(X,d,\ff)$.  Let $h$ be a homeomorphism of $X$ such that:
\begin{enumerate}
\item $d(h(x),x)<c$ for all $x{\in}X$ \item
$h(\rho(g)\fL)=\rho(g)(h(\fL))$ for any leaf $\fL$ of $\ff$ and
any $d{\in}D$, i.e. $h$ and $\rho$ commute as actions on leafs of
$\ff$.
\end{enumerate}
then $h(\fL)=\fL$ for every leaf $\fL$ of $\ff$.
\end{lemma}

\begin{proof}
Assume $h(\fL){\neq}\fL$.  Then there is a point $x{\in}{\fL}$
with $h(x){\notin}{\fL}$.  By our assumptions, there exists
$g{\in}D$ such that $d_S(\rho(g)(\ff_{h(x)}),\rho(g)\ff_x)>c$. But
then $d_S(h(\rho(g)(\ff_x)),(\rho(g)(\ff_x)))>c$ which contradicts
$(1)$ above.
\end{proof}

We define a subgroup $\Homeo(\base,\ff)$ of $\Homeo(\base)$ which
consists of all homeomorphisms which map each leaf of $\ff$ to
itself.

\begin{proposition}
\label{proposition:centralizer} Given a quasi-affine action $\rho$
of $\G'$ on $\base$ which lifts to $\cover$, any small enough
homeomorphism in the centralizer in $\Homeo(\base)$ of $\rho(\G')$
is an element of $\Homeo(\base,\ff)$.
\end{proposition}

\begin{proof}
If $f$ is a small homeomorphism commuting with $\rho(\G')$, the
there is a unique lift $\tilde f$ of $f$ to $\cover$ such that
$\tilde f$ is small as a homeomorphism of $\cover$.  Since $\tilde
f$ is small, for any small enough $\g{\in}\G'$ we have that
$[\tilde f,\tilde \rho(\g)]$ is a small homeomorphism of $\cover$
covering the identity on $\base$ and so $\tilde f$ and $\tilde
\rho(\g)$ commute. Since $\G'$ is finitely generated this implies
that $\tilde f$ commutes with $\G'$ on $\cover$.  Let $\tilde \ff$
be the lift to $\cover$ of the foliation $\ff$. Since Corollary
\ref{corollary:separatingpoints} implies that the $\tilde
\rho'(\G)$ action on $\cover$ is leafwise expansive, Lemma
\ref{lemma:commutingwithexpansive} implies that $\tilde f$ maps
each leaf of $\tilde \ff$ to itself. This then implies that $f$
maps each leaf of $\ff$ to itself.
\end{proof}

We need one more purely algebraic lemma. For simplicity, we will
write $\rho''(\g)$ for $\varphi{\inv}{\circ}\rho'(\g)\circ\varphi$
and state the lemma only in the form needed for our applications.
The lemma is true for any pair of homomorphisms from a group $D$
to a group $H$ which agree on a normal subgroup in $D$.

\begin{lemma}
\label{lemma:commuting} If $\g_0$ is in $\G$, then the
diffeomorphism $\rho(\g_0){\circ}\rho''(\g_0){\inv}$ commutes with
$\rho(\g)$ for all $\g{\in}\G'$.
\end{lemma}

\begin{proof}
For any $g{\in}\G'$ and $\g{\notin}\G$, we have
$\g_0\g\g_0{\inv}{\in}\G'$ which implies that
$\rho''(\g_0\g\g_0{\inv})=\rho(\g_0\g\g_0{\inv})$.  Expanding
gives:
$$\rho''(\g_0)\rho(\g)\rho''(\g_0{\inv})=\rho(\g_0)\rho(\g)\rho(\g_0{\inv})$$
which can be rearranged as
$$\big(\rho(\g_0){\inv}\rho''(\g_0)\big)\rho(\g)=\rho(\g)\big(\rho(\g_0{\inv})\rho''(\g_0)\big)$$
proving the lemma.
\end{proof}

We choose a set of coset representatives $\g_1,\ldots,\g_j$ for
$\G/\G'$ and assume that $\rho'$ is close enough to $\rho$ so that
$\rho(\g_i)\rho''(\g_i){\inv}$ is sufficiently $C^{0}$ small so
that Lemma \ref{lemma:commuting} and Proposition
\ref{proposition:centralizer} imply that
$\rho(\g_i)\rho''(\g_i){\inv}$ is close to the identity in
$\Homeo(\base,\ff)$.  This implies that
$\rho''(\g_i)\fL=\rho(\g_i)\fL$ or that
$\rho'(\g)\circ\varphi(\fL)=\varphi{\circ}\rho(g_i)(\fL)$ for all
$\g{\in}\G$.

We know that $\varphi=\phi{\circ}\psi$ where $\psi$ is the
homeomorphism constructed in Theorem \ref{theorem:leafconjugacy}.
This implies that $\phi$ is also a homeomorphism, which, by
construction is in $\Homeo(\base,\ff)$.  Combined with the
conclusion of the last paragraph, this implies that
$\rho'(\g){\circ}\psi(\fL)=\psi{\circ}\rho(\g_i)(\fL)$ for all
$\g{\in}\G$.  This suffices to allow us to apply Theorem
\ref{theorem:almostconjugacygen} to the entire $\G$ action in the
proof of Theorem \ref{theorem:candidateconjugacy}, rather than
just to the $\G'$ action.  This constructs a map $\phi$ such that
$\varphi=\phi{\circ}{\psi}$ is $\G$ equivariant and $\phi$
satisfies all the conclusions stated in Theorem
\ref{theorem:candidateconjugacy}.

\noindent {\bf Remarks:}\begin{enumerate}

\item We can now re-apply the arguments of subsections
\ref{subsection:contractingleaves} and \ref{subsection:continuity}
to show that $\varphi$ is a $C^0$ small homeomorphism.

\item It is not clear that the $\varphi$ constructed from the
$\G'$ action is actually equivariant for $\G$.  In applying
Theorem \ref{theorem:almostconjugacygen} to the $\G$ action, we
may be finding a different conjugacy.

\item Due to the arguments of this subsection, for the remainder
of this paper, we no longer assume that $\rho$ lifts to an action
of $G$ on $\cover$.
\end{enumerate}

\section{Smoothness along dynamical foliations}
\label{section:katokspatzier}

In this section, we adapt the method of Katok-Spatzier to show
that $\varphi$ is differentiable along certain special expanding
and contracting foliations by constructing transitive $C^k$ group
actions along those foliations.  All notations are as in the
previous section.

\subsection{Some other important dynamical foliations}
\label{subsection:weighffoliations}

In this subsection we define some additional important foliations
related to the group actions $\rho$ and $\rho'$.  These foliations
are the ones to which we will apply the method of Katok-Spatzier,
building transitive smooth actions of Lie groups, along the
leaves, that are intertwined by $\varphi$. First we define the
relevant foliations in the unperturbed setting. The exposition
here is similar to the exposition in section 5.1 of \cite{MQ}.

Recall that $G=J$ or $\Gamma$.  Let $S$ be a maximal $\mathbb R$
split torus in $J$ and $T$ be a maximal torus containing $S$. The
foliations we are now defining depend on the choice of $T$ and in
the case of $\G$ actions, we will make restrictions on that choice
below.

Recall from Definition \ref{definition:affine} that $\rho$ is a
skew product action on $\base$. More precisely, the action on
$\base$ is defined by an action on $H/{\Lambda}$ and a cocycle
$\iota:G{\times}H/{\Lambda}{\rightarrow}K$ over that action where
$K<\Isom(M)$.  Recall from Theorems
\ref{theorem:describingactionsG} and
\ref{theorem:describingactionsGamma} that, possibly after passing
to a subgroup of finite index when $G=\G$, the action on
$H/{\Lambda}$ is defined by a homomorphism $\pi:G{\rightarrow}L$
where $L=\Aut(H)^0{\ltimes}H$ is an algebraic group. Note that $H$
is normal in $L$, so $\fh$ is invariant under $\Ad{L}$. We have an
invariant splitting of the tangent bundle
$T(\base)=(H/{\Lambda}{\times}{\fh}){\times}TM$ and all elements
of $G$ are isometries along $TM$.  The derivative cocycle leaves
$\fh$ and $TM$ invariant, and, again after passing to a subgroup
of finite index if $G=\G$ the restriction to $\fh$ is given by the
representation $\sigma=Ad_L|_{\fh}{\circ}\pi_0$ of $G$ on $\fh$.
From now on when describing the action and the derivative cocycle,
we assume that if the acting group is $\G$ we have passed to a
finite index subgroup for which this description holds.  We recall
that $\pi_0=\pi_0^E\pi_0^K$ where $\pi_0^E$ is (the restriction
of) a representation of $J$, $\pi_0^K$ has bounded image, and the
images of $\pi_0^K$ and $\pi_0^E$ commute. Therefore we can write
$\sigma=\sigma^E\sigma^K$ where $\sigma^E$ is (the restriction of)
a representation of $J$, $\sigma^K$ has bounded image and the
images of $\sigma^K$ and $\sigma^E$ commute.

For $g{\in}T{\cap}G$, define the {\em Lyapunov exponents} of
$\sigma(g)$ as the log's of the absolute values of the eigenvalues
of $\sigma(g)$.  We obtain homomorphisms
$\chi:T{\cap}G{\rightarrow}{\mathbb R}$ which extend to
homomorphisms $\chi:T{\rightarrow}{\mathbb R}$.  The $\chi$ are
exactly the absolute values of the weights of the representation
$\sigma^E$ for the torus $T$, and we will refer to them as {\em
generalized weights}. There is a decomposition of $\fh$ into
generalized weight spaces $E_\chi$, $\fh={\oplus}_{\chi}E_{\chi}$.
Corresponding to this there is a decomposition of the tangent
bundle to $\base$ into invariant subbundles for the derivative
action,
$T(\base)=\big((H/{\Lambda}{\times}({\oplus}_{\chi}E_{\chi})\big){\times}TM$.
We call $\base{\times}E_{\chi}$ a {\em Lyapunov distribution} for
the $G{\cap}T$ action defined by $\rho$.

The set $\Omega$ of all generalized weights for $(\sigma,T)$ can
be decomposed into disjoint subsets $[\chi]$ such that
$\chi'{\in}[\chi]$ if and only if $\chi'=t\chi$ for some positive
real number $t$. We fix a set $\bar \Omega{\subset}\Omega$ of
representatives for the subsets $[\chi]$.  If $\chi=0$
identically, we call $\base{\times}E_{0,T}$ the {\em central
distribution} for the action $G{\cap}T$.  It is integrable, and we
denote by ${\sw}_{0,T}$ the corresponding foliation.  It is clear
that ${\cap}_TE_{0,T}=E_0$ defines an integrable distribution and
that $H/{\Lambda}{\times}E_0{\oplus}TM=E^c_{\rho(G)}$. Other
Lyapunov distributions may or may not be integrable, but
$\base{\times}E_{[\chi]}={\oplus}_{\lambda{\in}{[\chi]}}\base{\times}E_{\lambda}$
is always integrable. We denote the integral foliation on $\base$
by ${\sw}_{[\chi]}$ suppressing the dependence on $\rho$ and $T$
in the notation.

We now describe the choice of $T$ for the case $G=\G$. The reader
should note that if we choose $T$ so that $T{\cap}\G=e$, then all
of $\fh$ is in $E_0$.  To obtain a more useful set of Lyapunov
distributions, we need the following theorem which we derive from
results of Prasad and Rapinchuk from \cite{PRp}.  Without loss of
generality, we assume that $\Gamma$ is the direct product of a
finite number of irreducible lattices $\Gamma_i$.  We write
$J=\prod_IJ_i$ where $\G_i<J_i$ is irreducible, and for any
maximal torus $T<J$, we can write $T=\prod_IT_i$ where $T_i<J_i$
is a maximal torus.

\begin{theorem}
\label{theorem:prasadrapinchuk} Let $J$ and $\Gamma$ be as above.
Then there is a maximal torus $T$ in $J$ such that:
\begin{enumerate}
\item $T$ contains a maximal $\Ra$-split torus $S$ \item
$\Sigma=\Gamma{\cap}T$ is cocompact in $T$, \item For any $T_i$
there is no proper algebraic torus $T_i'<T_i$ such that
$T_i'{\cap}\Gamma$ is a lattice in $T_i'$.
\end{enumerate}
Furthermore if  $\mu:J{\rightarrow}GL_n(\mathbb R)$ is any linear
representation of $J$ and $\sigma{\in}\Sigma$ projects to an
infinite order element in each $\Gamma_i$, then $\xi(\sigma)$ is
not a root of unity for any nontrivial weight $\xi$ of $\mu$.
\item
\end{theorem}

\begin{proof} It suffices to prove the theorem for
$\G_i<J_i$.  The first two assertions follow from \cite[Theorem
1]{PRp} and the assertion immediately preceding the proof of that
theorem.  The third assertion is an immediate consequence of two
facts.  First $\G_i$ is arithmetic, and therefore $\G_i<\Ga(k)$
for some field $k$.  Combined with  \cite[Proposition
$1(ii))$]{PRp} this implies that any infinite order element of
$\G_i$ generates a Zariski dense subgroup of $T_i$.  The last
statement follows from \cite[Proposition $1(iii)$]{PRp} and the
fact that $\mu$ necessarily agrees with a homomorphism of $\G_i$
on subgroup of finite index and so must be defined over a finite
extension of $k$.
\end{proof}

\noindent From now on we assume that we have picked $T$ satisfying
the conclusion of Theorem \ref{theorem:prasadrapinchuk}.  We do
not use all the properties of $T$ here, but will need them in
subsection \ref{subsection:ergodiccomponents}.

  The
following lemma is analogous to \cite[Lemma 5.2]{MQ}.

\begin{lemma}
\label{lemma:tangentdecomposition} Fix a maximal torus $T$ as
above.  Let $E(T)$ be the sum of $E_{[\chi]}$ for all non-trivial
weights $\chi$ for $(\sigma,T)$.  Then there exists a finite
subset $\Psi{\subset}G$ such that:
$$T(\base)=TM{\times}\big({\sum}_{g{\in}\Psi}D{\rho}(g)(H/{\Lambda}{\times}E(T))\big){\times}E_0.$$
\end{lemma}

\begin{proof}
It suffices to prove that
$$T(H/{\Lambda})=\big({\sum}_{g{\in}\Psi}D{\rho}(g)(H/{\Lambda}{\times}E(T))\big){\oplus}E_0).$$
We know that the derivative of $\rho$ is given by linear
representation $\sigma:G{\rightarrow}\Ad(L)$.  The structure of
$\sigma$ implies that $\fh=\fh^E{\oplus}\fh^K$ where $\sigma^E$
(resp. $\sigma^K$) is trivial on $\fh^K$ (resp. $\fh^E$). From the
definitions, it suffices to see that
$$\fh^E=\sum_{g{\in}\Psi}D{\rho}(g)(H/{\Lambda}{\times}E(T)).$$
\noindent The inclusion of the right hand in the left hand side is
clear.  If the left hand side contains a subspace $V$ not
contained in the right hand side, then we have that $V$ is in the
zero weight space for $\sigma^E|_{G{\cap}T}$ with respect to our
choice of maximal torus $T$ for $J$.  For $G=J$, this is only
possible if the representation of $\sigma^E|V$ is trivial,
contradicting our assumptions.  For $G=\G$, the contradiction
follows since we have chosen a $T$ as in Theorem
\ref{theorem:prasadrapinchuk}.
\end{proof}

Fix a non-trivial generalized weight $\chi_0$ for $(\sigma^E,T)$.
Then there exists $g_0{\in}T{\cap}D$ such that $\chi_0(g_0)<0$. It
follows that for all $\chi'{\in}[\chi_0]$, $\chi'(g_0)<0$.  For
every $g{\in}T{\cap}G$ with $\chi_0(g)<0$, note that
${\oplus}_{\chi(g)<0}E_{\chi}$ is the stable distribution
$E^s_{\rho(a)}$.  This distribution is tangent to the foliation
${\sw}_{\rho(a)}^s$ . It is clear that $E_{\rho(a)}^s$ contains
$E_{[\chi_0]}$.  We call an element of $T$ {\it regular} if for
all non-trivial weights $\chi$ for $(\sigma^E,T)$,
$\chi(a){\neq}0$. Combined with the usual descending chain
arguments, this yields the following lemma.

\begin{lemma}
\label{lemma:weightdistributions} Let $\chi_0$ be a non-trivial
weight for $(\sigma^E,T)$.  Then
$$E_{[\chi_0]}={\bigcap}E^s_{\rho(a)}$$
where the intersection is taken over all regular $a$ with
$\chi_0(a)<0$.  Furthermore, there exist regular elements
$a_1,{\ldots}a_q{\in}T{\cap}G$ with $\chi_0(a_i)<0$ such that we
can take the intersection just over $E_{\rho(a_i)}^s$.
\end{lemma}

We now define a finite collection of foliations and distributions
which we will use below.  Fix a maximal torus $T$ and a set $\bar
\Omega{\subset}\Omega$ as above. Also fix a collection of elements
$a_1,\ldots,a_q{\in}T$ as in Lemma
\ref{lemma:weightdistributions}. Given $g{\in}G$ let
$E^g_{[\chi]}=\rho(g)E_{[\chi]},
\sw^g_{[\chi]}=\rho(g)\sw_{[\chi]},
E^{g,s}_{\rho(a_i)}=\rho(g)E^s_{\rho(a_i)}$ and
$\sw^{g,s}_{\rho(a_i)}=\rho(g)\sw^s_{\rho(a_i)}$. We can also
define $E^{g,s}_{\rho'(a_i)}=\rho'(g)E^s_{\rho'(a_i)}$ and
$\sw^{g,s}_{\rho'(a_i)}=\rho'(g)\sw^s_{\rho'(a_i)}$. We will show
that $\varphi$ is smooth along each $\sw^g_{[\chi]}$ for
$g{\in}\Psi$ and $\chi{\in}\bar \Delta$.  To do this, we first
need to identify the image of $\sw^g_{[\chi]}$ under $\varphi$.

\begin{proposition}
\label{proposition:specialfoliations} For every $m>0$, if $\rho'$
is sufficiently $C^1$ close to $\rho$, for every $x$ the
intersection $\tsw^g_{[\chi]}={\bigcap}\sw^{g,s}_{\rho'(a_i)}$ is
a $C^k$ submanifold tangent to the distribution ${\tilde
E}_{[\chi_0]}={\bigcap}{\tilde E}^{s,g}_{a_i}$.  Furthermore,
$\varphi(\sw_{[\chi]}^g(x))=\tsw^g_{[\chi]}(\varphi(x))$ for every
$x{\in}X$.
\end{proposition}

\begin{proof}
It suffices to consider the case where $g$ is the identity, since
other cases follow by translation.  We will show that the
intersection is $s$-transverse and then apply Lemma
\ref{lemma:stransverse}$(1)$. Since dimension of intersection of
the distributions $E^s_{\rho'(a_i)}$ can only be smaller than for
the corresponding intersection of $E^s_{\rho(a_i)}$ , it suffices
to show that dimensions of the intersections of the foliations
does not decrease. This is immediate from the fact that $f$
projects $\sw^{s}_{\rho'(a_i)}$ homeomorphically onto $\sw^s_{\bar
\rho(a_i)}$ by Lemma \ref{lemma:projectionsofstableleaves} and $p$
projects $\sw^{s}_{\rho(a_i)}$ homeomorphically onto $\sw^s_{\bar
\rho(a_i)}$ by definition.  So dimensions of all intersections of
foliations are equal in the perturbed and unperturbed cases and
therefore the intersection ${\bigcap}\sw^{g,s}_{\rho'(a_i)}$ is
$s$-transverse.

The final claim follows from the fact that
$\varphi(\sw_{\rho(a_i)}^{g,s}(x))=\sw^{g,s}_{\rho'(a_i)}(\varphi(x))$
which is true for each $a_i$ for almost every $x$ by Proposition
\ref{proposition:stable}, and therefore also true for all $x$ by
Theorem \ref{theorem:continuity} and Lemma
\ref{proposition:continuousalongleaves}.  It is also possible to
prove that
$\varphi(\sw_{\rho(a_i)}^{g,s}(x))=\sw^{g,s}_{\rho'(a_i)}(\varphi(x))$
for all $x$ directly by a dynamical argument.
\end{proof}

The rest of this section describes a variant of the method of
Katok-Spatzier which we use to show smoothness of $\varphi$ along
the foliations ${\sw}^g_{[\chi]}$ for $g{\in}\Psi$ and
$\chi{\in}\bar \Delta$.  The outline here is close to that of
\cite{KS} or \cite{MQ} , but there are two additional
difficulties. First, we need to have estimates on the $C^k$ size
of various maps in both the normal form theory of Guysinsky-Katok
\cite{GK,Gu} and in the work of Montgomery-Zippin \cite{MZ}. In
both cases, these estimates follow from examination of the
existing proofs, as is explained below. Secondly, we will need to
show that ergodic components consist of entire leaves of
$\sw^g_{[\chi]}$ for a more general class of actions than those
considered in \cite{MQ}.

The precise statement we prove is:

\begin{theorem}
\label{theorem:smoothalongleaves} We use the notation introduced
before.  Let $n$ be the dimension of ${\sw}^g_{[\chi]}$. Let
$\theta:D^n{\times}D^{m-n}{\rightarrow}X$ be a smooth foliation
chart for ${\sw}^g_{[\chi]}$.  Then there is a number $k_1$
depending only on $\rho$ such that for all $k{\geq}k_1$ and all
$\rho'$ sufficiently $C^k$ close to $\rho$:
\begin{enumerate}
\item  the map
$\varphi:{\sw}_{[\chi]}(x){\rightarrow}{\tsw}_{[\chi]}(x)$ is a
$C^k$ injective immersion \item the map
$Q:D^{m-n}{\rightarrow}\Emb^k(D^n,X)$ given by
$\phi{\circ}\theta(\cdot,y)$ is continuous and $C^0$ close to the
map $Q_0$ induced by the identity on $\base$.
\end{enumerate}
Furthermore if $\rho'$ is $C^l$ for some $l{\geq}k$ then
\begin{enumerate}
\item  the map
$\varphi:{\sw}_{[\chi]}(x){\rightarrow}{\tsw}_{[\chi]}(x)$ is a
$C^l$ injective immersion \item the map
$Q:D^{m-n}{\rightarrow}\Emb^l(D^n,X)$ given by
$\phi{\circ}\theta(\cdot,y)$ is continuous.
\end{enumerate}
\end{theorem}

\subsection{\bf Theory of non-stationary normal forms}
\label{subsection:normalforms} Before giving the construction of
the groups acting transitively on foliations, we outline the
theory of non-stationary normal forms that will be used to show
smoothness of the group actions on the leaves for the perturbed
action. The theorems we use are due to Guysinsky and Katok and the
references are \cite{Gu,GK}. Some of our definitions are slightly
different from those of \cite{Gu,GK}.

Consider a continuous extension $\sF$ of a homeomorphism $f$ of a
compact connected metric space $X$ to a vector bundle $V$ over $X$
which is smooth along the fibers and preserves the zero section.
Let $F=D{\sF}_0$ where the derivative is taken at the zero section
in the fiber direction.  Fix a continuous family of Riemannian
metrics on the fiber of $V$.  Consider the induced operator
$\sF^*$ on the Banach space of continuous sections of $V$ endowed
with the uniform norm, given by $\sF^*v(x)=F(v(f{\inv}(x)))$. For
$i=1,{\ldots}l$, let $\Delta_i=[\lambda_i,\mu_i]$ be a finite set
of disjoint intervals on the negative half line with
$\lambda_{i+1}>\mu_i$.  Assume that $V$ splits as a sum of
subbundles $V=V_1{\oplus}{\cdots}{\oplus}V_k$ such that the
spectrum of $\sF^*$ on the space of sections of $V_i$ is contained
in the annulus with inner radius $\exp(\lambda_i)$ and outer
radius $\exp(\mu_i)$. If $\mu_l<0$, then the map $\sF$ is a
contraction with respect to the continuous family of Riemannian
metrics chosen above.

{\bf Remark:}When $\sF$ is a contraction it also makes sense to
consider $\sF$ which is only defined in a neighborhood of the zero
section in $V$.  Theorem \ref{theorem:normalform} below holds in
this generality, and with some care a version of Theorem
\ref{theorem:centralizernormalform} can be stated in this context
as well.

We say that $\sF$ has {\em narrow band spectrum} if
$\mu_i+\mu_l<\lambda_i$ for all $i=1,{\ldots},l$.

We call two extensions $C^k$ conjugate if there exists a
continuous family of $C^k$ diffeomorphisms of the fibers $V(x)$,
preserving the origin which transforms one extension into the
other.  The following two theorems on normal forms and
centralizers are from \cite{Gu,GK}.  We remark that to avoid
unnecessary definitions we did not state the theorems in their
full generality, but these are sufficient for our applications.

\begin{theorem}
\label{theorem:normalform} Let $f$ be a homeomorphism of a metric
space $X$ and suppose that $\sF$ is $C^l$ extension of $f$ which
is a contraction, that the linear extension $D{\sF}_0$ has narrow
band spectrum determined by the vectors
$\lambda=(\lambda_1,{\ldots},\lambda_l)$ and
$\mu=(\mu_1,{\ldots},\mu_l)$, and that $\sF$ is $C^k$ close to
$D{\sF}_0$ in a neighborhood of the zero section. There exists a
constant $k_1=k_1(\lambda,\mu)$ such that if $k{\geq}k_1$ there
exist
\begin{enumerate}
\item a finite dimensional Lie group $G_{\lambda,\mu}$ which is a
subset of all polynomial maps from $\Ra^m$ to $\Ra^m$ of degree
less than or equal to $d$ for some $d<\infty$,

\item an extension $\tilde \sF$ such that for every $x{\in}X$, the
map
$$\tilde \sF|_{V(x)}:V(x){\rightarrow}V(f(x))$$
is an element of $G_{\lambda,\mu}$;

\item a $C^l$ conjugacy $H$ between $\tilde \sF$ and $\sF$ which
is  $C^k$ small.
\end{enumerate}
\end{theorem}

\begin{theorem}
\label{theorem:centralizernormalform} Suppose $g$ is a
homeomorphism of the space $X$ commuting with $f$ and $\sg$ is an
extension of $g$ by $C^l$ diffeomorphisms of the fibers commuting
with the extensions $\sF$ and that $\sF$ satisfies all of the
hypotheses of Theorem \ref{theorem:normalform} and $k{\geq}k_1$.
Then $H$ conjugates $\sg$ to a map of the same form, i.e. one
where:
$$\tilde \sg|_{V(x)}:V(x){\rightarrow}V(f(x))$$
is a polynomial of degree at most $d$ and is in fact in the group
$G_{\lambda,\mu}$ from Theorem \ref{theorem:normalform}.
\end{theorem}

\begin{proof}
The only statement which is not justified explicitly in the proofs
of \cite{GK,Gu} is the bound on the size of $H$.  Recall that $H$
is constructed in two steps.  First one constructs a conjugacy
between $\sF$ and an extension $\sF'$ of $f$ such that $\sF'$ is
$C^k$ tangent to ${\tilde \sF}$ at the zero section.  In this step
one proceeds by solving an iterative equation for the conjugacy,
see \cite[Proof of Theorem $1$,Step $1$, page 851]{Gu}.  It is
clear from the formula that if $\sF$ and $D{\sF}_0$ are close,
then this conjugacy is small and $\sF'$ is also $C^k$ close to
$D{\sF}_0$. (We note that in \cite{GK}, this step is broken into
two steps, first finding the Taylor series of the conjugacy at the
zero section, and then proving that one can find a conjugacy with
this Taylor series.) In the second step, one constructs an action
$\bar F$ on a set of local changes of coordinates, and applies a
contraction mapping argument to find the conjugacy $H$ between
$\sF'$ and and $\tilde {\sF}$. If $\sF'=D{\sF}_0$ it is clear from
the construction that this contraction $\bar F$ has unique fixed
point the identity map, and that if $\sF'$ is $C^k$ close to
$D{\sF}_0$ then this unique fixed point of $\bar F$ will be $C^k$
close to the identity.
\end{proof}

\noindent {\bf Remark:} The number $k_1$ is explicitly computable
in terms of the spectrum of the contraction $\sF$, see
\cite{Gu,GK} for details.  The computation yields that Theorem
\ref{theorem:normalform} and \ref{theorem:centralizernormalform}
are true for
$$k_1{\geq}|\frac{\lambda_1\lambda_2{\cdots}\lambda_{l-1}}{\mu_2\mu_3{\cdots}\mu_l}|.$$
For some special choices of $\Delta_i$ it is possible to achieve
much lower values of $k_1$.

\subsection{\bf Smoothness along contracting foliations}
\label{section:contractingsmoothness}

In this subsection we retain the number $k_1$ is as in the last
subsection and $\varphi$ is the map constructed in Theorem
\ref{theorem:candidateconjugacy}, which we know to be a $C^0$
small homeomorphism and a conjugacy between the unperturbed and
perturbed actions. We first show that to prove Theorem
\ref{theorem:smoothalongleaves}, it suffices to verify the
following lemma.

\begin{lemma}
\label{lemma:groupsonleaves} For $k{\geq}k_1$ there is a connected
Lie group $\Delta$ and, for each $x{\in}{\cover}$, there is an
open set $U_0{\subset}{\cover}$ which contains $x$ and is the
union of leaves of $\sw_{[\chi]}$, such that:

\begin{enumerate}
\item[S1] there is a locally free $C^{\infty}$ action
$d:\Delta{\times}U_0{\rightarrow}U_0$ such that
$\delta{\sw}_{[\chi]}(y)={\sw}_{[\chi]}(y)$ and $\Delta$ acts
transitively on $\sw_{[\chi]}(y)$ for all $y{\in}U_0$;

\item[S2] the set $\varphi(U_0)={U_0}'$ is the union of leaves of
$\tsw_{[\chi]}$ and there is a locally free $C^0$ action
$d':\Delta{\times}{U_0}'{\rightarrow}{U_0}'$ such that
$\delta{\tsw}_{[\chi]}(z)={\tsw}_{[\chi]}(z)$ for all
$z{\in}{U_0}'$.

\item[S3] for all $\delta{\in}\Delta$, the map
$d'(\delta):{U_0}'{\rightarrow}{U_0}'$ is $C^l$ when restricted to
every leaf of $\tsw_{[\chi]}$, and all partial derivatives along
the leaf are globally continuous. Furthermore, the $k$-jet of
$d'(\delta)$ along leaves of $\tsw_{[\chi]}$ tends to the $k$-jet
of $d(\delta)$ along $\sw_{[\chi]}$ as $\rho'{\rightarrow}{\rho}$.

\item[S4] $\varphi$ is a $\Delta$ equivariant map from ${U_0}$ to
${U_0}'$.
\end{enumerate}
\end{lemma}

We temporarily defer the proof of Lemma \ref{lemma:groupsonleaves}
and first show how it implies Theorem
\ref{theorem:smoothalongleaves}.  We state a variant of the
results of Bochner and Montgomery which we will use in the proof.

\begin{theorem}[Bochner and Montgomery]
\label{theorem:bochnermontgomery} Let $\kappa$ be a continuous
action of a Lie group $Q$ on a manifold $N_1{\times}N_2$ such that
the action is trivial in the second factor. If for each
$q{\in}Q$,and $n_2{\in}N_2$, the map
$\kappa(q):N_1{\times}n_2{\rightarrow}N_1{\times}{n_2}$ is $C^l$
with all derivatives continuous in $N_1{\times}N_2$, then the map
$\kappa:Q{\times}N_1{\times}\{n_2\}{\rightarrow}N_1{\times}\{n_2\}$
is $C^l$ for each $n_2$ in $N_2$ and depends continuously on $n_2$
in the $C^l$ topology. Furthermore if $\kappa$ and $\kappa'$ are
two such actions which are $C^0$ close, such that $\kappa(q)$ and
$\kappa'(q)$ are $C^k$ close as maps of $N_1$ for all $g{\in}Q$,
then the actions $\kappa$ and $\kappa'$ are $C^k$ close as actions
on $N_1{\times}\{n_2\}$ for any $n_2{\in}N_2$.
\end{theorem}

\begin{proof}
All statements follow from the proof of the results of Bochner and
Montgomery given in \cite[Chapter V, sections 1 and 2.]{MZ}. The
possibility of adding the $N_2$ factor along which the action is
trivial is already noted in \cite[Proof of Lemma $5.12$]{MQ}. That
the actions are actually $C^k$ close follows from the explicit
formulas for derivatives of $\kappa$ along $Q$ given in
\cite[V.2.1]{MZ}.
\end{proof}

\begin{proof}[Proof of Theorem \ref{theorem:smoothalongleaves} from Lemma \ref{lemma:groupsonleaves}]
Possibly after shrinking $U_0$, we can assume that $U_0$ is a
product of a leaf $\tsv_{[\chi]}$ and a small transverse
neighborhood $V_0$.  (We will in fact construct the $\Delta$
action on such a neighborhood.)  The hypotheses $S1-S4$ imply that
the map $\varphi$ intertwines two actions $d$ and $d$ of the group
$\Delta$ such that orbits of $d(\Delta)$ (resp. orbits of
$d'(\Delta)$) are leaves of $\sw_{[\chi]}$ (resp.
$\tsw_{[\chi]}$.) and such that, for each $\delta{\in}\Delta$ and
$v{\in}\varphi(V_0)$, the map
$\delta'(d):\tsw_{[\chi]}{\times}\{v\}{\rightarrow}\tsw_{[\chi]}{\times}\{v\}$
is $C^k$ with all derivatives depending continuously on $v$.
Combined with Theorem \ref{theorem:bochnermontgomery}, this
implies $d':{\Delta}{\times}
\tsw_{[\chi]}{\times}\{v\}{\rightarrow}\tsw_{[\chi]}{\times}\{v\}$
is a $C^l$ action depending continuously on $v{\in}\varphi(V_0)$,
 which suffices
to prove Theorem \ref{theorem:smoothalongleaves}$(1)$. Since $S3$
and Theorem \ref{theorem:bochnermontgomery} also imply that $d$
and $d'$ are close as actions which are $C^k$ along orbits with
all derivatives transversely continuous and that $d'$ tends to $d$
in the natural topology on such actions as
$\rho'{\rightarrow}\rho$, Theorem
\ref{theorem:smoothalongleaves}$(2)$ follows as well.
\end{proof}

The remainder of this section is devoted to a proof of Lemma
\ref{lemma:groupsonleaves}. We begin by constructing the group
$\Delta$ and it's actions on $U_0$ and ${U_0}'$. Recall that we
have identified the tangent bundle to $\base$ with
$H/{\Lambda}{\times}{\fh}{\times}TM$.  We note that $E_{[\chi]}$
is a nilpotent Lie subalgebra of $\fh$. Let $F_{[\chi]}$ denote
the corresponding Lie subgroup of $H$. The smooth foliation of
$\cover$ by left cosets for $F_{[\chi]}$ is tangent to
${\cover}{\times}E_{[\chi]}$ and the projection of this foliation
to $\base$ is exactly ${\sw}_{[\chi]}$.

Fix $x{\in}H{\times}M$ and an orthogonal complement
$E^{\perp}_{[\chi]}$ to $E_{[\chi]}$ in $TX_x$. Note that we can
pick $E^{\bot}_{[\chi]}$ to be a direct sum of a subspace of
$E^{\bot}_{\fh}{\subset}\fh$ and $TM$.  Let $O$ be a small open
disc in $E^{\bot}_{\fh}$. Let $U_1=exp_x(O){\subset}H$. For any
point $x{\in}{\cover}$ we let $x=(x_1,x_2)$ be coordinates for the
product structure and chose a small open neighborhood $U_2$ of
$x_2$ in $M$. Then $U_0=F_{[\chi]}U_1{\times}U_2$ is an open
subset in $\cover$ containing $x$.  If $O$ is small enough then
each leaf of ${\sw}_{[\chi]}$ contained in $U_0$ has a unique
expression as $(F_{[\chi]}ux_1,x_2')$ where $u{\in}U_1$ and
$x_2'{\in}U_2$. We then let $\Delta=F_{[\chi]}$ and let $\Delta$
act on $U_0$ via
$(\delta,f_0ux){\rightarrow}(f_0{\delta}{\inv}ux).$ This clearly
defines a $\Delta$ action on $U_0$ which is $C^{\infty}$, free and
transitive along the leaves of ${\sw}_{[\chi]}$ and establishes
$S1$.

We need to understand the derivative of the $\Delta$ action.  Note
that $T(\base)|_{U_0}=F_{[\chi]}U_1{\times}{\fh}{\times}
TM|_{U_2}$. Also note that the $\Delta$ action is trivial on the
second factor. The fact that we identify the tangent space to $H$
with right invariant vectors and that $d(\delta)$ acts on the
right on $F_{[\chi]}$ orbits, implies the following:

\begin{lemma}
\label{lemma:trivialderivative} For all $f_0{\in}F_{[\chi]},
u{\in}U_1, v{\in}{\fh}, m{\in}M$, and $w{\in}TM_m$ we have
$Dd({\delta})(f_0u,v,m,w)=(f_0{\delta}{\inv}u,v,m,w)$
\end{lemma}

The following lemma records the fact that leaves of $\sw_{[\chi]}$
are injectively immersed in $\base$ and remain injectively when
lifted to $\cover$ and projected to $\quotient$.

\begin{lemma}
\label{lemma:injectivity} The projection of $F_{[\chi]}x$ to
$\base$, $H/{\Lambda}$ or $\quotient$ is an injective immersion.
\end{lemma}

\begin{proof}
This is immediate from the fact that $F_{[\chi]}<F^s_{\rho(a_i)}$
and that the leaves of $\sw^s_{\rho(a_i)}$ are injectively
immersed in $\base$, $H/{\Lambda}$ or $\quotient$ by Proposition
\ref{proposition:explicitfoliations}.
\end{proof}

We choose a lift $\tilde \varphi$ of $\varphi$ to a map from
$H{\times}M$ to $H{\times}M$ and let $U_0'=\tilde \varphi(U_0)$.
Let $\Delta$ act on ${U_0}'$ by letting
$d'(\delta)(x)=\phi(d(\delta)(\phi{\inv}(x)))$ for every
$x{\in}{U_0}'$ and every $\delta{\in}\Delta$.  The properties $S2$
and $S4$ are immediate from this definition.

We now show how to realize the action $\Delta$ differently, in a
way that will allow us to use Theorem
\ref{theorem:centralizernormalform} to prove $S3$.

We first explain why it suffices to consider the case of $\G$
actions. In the case when $\rho$ is a $J$ action, we fix a lattice
$\G<J$. As a consequence of Theorem \ref{theorem:prasadrapinchuk}
and Zariski density of $\G$ in $J$ we can choose the elements
$\Psi$ from Lemma \ref{lemma:tangentdecomposition} and the
elements $a_1,{\ldots},a_q$ in Lemma
\ref{lemma:weightdistributions} to be in $\G$ even when $\rho$ is
a $J$ action.  For the remainder of this section, we can therefore
restrict our attention to the case of $G=\G$.

Let $T$ be a torus as given by Theorem
\ref{theorem:prasadrapinchuk} and restrict $\rho$ and $\rho'$ to
$\Sigma$ actions.  We then form the induced $T$ actions
$\rho_{in}$ and ${\rho_{in}}'$ on $(T{\times}\base)/{\Sigma}$. The
map $\varphi$ can be extended to an equivariant map which we
denote $\varphi_{in}$ from $X_{in}=((T{\times}base)/\Sigma,
\rho_{in})$ to ${X_{in}}'=((T{\times}\base)/{\Sigma},
{\rho_{in}}')$. Let $\ft$ be the Lie algebra of $T$ and identify
$T((T{\times}\base)/\Sigma)$ with
$((T{\times}H{\times}TM)/{\Lambda})/\Sigma{\times}\fh{\times}\ft$.
There is a $\rho_{in}(T)$ invariant smooth foliation of $X_{in}$
given by $\svc[t,x]=[t,\swc(x)]$ tangent  to a $\rho_{in}(T)$
invariant distribution $V_{[\chi]}$ which is $E_{[\chi]}$ viewed
as a subbundle of $T(T{\times}\base)/\Sigma)$. Furthermore there
are analogously defined foliations and distributions $\tsvc$ and
$\tilde V$ on ${X_{in}}'$ and $\psi_{in}$ maps every leaf
$\svc[t,x]$ to the leaf $\tsvc({\varphi_{in}}(t,x))$.

As in \cite{KS} and \cite{MQ}, we will verify $S3$ by verifying it
for these induced actions. It is clear that the actions of
$\Delta$ on $U_0$ and ${U_0}'$ defined above can be extended to
neighborhoods in the universal cover of $X_{in}$ and ${X_{in}}'$
simply be taking the trivial action on the first factor. It is
also clear that individual leaves of $\svc$ and $\tsvc$ are still
injectively immersed in $X_{in}$ and ${X_{in}}'$.

\begin{lemma}
\label{lemma:s3byergodicity} Let $c_t$ be any non-trivial one
parameter $\mathbb R$-split subgroup in $S$ that is in the kernel
of $\chi$ and which has noncompact image when projected to any
simple quotient of $J$. To prove $S3$ it is sufficient to prove
that any leaf of the foliation $\svc$ is contained in the support
of an ergodic component of $c_t$ acting on $X_{in}$.
\end{lemma}

\begin{proof}
We proceed by giving a different description of the group $\Delta$
acting on $X_{in}$. Since $c_t$ is in the kernel of $\chi$, it
follows that the maps $\rho(c_t):\svc(x){\rightarrow}\svc(c_t(x))$
are isometries with respect to the metric on the leaves. Since
every ergodic component consists of entire leaves of $\svc$, for
any point $y{\in}\svc(x)$ there exists a sequence $t_i$ such that
$\lim_{i{\rightarrow}{\infty}}\rho(c_{t_i})x=y$. By passing to a
subsequence, we may assume that
$\rho(c_t):\svc(x){\rightarrow}\svc(\rho(c_{t_i})x)$ converges to
an isometry $\tilde \delta:\svc(x){\rightarrow}{\svc(x)}$ which
takes $x$ to $y$.  The group $\Delta_x$ generated by such limits
is clearly transitive on $\svc(x)$.  Note also that the tangent
map $D\rho(c_{t_i})(x,v)=(\rho(c_{t_i})x,v)$ and so the tangent
map $D\delta(x,v)=({\delta}x, v)$. Since $\Delta$ and $\Delta_x$
both act by isometries with trivial derivative on $\svc(x)$ it is
clear that they are equal.  (For further discussion of this
construction see \cite{MQ} or \cite{KS}.)

We note that, by equivariance,

$$d'(\delta)x=\lim_{i{\rightarrow}\infty}\varphi_{in}{\circ}\rho_{in}(c_{t_i}){\circ}\varphi_{in}{\inv}(x)$$
$$=\lim_{i{\rightarrow}\infty}\rho'_{in}(c_{t_i})x.$$

Therefore,
$d'(\delta)=\lim_{{t_i}{\rightarrow}{\infty}}{\rho_{in}'}(c_{t_i})$.
We let $s{\in}T$ be an element with $\chi(s)<0$, so the natural
extension of $\rho'(s)$ to the tangent bundle of $\tsvc$ is a
contraction with narrow band spectrum which is close to it's
linear part. This follows because this contraction is close to the
one defined by $\rho(s)$ which is linear.  By Theorem
\ref{theorem:normalform} as long as $k{\geq}k_1$ there is a number
$d$ depending only on the dynamics of $\rho(t_i)$, and continuous,
$C^k$ small along fibers, conjugacy $H$ between $\rho'(s)$ and a
polynomial of order less than $d$. Furthermore, since $\rho'(c_t)$
commutes with $\rho'(s)$, by Theorem
\ref{theorem:centralizernormalform} the conjugacy $H$ conjugates
each $\rho'(c_t)$ to a polynomial map of order less than $d$.
Since
$d'(\delta)=\lim_{{t_i}{\rightarrow}{\infty}}{\rho_{in}'}(c_{t_i})$,
it follows that in the coordinates along a fiber given by $H$, the
map $d'(\delta)$ is given by a polynomial of order less than $d$.
Identifying leaves with fibers we see that $d'(\delta)$ is $C^k$
along each leaf, that the $k$-jet of the $d'(\delta)$ depends
continuously on the leaf,  and that the $k$-jet is close to the
one for $d(\delta)$, since it is given by composing a map close to
$d(\delta)$ with a change in coordinates which is $C^k$ small.

We also see that if $\rho$ is $C^l$ for some $l{\geq}k$, then
$d'(\delta)$ is $C^l$ since it is a uniform limit of maps $C^l$
conjugate to polynomials.
\end{proof}

\subsection{\bf Ergodic components and dynamical foliations}
\label{subsection:ergodiccomponents}

We retain all notation from the previous subsection. To prove
Theorem \ref{theorem:smoothalongleaves} it now suffices to prove
the following:

\begin{proposition}
\label{proposition:ergodiccomponents} Let $T$ be the torus
described in the last section and $r$ in $T$ a regular element for
the representation $\sigma$.  For any one parameter subgroup $c_t$
of $S$ which is in the kernel of $\chi$ and which projects to a
noncompact subgroup of each simple factor of $J$, the ergodic
components of $\rho_{in}(g_t)$ consist of entire leaves of
$\sw^s_{\rho_{in}(r)}$ and therefore of entire leaves of
$\sv_{[\chi]}$.
\end{proposition}

We first note an alternate description of $\rho_{in}$. Throughout
this subsection, we assume that we have passed to a finite index
torsion free subgroup of $\G$.  We will also need to pass to
further finite index subgroups of $\G$, but will abuse notation by
retaining the notation $\G$ for each of these successive
subgroups. We recall some facts from \cite{FM1}. First by
\cite[Theorem 6.5]{FM1},  the homomorphism
$\pi:\G{\rightarrow}\Aut(H){\ltimes}H$ defining the action $\rho$
on a subgroup of finite index is a product of two homomorphisms
$\pi_A:\G{\rightarrow}\Aut(H)$ and $\pi_H:\G{\rightarrow}H$ whose
images commute. It follows from the proof of \cite[Theorem
6.5]{FM1} that after changing the algebraic structure on $H$ as in
Proposition \ref{proposition:Hreplacement} and passing to a
further subgroup of finite index, that $\pi(\G)$ is actually
contained in $\Aut(U)$ where $U$ is the unipotent radical of $H$.
Fixing a Levi complement $L$ for $U$ in $H$ and letting $L=ZM$
where $Z$ is a central torus and $M$ is semisimple, the
superrigidity theorems imply that (again after passing to a
subgroup of finite index) $\pi_A(\G)<M$.  After passing to another
finite index subgroup, the restriction of $\pi_A$ to
$\Sigma=\G{\cap}T$ extends to a homomorphism
$\pi_A^T:T{\rightarrow}\Aut(U)$ and the restriction of $\pi_H$ to
$\Sigma$ extends to a homomorphism $\pi_H^T:T{\rightarrow}M$.
These homomorphisms are not quite canonical, but suffice for our
purposes.  It is clear that the images of these homomorphisms
commute and so we can define a homomorphism
$\pi^T(t)=\pi^A_T\pi^H_T(t)$.

We can now give a simple description of a finite cover of the
induced action. In the description in the last paragraph, we
passed to a finite index subgroup of $\G$, which also causes us to
pass to a finite index subgroup $\Sigma'<\Sigma$.  The map
$(t,h){\rightarrow}(t,\pi^T(t)h)$ descends to a map from
$(T{\times}H/{\Lambda})/{\Sigma'}$ to
$(T{\ltimes}H)/(\Sigma'{\ltimes}\Lambda)$ where the semidirect
product is defined by $\pi^A_T$ and $\pi^A_T(\Sigma)$ normalizes
$\Lambda$ by definition.  This map conjugates the induced action
to an action defined by $\rho_T(t_0)[t,h]=[t_0t, \pi_H^T(t)h]$. We
summarize this discussion with the following Proposition.

\begin{proposition}
\label{proposition:describinginducedactions} \begin{enumerate}
\item If $\rho$ is an affine action, then there is a finite index
subgroup $\Sigma'$ such that the lift of the action $\rho_{in}$ to
$(T{\times}H/{\Lambda})/\Sigma'$ is smoothly conjugate to a left
translation action $\rho_T$ of $T$ on
$(T{\ltimes}H)/(\Sigma'{\ltimes}\Lambda)$ as described above.
 \item
 If $\rho$ is a
quasi-affine action, then there is a finite index subgroup
$\Sigma'$ in $\Sigma$ and a left translation action $\rho_T$ of
$T$ on $(T{\ltimes}H)/(\Sigma'{\ltimes}\Lambda)$ as above and a
cocycle
$\iota:T{\times}(T{\ltimes}H)/(\Sigma{\ltimes}\Lambda){\rightarrow}\Isom(M)$,
such that the lift of $\rho_{in}$ to
$(T{\times}H/{\Lambda}{\times}M)/\Sigma'$ is smoothly conjugate to
the skew product action over $\rho_T$ defined by $\iota$.
\end{enumerate}
\end{proposition}

We begin by showing that, even for quasi-affine actions, it
suffices to consider the action on
$(T{\ltimes}H)/(\Sigma'{\ltimes}\Lambda)$.

\begin{lemma}
\label{lemma:zimmerthesis} For any one parameter subgroup $g_t$ of
$T$, and any regular element $r$ in $T$, if the ergodic components
of the left translation action $g_t$ on
$(T{\ltimes}H)/(\Sigma'{\ltimes}\Lambda)$ consist of entire leaves
of $\sw^s_{\rho_{T}(r)}$ then the ergodic components of
$\rho_{in}(g_t)$ consist of entire leaves of
$\sw^s_{\rho_{T}(r)}$.
\end{lemma}

\begin{proof}
This follows two facts.  The first is one of the main results of
Zimmer's thesis. This says that if $K$ is a compact group and $K$
acts on a standard probability measure space $(Y,\nu)$, and $\rho$
is an action of locally compact group $G$ by measure preserving
transformations on a standard measure space $(X,\mu)$ and
$\iota:G{\times}X{\rightarrow}K$ is a cocycle, then the ergodic
components of the skew-product action of $G$ on
$(X{\times}Y,\mu{\times}\nu)$ are of the from $E{\times}L{\cdot}y$
where $E$ is an ergodic component of $X$, $y$ is a point in $Y$
and $L$ is a subgroup of $K$ such that $\iota$, restricted to $E$,
is cohomologous to a cocycle taking values in $L$.

The second fact describes dynamical foliations for skew product
extensions.  Again, let $K$ be a compact group.  Let $X$ be a
smooth compact manifold and $Y$ be an associated bundle to a
principal $K$ bundle over $X$.  Assume $G$ is a locally compact
group and that $\rho:G{\times}X{\rightarrow}X$ and $\tilde
\rho:G{\times}Y{\rightarrow}Y$ be two actions which commute with
the bundle projection $\pi:Y{\rightarrow}X$.  Then the $G$ action
on $Y$ is measurably isomorphic to a skew product extension as
described in the previous paragraph and for any $g{\in}G$ which is
partially hyperbolic and normally hyperbolic to a central
foliation on both $X$ and $Y$, the map $\pi$ is a diffeomorphism
from each leaf of $\sw^s_{\tilde \rho(g)}$ onto a leaf of
$\sw^s_{\rho(g)}$. This follows from the dynamical
characterization of $\sw^s_{\rho(g)}$ in \cite[Theorem 6.8e]{HPS}

To prove the lemma, we apply these facts twice, first to the $g_t$
action on $(T{\times}\base)/{\Sigma'}$ covering the
$\rho_{\in}(g_t)$ action on $(T{\times}\base)/{\Sigma}$ and second
to the actions $\rho_{in}$ on
$(T{\ltimes}H)/(\Sigma'{\ltimes}\Lambda){\times}M$ which is a skew
product action over the action $\bar \rho_{in}$ on
$(T{\ltimes}H)/(\Sigma'{\ltimes}\Lambda)$.
\end{proof}

We are now reduced to identifying ergodic components for left
translation actions on homogeneous spaces.  To do this, we will
use work of Brezin and Moore \cite{BM}.  Following that paper, we
note that for any Lie group $L$, any finite volume homogeneous
space $L/\Delta$ has two special quotients an {\em maximal toral
quotient} and a {\em maximal semisimple quotient}. An affine
quotient of the space $L/{\Delta}$ is one of the from
$P/{\phi(\Delta)}$ where $\phi:L{\rightarrow}P$ is a surjective
homomorphism.  The maximal toral quotient is the maximal affine
quotient of $L/{\Delta}$ which is a torus and the maximal
semisimple quotient is the maximal affine quotient of $L/{\Delta}$
where $P$ is semisimple.  Given a one parameter subgroup $l_t$ in
$L$, we can project $l_t$ to either the torus or $M$ and this
defines a quotient of the left translation flow of $l_t$ on
$L/{\Delta}$. Let $\phi_1(l_t)$ be the quotient action on the
maximal toral quotient and let $\phi_2(l_t)$ be the quotient
action on the maximal semisimple quotient. The following is a
restatement of \cite[Theorem 6.1]{BM}.

\begin{theorem}
\label{theorem:brezinmoore} Let $l_t$ be a one parameter subgroup
of $L$ acting by left translation on a finite volume homogeneous
space $L/{\Delta}$ for a Lie group $L$.  Then the action of $l_t$
is ergodic if and only if both $\phi_1(l_t)$ and $\phi_2(l_t)$ are
ergodic.
\end{theorem}

To prove Proposition \ref{proposition:ergodiccomponents}, we
require an additional lemma which is an immediate consequence of
Theorem \ref{theorem:prasadrapinchuk}$(3)$.

\begin{lemma}
\label{lemma:torusergodic} Let $c_t$ be a one parameter subgroup
of $T$ which is in the kernel of $\chi$ and projects to a
non-compact subgroup in each $J_i$.  Then the action of $c_t$ on
$T/\Sigma'$ is ergodic.
\end{lemma}

\begin{proof} This is immediate since an ergodic component of the
action is necessarily of the form
$\prod_IT_i'/(\Sigma'{\cap}T_i')$ where $T_i'<T_i$ is a subtorus
and $\Sigma'{\cap}T_i'$ is a lattice in $T_i'$.  This forces
$T_i'$ to be the Zariski closure of a subgroup of $\Sigma'$ and
therefore to be algebraic.  Theorem
\ref{theorem:prasadrapinchuk}$(3)$ then implies that $T_i'=T_i$.
\end{proof}

\begin{proof}[Proof of Proposition \ref{proposition:ergodiccomponents}]
By Proposition \ref{proposition:describinginducedactions} and
Lemma \ref{lemma:zimmerthesis} we are reduced to showing that
ergodic components of the $c_t$ action on
$(T{\ltimes}H)/(\Sigma'{\ltimes}\Lambda)$ consist of entire leaves
of  $\sw_{\rho_L(r)}^s$.  We do this by explicitly identifying
ergodic components, or rather explicitly identifying ergodic
components modulo finite extensions.

 Note that arguments as in Lemma
\ref{lemma:zimmerthesis} shows that there is no loss of generality
in passing to finite covers, so for simplicity we pass to a finite
cover of $H$ such that:
\begin{enumerate}
\item $\Lambda$ is torsion free, \item The Levi complement of $L$
is the direct product of $Z$ and $M$, \item $M$ is the direct
product of its simple factors and $\Lambda$ does not intersect the
center of $M$, \item $Z$ is a direct product of copies of $S^1$
and copies of $\Ra^*$.
\end{enumerate}

Let $\sigma_M$ be restriction to $M$ of the map from $L$ to
$\Aut(U)$ defining the semi-direct product structure of $H$ as
$L{\ltimes}U$.  Let $M^U$ be the kernel of $\sigma_M$ and let
$M^U_K$ be the maximal normal connected compact subgroup of $M^U$.
Note that our assumptions imply that $M^U_K$ is normal in $H$ and
that $H=H{\times}M^U_K$ where $H'$ is isomorphic to $H/M^U_K$.

We write $M'$ as $M''{\times}C$ where $C$ is the maximal connected
normal compact subgroup of $M'$. It from the proof of
\cite[Theorems $1$ and $2$]{A} that in $H$,
\begin{enumerate}
\item $\Lambda{\cap}M^U_K{\times}U=\Lambda_U$ is a lattice in
$M^U_K{\times}U$ and projects to a lattice in $U$, \item the
projection of $\Lambda$ to $M''$ is a lattice in $M''$.
\end{enumerate}

 Since $\pi_{H}$ is a homomorphism of $\G$ the superrigidity
 theorems imply that there is a homomorphism
 $\pi_{H}^E:J{\rightarrow}H$ and a homomorphism
$\pi_{H}^K:\G{\rightarrow}H$  with bounded image such that the
images commute and $\pi_{H}(\g)=\pi_{H}^E(\g)\pi_H^K(\g)$. Note
that, after passing to a further finite index subgroup, $\pi_{H}$
necessarily takes values in $M'$ and $\pi_{H}^E$ necessarily takes
values in $M''$. Using that $\Aut(U)$ is an algebraic group, we
can also write $\pi_A$ as a product of
$\pi_A^E:J{\rightarrow}\Aut(U)$ and
$\pi_A^K:\G{\rightarrow}\Aut(U)$.  We write $M'$ as a direct
product $M_1M_2$ where $M_1$ is the minimal product of simple
factors of $M'$ such that $\pi_H^E$ takes values in $M_1$ and the
projection of $\Lambda$ to $M_1$ is a lattice.  This implies that
$M_1$ is a direct product of semisimple groups $M_1^i$ where the
projection of $\Lambda$ to $M_1$ is commensurable to a product of
irreducible lattices $\Lambda_{M_1}^i<M_1^i$ and such that the
projection of $\pi^T_H(c_t)$ is non-compact in each $M_1^i$.  This
implies that the left translation action of $c_t$ on
$M_1/{\Lambda_{M_1}}$ defined by $\pi^T_H$ is ergodic.  Since
$\pi_{H'}^E$ has non-trivial image in each $M_1^i$, it follows
that each $M_1^i$ has real rank at least $2$ and so there is a
compact connected normal subgroup $M_1^K<M'$ such that
$(M_1{\times}M_1^K){\cap}{\Lambda}=\Lambda_{M_1}$ at least after
replacing $\Lambda$ by a subgroup of finite index.  We write
$M_1'$ for $M_1{\times}M_1^k$.  It is also easy to see that the
product $T{\cdot}M_1$ is a subgroup of $L$ and that
$T{\cdot}M_1{\cap}(\Sigma'{\ltimes}\Lambda)=\Sigma'{\times}\Lambda_{M_1}$.

We now construct a subgroup of $U$.  The fact that $\Lambda_U$
projects to a lattice in $U$ defines a rational structure on $U$
and $\fu$. We let let $\mu_1$ be the composition of $\pi_H^E$ with
the restriction of $\Ad_{T{\times}H}$ to $\fu$ and let $\fu_H$ be
the minimal Lie subalgebra containing all non-trivial root
subspaces of $\mu_1$ and invariant under $\pi^T_H(T)$. We let
$\mu_2$ be the composition of $\pi^E_A$ with the representation of
$\Aut(U)$ on $\fu$ and let $\fu_A$ be the minimal Lie algebra
containing all non-trivial root subspaces of $\mu_2$ and invariant
under $\pi_A^T$.  Finally we let $\fu_0$ be the minimal rational
Lie subalgebra of $\fu$ containing both $\fu_A$ and $\fu_H$ and
invariant under $T{\cdot}M_1$. Let $U_0<U$ be the Lie subgroup
with Lie algebra $\fu_0$. We let $K_0$ be the closure of the
projection of $\Lambda_U$ to $M^U_K$. Form the semidirect product
$N=(T{\times}K_0{\times}M'_1){\ltimes}U_0$. By construction it is
clear that:
\begin{enumerate}
\item $N{\cap}{\Lambda}=\Lambda_N$ is a lattice in $N$ \item For
any regular elements $r$ in $T$, the space $E_{\rho^{in}}^s(r)$ is
a subspace of $\fn$ \item the maximal semisimple quotient of
$N/{\Lambda_N}$ is $M_1/{\Lambda_{M_1}}$ and \item the maximal
toral quotient of $N/{\Lambda_N}$ is $T/{\Sigma'}$.
\end{enumerate}
Together with Lemma \ref{lemma:torusergodic} and the definition of
$M_1$ this implies that any ergodic component of the action of
$c_t$ on $(T{\ltimes}H)/(\Sigma'{\ltimes}\Lambda)$ contains a
translate of $N/{\Lambda_N}$ in
$(T{\ltimes}H)/(\Sigma'{\ltimes}\Lambda')$ which suffices to prove
the proposition.

\end{proof}

\section{\bf Final arguments}
\label{section:finalarguments}

\subsection{Elliptic operators and global regularity}
\label{subsection:isdiffeo}

In this section, we prove that $\varphi$ is a diffeomorphism. Here
$\rho,G$ and $\base$ are as in the remarks at the end of section
\ref{section:conjugacy}. The number $k_0$ is the smallest number
that allows us to apply the techniques of subsection
\ref{subsection:normalforms} to show that $\varphi$ is a
diffeomorphism along foliations of the type $\sw^g_{[\chi]}$. In
keeping with the statement of Theorem \ref{theorem:main} we let
$n=\frac{\dim(\base)}{2}+3$.  We now prove:

\begin{theorem}
\label{theorem:diffeomorphism} There is a neighborhood $V$ of
$\rho$ in $\Hom(G,\Diff^k(X))$ such that if $\rho'{\in}V$, the map
$\varphi$ constructed above is a $C^{k-n}$ small $C^{k-n}$
diffeomorphism which is conjugacy between $\rho$ and $\rho'$.
Furthermore
\begin{enumerate}
\item $\varphi{\rightarrow}\Id$ as $\rho'{\rightarrow}\rho$ and,
\item given $l{\geq}k$, we can choose $V$ so that if $\rho'$ is
$C^{\infty}$ and $\rho'{\in}V$, then, the map $\varphi$ is $C^l$.
\end{enumerate}
\end{theorem}

\noindent{\bf Remark:}  The proof below uses only standard facts
concerning elliptic operators and is straightforward.  The result
stated also follows from the main theorem in \cite{KS2}, but as
that article relies on much deeper and harder results concerning
hypo-elliptic operators, we give the proof below.

\begin{proof} We choose a finite cover of $\base$ by open sets $U_m$ such that:

\begin{enumerate}
\item each $U_m$ is contained in a neighborhood $W_m$ which is
coordinate chart on $\base$

\item for each $U_m$ we have $\varphi(U_m){\subset}W_m$

\item each $U_m$ is a foliation chart for $\ff$ and
$\sw_{[\chi]}^g$ for all $g{\in}\Phi$ and $\chi$ in $\bar \Sigma$
as defined in subsection \ref{subsection:weighffoliations},

\item each $U_m$ is of the form $U_{1m}{\times}U_{2m}$ where
$U_{1m}$ is an open set in $H/{\Lambda}$ and $U_{2m}$ is an open
set in $M$.
\end{enumerate}

For convenience, we denote $\ff$ by $\ff_0$ and fix an order on
the $\sw_{[\chi]}^g$ and relabel them $\ff_1,{\ldots},\ff_q$. We
choose a basis $X_{ij}$ of $\fh$ where $X_{0j}$ for
$1{\leq}j{\leq}\dim(Z)$ is a basis for $\fz$ and $X_{ij}$ is a
basis for $\ff_i$ with $1{\leq}i{\leq}q$. For $U_{2j}$ we choose
an explicit identification with an open ball in $\Ra^n$ and choose
a basis of constant vector fields $X_{0j}$ where
$\dim(Z)+1{\leq}j{\leq}\dim(Z)+\dim(M)$.

Identifying each $W_m$ with a subset of $\Ra^n$, we can write
$\varphi_m=\varphi|_{U_m}=\Id+h_m$ where
$h_m:U_m{\rightarrow}\Ra^n$ is $C^0$ small.

For any $\eta>0$ and any $l'>k$, by choosing $V$ small enough, and
applying Theorem \ref{theorem:smoothalongleaves} and Theorem
\ref{theorem:candidateconjugacy}, we have that each $h_m$ is
\begin{enumerate}
\item $C^k$ along each leaf of $\ff_l$ for $1{\leq}l{\leq}q$ with
$\supp_{U_m}X_{lj}^n(h_m)<\eta$ for any $0{\leq}n{\leq}k$ and any
$1{\leq}j{\leq}\dim(\ff_l)$,

\item $C^{k-2}$ along almost every leaf of $\ff$ with
$\int_{U_m}\|X_{0j}^{n}(h_m)(x)\|^2d\mu<\eta$ for any
$0{\leq}n{\leq}k-2$ and any $0\leq{j}\leq\dim(Z)+\dim(M)$,

\item if $\rho'$ is $C^{\infty}$ then $h_m$ is $C^{l'}$ along
$\ff_0$ and $C^{\infty}$ along each $\ff_i$ for $1{\leq}i{\leq}q$.
\end{enumerate}

We construct an elliptic operator as follows. Let $c$ be the least
even integer less than or equal to $k-2$, then the operator
$$\Delta=\sum_{i=1}^q\sum_{j=1}^{\dim(\ff_i)}X_{ij}^c$$

\noindent is elliptic with smooth coefficients on each $U_l$.
Standard estimates, see e.g. \cite[Section 6.3]{Z3}, imply that:

$$\|u\|_{2,k-2}<C(\|\Delta(u)\|_2+\|u\|_2)$$

\noindent for any $u$ in the Sobolev space $W^{2,k-2}(U_l)$, where
$W^{2,k-2}(U_l)$ is the Sobolev space of functions with $k$ weak
derivatives in $L^2$,  $\|{\cdot}\|_{2,k-2}$ is the Sobolev norm
and $\|{\cdot}\|_2$ is the $L^2$ norm.  (For $k-2$ odd, the
standard inequality involves the $\|{\cdot}\|_{2,1}$ for both
terms on the right hand side, but we will not need this.)

We want to apply this estimate to $h_m$, but $h_m$ is not a priori
in $W^{2,k-2}(U_m)$.  We let $U_m^{\varepsilon}$ be the set of
points $x$ in $U_m$ such that $B(x,\varepsilon){\subset}U_m$. To
complete the argument, we use mollifiers $J_{\varepsilon}$ such
that
\begin{enumerate}
\item $X_{ij}{J_{\varepsilon}}=J_{\varepsilon}X_{ij}$ for $X_{ij}$
above and,

 \item $J_{\varepsilon}u$ is defined on $U_m^{\varepsilon}$ and $J_{\varepsilon}$ maps $L^{1,loc}(U_m)$
to $C^{\infty}(U^{\varepsilon}_m)$

\item $J_\varepsilon$ are uniformly bounded on
$W^{2,k-2}(U_m){\subset}L^1_{loc}(U_m)$.

\item $J_{\varepsilon}$ converges uniformly to the identity on
$L^{1,loc}$ as $\varepsilon{\rightarrow}0$.

\end{enumerate}

\noindent  We briefly describe the operators $J_{\varepsilon}$
which are convolution operators for a family of functions
$f_{\varepsilon}$.  We write the function
$f_{\varepsilon}=f_{1\varepsilon}f_{2\varepsilon}$ where
$f_{i\varepsilon}$ is a function on $U_{im}$ for $i=1,2$. The
function $f_{1\varepsilon}$ is a standard mollifier and define
$J_{1\varepsilon}$ by standard convolution. We define $f_2$ by
taking a standard mollifier on a small neighborhood of zero in
$\fh$ and pulling back to $H$ via the inverse of the exponential
map.  We identify $U_{1m}$ with a small neighborhood in $H$ and
define
$J_{2\varepsilon}(u)=\int_{H}f_{2\varepsilon}(h)u(xh,m)d\mu$. The
fact that we act on the right on $x$ in the formula is necessary
to guarantee condition $(1)$ above.   We then let
$J_{\varepsilon}=J_{1\varepsilon}J_{2\varepsilon}$. It is easy to
see that $J_{1\varepsilon}$ and $J_{2\varepsilon}$ commute, and
that $J_\varepsilon$ satisfies $1-4$ above.

Letting $\varepsilon_n=\frac{1}{n}$, we have that
$$\|(J_{\varepsilon_n}-J_{\varepsilon_{n+1}})h_m\|_{2,k-2}<C(\|(J_{\varepsilon_n}-J_{\varepsilon_{n+1}})\Delta(h_m)\|_2
+\|(J_{\varepsilon_n}-J_{\varepsilon_{n+1}})h_m\|_2).$$
 The right
hand side converges  to zero, which implies that
$\{J_{\varepsilon_n}h_m\}_n$ is a Cauchy sequence in
$W^{2,k-2}(U_i)$. Since $\{J_{\varepsilon_n}h_m\}_n$ converges in
$L^{1,loc}$ to $h_m$, this implies that $h_m{\in}W^{2,k-2}$ and
$$\|h_m\|_{2,k-2}<C(\|\Delta(h_m)\|_2
+\|h_m\|_2)$$ which by the properties $(1)$ and $(2)$ of $h_m$
described above, implies that $\|h_m\|_{2,k-2}<C\eta$ for a
constant $C$ not depending on $\rho'$.  By the Sobolev embedding
theorems, this implies that $h_m$ is $C^{k-n}$ small where $n$ is
$\frac{\dim(\base)}{2}$.

This then implies that $\varphi$ is $C^{k-n}$ close to the
identity, which implies that $\varphi$ is a diffeomorphism, since
there is a neighborhood of the identity in the space of $C^{k-n}$
maps which consists of diffeomorphisms.

To show that $h_m$ is $C^l$ when $\rho'$ is $C^{\infty}$ follows a
similar outline.  We choose $V$ such that $h_m$ is $C^{\infty}$
along each $\ff_j$ for $1{\leq}j{\leq}q$ and $C^{l'}$ along
$\ff_0$ where $l'{\geq}l+\frac{\dim(X)}{2}+3$ is even.  The same
argument with $c=l'$ in the construction of the elliptic operator
shows that $h_m$ is in $W^{2,l'}$ and therefore is $C^l$. Note
that since we do not have a good bound on the $W^{2,l'}$ norm of
$\varphi$ along $\ff$ in Theorem \ref{theorem:candidateconjugacy}
or of the $C^{l'}$ norm of $\varphi$ along $\sw^g_{[\chi]}$ in
Theorem \ref{theorem:smoothalongleaves}, we do not obtain a bound
on the $C^l$ size of $\varphi$.
\end{proof}

\subsection{Smooth perturbations, smooth conjugacy, and
iterations} \label{subsection:smoothcase} We keep all the notation
from the previous subsection.

 For notational convenience
in the proof of the $C^{\infty}$ case of Theorem
\ref{theorem:main}, it is convenient to fix right invariant
metrics $d_l$ on the connected components of $\Diff^l(X)$ with the
additional property that if $\varphi$ is in the connected
component of $\Diff^{\infty}(X)$, then $d_l(\varphi,
\Id){\leq}d_{l+1}(\varphi,\Id)$.  To fix $d_l$, it suffices to
define inner products $<,>_l$ on $\Vect^l(X)$ which satisfy
$<V,V>_l{\leq}<V,V>_{l+1}$ for $V{\in}\Vect^{\infty}(X)$. As
remarked in \cite[Section 6]{FM2}, after fixing a Riemannian
metric $g$ on $X$, it is straightforward to introduce such metrics
using the methods of \cite[Section 4]{FM2}.

Once we have fixed the family of metrics $d_l$ and fix a
generating set $K$ for $G$, it is possible to rephrase parts of
Theorem \ref{theorem:diffeomorphism} more quantitatively as
follows:

\begin{corollary}
\label{corollary:diffeomorphismwithestimates} In the setting of
Theorem \ref{theorem:diffeomorphism}, given $k{\geq}k_0$ and
$l{\geq}k$, for every $\varepsilon>0$ there exists $\delta>0$ such
that if $\rho'$ is an action of $G$ on $\base$ with
$d_k(\rho'(g)\rho(g){\inv},\Id)<\delta$ for all $g{\in}K$ then
there exists a $C^l$ conjugacy $\varphi$ between $\rho$ and
$\rho'$ such that $d_{k-n}(\varphi,\Id)<\epsilon$.
\end{corollary}

{\noindent}{\bf Remark:} The algorithm presented below for proving
smoothness is essentially contained in the proof of the
$C^{\infty,\infty}$ case of \cite[Theorem $1.1$]{FM2}.

\begin{proof}[Proof of $C^{\infty,\infty}$ local rigidity in
Theorem \ref{theorem:main}] If $\rho'$ is a $C^{\infty}$
perturbation of $\rho$, then there exists some $k>1$, such that
$\rho'$ is $C^k$ close to $\rho$ and we can assume that
$k{\geq}k_0$. We fix a sequence of positive integers
$k=l_0<l_1<l_2<\cdots<l_i<\ldots$ with $l_{i+1}-l_i>n+3$ for each
$i$.  We construct a sequence of $C^{\infty}$ diffeomorphisms
$\phi_i$ such that the sequence
$\{\phi_n{\circ}{\ldots}{\circ}\phi_1\}_{n{\in}\Na}$ converges in
the $C^{\infty}$ topology to a conjugacy between $\rho$ and
$\rho'$.

We let $\phi^i=\phi_i{\circ}{\ldots}{\circ}\phi_1$ and
$\rho_i=\phi^i{\circ}\rho'{\circ}(\phi^i){\inv}$ and construct
$\phi_i$ inductively such that
\begin{enumerate}
\item $\rho_i$ is sufficiently $C^{l_i-n-3}$ close to $\rho$ to
apply Corollary \ref{corollary:diffeomorphismwithestimates} to
$\rho_i$ and $\rho$ with $l=l_{i+1}$ and
$\epsilon=\frac{1}{2^{i+2}}$,

\item $d_{l_i-n-3}(\phi_i,\Id)<\frac{1}{2^i}$ and,

\item
$d_{l_{i}-n-3}(\rho_i(\g){\circ}\rho(\g){\inv},\Id)<\frac{1}{2^i}$
for every $\g{\in}K$.

\end{enumerate}

To construct $\phi_{i+1}$, we assume that $\rho_i$ is close enough
to $\rho$ in the $C^{l_i}$ topology to apply Corollary
\ref{corollary:diffeomorphismwithestimates} with $l=l_{i+1}$ and
$\varepsilon=\frac{1}{2^{i+1}}$. Then we have a $C^{l_{i+1}-n-3}$
diffeomorphism $\psi_{i+1}$ with
$\psi_{i+1}{\circ}\rho_i{\circ}\psi_{i+1}{\inv}=\rho$ and
$d_{l_{i+1}-n-3}(\psi_{i+1},\Id)<\frac{1}{2^{i+2}}$. Using
standard approximation theorems, we can choose a $C^{\infty}$
diffeomorphism $\phi_{i+1}$ with
$d_{l_{i+1}-n-3}(\phi_{i+1},\Id)<\frac{1}{2^{i+1}}$ and
$\rho_i=\phi_{i+1}{\circ}\rho_{i}{\circ}\phi_{i+1}{\inv}$ close
enough to $\rho$ in the $C^{l_{i+1}-n-3}$ topology to apply
Corollary \ref{corollary:diffeomorphismwithestimates} with
$l=l_{i+2}$ and $\varepsilon=\frac{1}{2^{i+3}}$ and so that $(3)$
above is satisfied.

To start the induction it suffices that $\rho'$ is sufficiently
$C^k$ close to $\rho$ to apply Theorem
\ref{theorem:diffeomorphism}$(2)$ with $l=l_1$ and
$\epsilon=\frac{1}{2}$.

It remains to show that the sequence
$\{\phi_n{\circ}{\ldots}{\circ}\phi_1\}_{n{\in}\Na}$ converges in
the $C^{\infty}$ topology to a conjugacy between $\rho$ and
$\rho'$.  Combining condition $(2)$ with the fact that
$d_{l_i}(\phi_i,\Id){\leq}d_j(\phi_i,\Id)$ for all $j{\geq}l_i$,
and the fact that $d_{l_i}$ is right invariant implies that
$d_{l_i-n-3}(\phi_j,\Id)=d_{l_i-n-3}(\phi^j,\phi^{j-1}){\leq}\frac{1}{2^j}$
for all $j{\geq}i$.  This implies that $\{\phi^j\}$ is a Cauchy
sequence in $\Diff^{l_i-n-3}(X)$ for all $i$, and so that
$\{\phi^j\}$ converges in $\Diff^{\infty}(X)$.  Similarly,
condition $(3)$ implies $\rho_i$ converges to $\rho$ in the
$C^{\infty}$ topology.
\end{proof}

\bigskip
\noindent
David Fisher\\
Department of Mathematics and Computer Science\\
Lehman College - CUNY\\
250 Bedford Park Boulevard W\\
Bronx, NY 10468\\
david.fisher@yale.edu\\

{\noindent}Gregory Margulis\\
Department of Mathematics\\
Yale University\\
P.O. Box 208283\\
New Haven, CT 06520\\
margulis@math.yale.edu\\

\end{document}